\documentclass[12pt,amsfonts]{article}

\usepackage{hyperref}
\usepackage{bm}
\usepackage{calrsfs}
\usepackage{calligra}
\usepackage[T1]{fontenc}
\usepackage{latexsym}
\usepackage{amsbsy}
\usepackage{amsfonts} 
\usepackage{amssymb}
\usepackage{mathrsfs} 

\usepackage{tikz}
\usetikzlibrary{math} 
\usetikzlibrary{calc} 

\usepackage{xcolor}

\DeclareMathAlphabet{\mathcalligra}{T1}{calligra}{m}{n}

\newtheorem{theorem}{Theorem}[section]
\newtheorem{lemma}[theorem]{Lemma}%[section]
\newtheorem{corollary}[theorem]{Corollary}%[section]
\newtheorem{proposition}[theorem]{Proposition}%[section]

\newtheorem{definition}[theorem]{Definition}
\newcommand{\bd}[1]{\begin{definition}\label{#1}\rm}
\newcommand{\ed}{\end{definition}}
\newcommand{\bt}[1]{\begin{theorem}\label{#1}}
\newcommand{\vet}{\end{theorem}}
\newcommand{\bprop}[1]{\begin{proposition}\label{#1}}
\newcommand{\eprop}{\end{proposition}}
\newcommand{\bcor}[1]{\begin{corollary}\label{#1}}
\newcommand{\ecor}{\end{corollary}}

\newcommand{\lra}{\longrightarrow}

\newcommand{\stack}[2]{\raisebox{-2pt} 
{\renewcommand{\arraystretch}{.01} 
\begin{tabular}{c} 
$#2$\\$\scriptscriptstyle #1$ 
\end{tabular} 
}}

\renewcommand{\l}{{\it leb}}

\newcommand{\vp}{\varphi}
\newcommand{\ve}{\varepsilon}
\newcommand{\nid}{\noindent}
\newcommand{\qed}{\hfill$\Box$} 

\def\1{\, {\rm I}\mskip-10mu 1} 

\def\Nab{\, \nabla{\hspace{-3mm}\nabla}}

\def\h{\, {h}{\hspace{-2.85mm}{h}}}
\def\sh{\, {h}{\hspace{-2.05mm}{h}}}
\def\bh{\boldsymbol{h}}
\def\c{\, {c}{\hspace{-2.2mm}{c}}} 
\def\d{\, {d}{\hspace{-2.6mm}{d}}}

\def\+{\, \ensuremath{\boldsymbol{\large\pmb{+}}}\, } 
\def\plus{\begin{array} {c}{{\Large{+}}} \\ \vspace{-9.7mm} \\ {\hspace{1mm}{\large{+}}}\end{array}}
\def\g{\, {\gamma}\mskip-11.8mu \gamma} %\,

\renewcommand{\t}[1]{\tilde{#1}} 
\newcommand{\wt}{\widetilde}

\begin{document}
\title{Elements of the Stochastic Calculus for a Class of Boltzmann type processes  and Application to Regularity of Densities} 
\par
\author{J\"org-Uwe L\"obus
%%\footnotemark[1]
\\ Matematiska institutionen \\ 
Link\"opings universitet, %\\ 
% SE-581 83 Link\"oping, \\ 
Sverige  \\ 
{\tt jorg-uwe.lobus@liu.se} \\ 
{\tt https://orcid.org/0000-0001-5646-1277}
}
\date{}
\maketitle
{\footnotesize
\noindent
\begin{quote}
{\bf Abstract} For a class of piecewise deterministic random processes we introduce a stochastic calculus which is a certain non-Gaussian counterpart to the classical Malliavin calculus. As an application we investigate the regularity of densities of $N$-particle Boltzmann type processes at time $t>0$. 
\noindent 

{\bf AMS subject classification (2020)} primary 76P05, secondary 60H07

\noindent
{\bf Keywords} Piecewise deterministic random processes, $N$-particle Boltzmann type processes, Non-Gaussian stochastic calculus, Regularity of densities

\end{quote}
}

\section{Introduction}\label{sec:1}
\setcounter{equation}{0} 

Boltzmann type processes in the terminology of e.g.\ \cite{RW05}, Section 2, form a certain class of piecewise deterministic random processes. These processes model a collection of particles (molecules), each given by its position, velocity, and mass. In a simplified setting we just consider position as well as velocity and keep the number of particles constant $N$. The particles move freely through the physical state space until they hit the boundary or collide with another particle. In both instances they change velocity including direction according to some specified random mechanism. 

We are interested in the clustering and the  regularity of the $N$-particle density at a fixed time $t>0$. Figures 1-3 below illustrate this objective by schematically following the position density of one fixed particle at time zero and two more small times.  
\bigskip

\vspace*{-.8cm}\hspace*{-3cm}\includegraphics[scale=1.0]{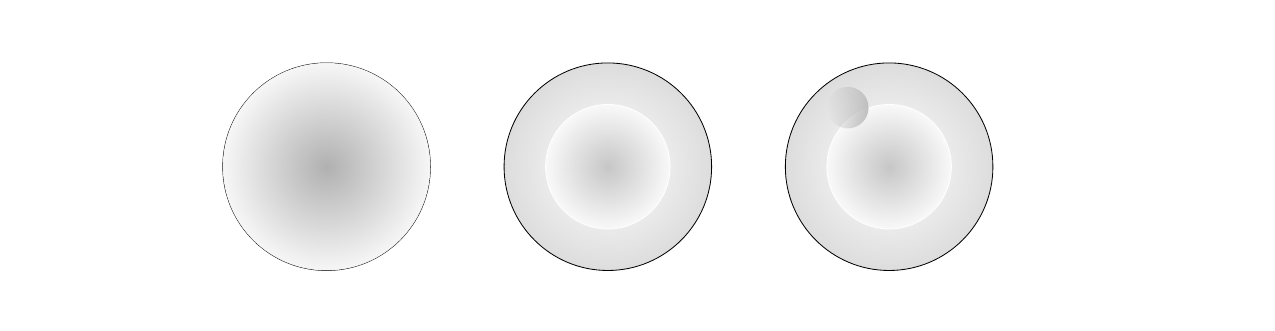} 

\vspace*{-1cm}

\noindent
{\hspace{1cm}\small{{\bf Fig.~1}\hspace{3.8cm}{\bf Fig.~2}\hspace{3.8cm}{\bf Fig.~3}\\ \\ 
Figure 1 displays schematically the distribution of the particle at time $t=0$. We shall assume that this distribution is given by a density which is positive on the whole physical state space. Provided that all $N$ particles are initially located at different places, with positive probability the first change of velocity is caused by a diffusive (i.e. random) reflection at the boundary. Figure 2 shows schematically the distribution of the particle at a sufficiently small time $t_1>0$ before the first collision with one of the other $N-1$ particles occurs. Figure 3 displays schematically the distribution of the fixed particle at a time $t_2>t_1$ shortly after the first collision with one of the other particles took place. In Figure 3 it is assumed that $t_2-t_1$ is sufficiently small. }} 
\bigskip

{\bf Related previous research.} Piecewise deterministic processes are models of (particle) systems that evolve deterministically between consecutive random times of random changes in the system. Such systems have been studied from different  perspectives for several decades. For some recent results on existence of densities in particular situations see e.g. \cite{Lo18} and \cite{CHW21}. However, to the author's best knowledge, the problem of clustering and regularity of the density of certain $N$-particle Boltzmann type processes has never been addressed before. 
\medskip

In addition to research on the regularity of densities of $N$-particle processes, there is intensive research on the regularity of solutions to classes of PDEs, including Boltzmann type equations. For the latter we refer to e.g. \cite{BF11}, \cite{Fo15}. These papers use, among other things, the stochastic calculus of Poisson jump processes as a mathematical tool.  

Special attention has also been given to regularity of solutions to the Navier-Stokes equation, for example in \cite{DR14}, \cite{Mo05}. An overview on actual regularity results relative to solutions to a number of other PDEs in mathematical physics is provided in \cite{Ko19}. 
\medskip

From the technical side, certain forms of integration by parts or corresponding approximations can be crucial in order to investigate regularity of densities of a variety of stochastic variables or processes or regularity of solutions to classes of PDEs. For this we refer to \cite{Sa05}, \cite{Nu06}, \cite{BC11}, \cite{BC16}, \cite{BC17}, and \cite{Ba19}. 
\bigskip
 
{\bf $N$-particle Boltzmann type processes. } Throughout the paper, let $|\cdot|$ denote the Euclidean norm and let $\langle \cdot ,\cdot \rangle$ denote the corresponding scalar product where the dimension is always clear from the context.  The state space is a domain $(D\times V)^N \subseteq ({\mathbb R}^d\times {\mathbb R}^d)^N$ where $d\in\{2,3\}$, $D$ is convex, and the boundary $\partial D$ of $D$ is assumed to be smooth. The closure of the set $D$ is the set of all possible locations of a particle and $V=\{v:v_{min}<|v| <v_{max}\}$ for some $0\le v_{min}<v_{max}<\infty$ is the set of all its possible velocities. We are interested in stochastic processes of the form $Y_t=(x(t),v (t))\equiv ((x_1 (t), v_1 (t)),\ldots ,(x_N(t),v_N(t)))$, $t\ge 0$, with $x_i\in D$ and $v_i\in V$, $1\le i\le N$, on some probability space $(\Omega, {\cal F},P)$. For a fixed initial configuration randomness occurs, among others, at random times of two types 
\begin{eqnarray*}
0<\sigma_1(i,j)<\sigma_2(i,j)<\ldots\, ,\quad 1\le i<j\le N, 
\end{eqnarray*}
and 
\begin{eqnarray*}
0<\tau_1(i)<\tau_2(i)<\ldots\, ,\quad 1\le i\le N. 
\end{eqnarray*}
At times $\sigma_k(i,j)$ two of the $N$ particles, namely $i$ and $j$, change velocity due to collision. In other words, the velocities $v_i$ and $v_j$ perform a right-continuous random jump at time $\sigma_k(i,j)$. At times $\tau_l (i)$ one of the $N$ particles, namely $i$, hits the boundary $\partial D$ and is diffusively reflected. Thus, the velocity $v_i$ performs a right-continuous random jump at time $\tau_l(i)$. Between two successive times of $0<\sigma_1(i,j)<\sigma_2(i,j)<\ldots \, $, $1\le i<j\le N$, and $0<\tau_1(i)< \tau_2(i)<\ldots\, $, $1\le i\le N$, the velocities of all particles are constant and, consequently, their positions change linearly. In addition, jump times are assumed not to converge. Such processes $Y_t=(x(t),v (t))\equiv ((x_1(t),v_1 (t)),\ldots ,(x_N(t), v_N(t)))$, $t\ge 0$, we will call {\it $N$-particle Boltzmann type processes}. 

Our model is a slight modification of the ones introduced in \cite{CPW98} and \cite{RW05}. We assume particles to be small $d$-dimensional balls with diameter $\beta >0$ represented by their centers $x_1,\ldots ,x_N\in\mathbb{R}^d$. Let us explain the terms reflection and collision in our model. 
\medskip

{\it Reflection} of particle $i$ is, according to (2.24) of \cite{CPW98}, performed 
precisely at the time when $x_i$ hits the boundary $\partial D$. The new velocity 
including the new direction is randomly chosen according to a given velocity redistribution kernel, see Section \ref{sec:2}.  
\medskip

{\it Collision} is more subtle. We model particles as soft spheres such that any two particles $i$ and $j$ with $i<j$ can mutually penetrate each other upon collision. Assuming that $|x_i(s)-x_j(s)|=\beta$, $|x_i(s)+(t-s)v_i(s)-(x_j(s) +(t-s) v_j(s))|=\beta$, and $|x_i(s)+(u-s)v_i(s)-(x_j(s)+(u-s)v_j(s))|<\beta$ for some $0<s<t$ and all $s<u<t$ the actual time $\sigma(i,j)\in (s,s+\frac12(t-s))$ the particles $i$ and $j$ change velocity, and hence direction, is determined by a random control variable $\gamma(i,j)\in (0,1)$ triggering the velocity change. Let us suppose that, in this situation, $\sigma(i,j)=s+\frac12\, \gamma(i,j)\cdot(t-s)$ and that the distribution of $\gamma(i,j)$ is given by a probability density $g_\gamma$ on $(0,1)$ independent of $i$ and $j$. In this way we have implicitly introduced random collision places, namely the centers $x_i$ and $x_j$ of the particles $i$ and $j$ at the time $\sigma(i,j)$ of the change of their velocities. 

Furthermore, given the pre-impact velocities $v_i(\sigma(i,j)-)=v$ and $v_j(\sigma 
(i,j)-)=v'$, the post-impact velocities $v_i(\sigma(i,j))$ and $v_j(\sigma(i,j))$ 
are random. The distributions of $v_i(\sigma(i,j))$ and $v_j(\sigma(i,j))$ are 
expressed by the probability density $B_{v,v'}(e)$ of a certain random impact 
variable $\ve\equiv\ve(i,j)$ with outcomes 
\begin{eqnarray*}
e\in S^{d-1}_+(v-v'):=\left\{\t e\in S^{d-1}:\langle\t e,v-v'\rangle>0\right\}  
\end{eqnarray*} 
where $S^{d-1}$ denotes the unit sphere in $\mathbb{R}^d$. In fact, given outcomes $v$ of $v_i(\sigma(i,j)-)$ and $v'$ of $v_j(\sigma(i,j)-)$, for the outcomes $v^\ast$ of $v_i(\sigma(i,j))$ and $v'^\ast$ of $v_j(\sigma(i,j))$ we shall suppose 
\begin{eqnarray*}
v^\ast :=v-\langle e,v-v'\rangle\, e\quad\mbox{\rm and}\quad v'^\ast:=v'+\langle e, 
v-v'\rangle\, e\, . 
\end{eqnarray*} 

Let us set $(D^N)_\beta:=\left\{(x_1,\ldots ,x_N)\in D^N:x_i\in D,\ |x_i- x_j|> \beta\ \mbox{for}\ 1\le i<j\le N\right\}$.
In the paper, we shall also suppose that the random initial configuration $Y_0= (x_0, v_0)$ is given by a probability density with respect to the Lebesgue measure on $(D^N)_\beta\times V^N$, where $x_0\in (D^N)_\beta$ and $v_0\in V^N$. Let us also suppose there is a strict order between all jump times $\sigma_k(i,j)$ and $\tau_1(i')$, $1\le i <j\le N$, $1\le i'\le N$, $k,l\in\mathbb{N}$. We will therefore assume that on $(\Omega,{\cal F}, P)$ $P$-a.e. all outcomes display this property. 

It is standard in the literature to consider both random phenomena, collision places and post-impact velocities. While the post-impact velocities are usually modeled by means of a certain collision kernel, the treatment of the collision places varies in the literature. For example, in \cite{CPW98} and references therein a so-called spatial smearing function is introduced and in \cite{RW05} a certain scattering mechanism is explained. 
\medskip

{\bf Main result. } Under certain assumptions on the regularity of the characterizing terms, namely the initial density $p_0$, the collision kernel $B$, the random control variable $\gamma$, and the velocity redistribution kernel after hitting the boundary $M$, the main result of the paper can be formulated as follows.
\medskip

\nid
{\bf Theorem. } Fix $t>0$. \\ 
(a) There exists a finite collection $\mathcal{C}=\{\mathcal{C}_i\}$ of disjoint open subsets of $V^N$ satisfying 
\begin{eqnarray*}
\bigcup_{\mathcal{C}_i\in \mathcal{C}}\mathcal{C}_i\subseteq\{v(\omega) (t):\ \mbox{\rm all}\ \omega\} \subseteq\bigcup_{\mathcal{C}_i\in \mathcal{C}}\overline{\mathcal {C}_i} 
\end{eqnarray*}
such that the random variable $v(t)\equiv v(\omega)(t)$ possesses on every $\mathcal{C}_i\in \mathcal{C}$ a not necessarily bounded infinitely differentiable density $p_{v,i}$ with respect to the Lebesgue measure. \\ 
(b) There exists a finite collection $\mathcal{D}=\{\mathcal{D}_i\}$ of disjoint open subsets of $D^N$ satisfying 
\begin{eqnarray*}
\bigcup_{\mathcal{D}_i\in \mathcal{D}}\mathcal{D}_i\subseteq\{x(\omega) (t): \ \mbox{\rm all}\ \omega\} \subseteq\bigcup_{\mathcal{D}_i\in \mathcal{D}}\overline{\mathcal {D}_i} 
\end{eqnarray*}
such that the random variable $x(t)\equiv x(\omega)(t)$ possesses on every $\mathcal{D}_i\in \mathcal{D}$ a not necessarily bounded infinitely differentiable density $p_{x,i}$ with respect to the Lebesgue measure. 
\medskip

{\bf Organization of the paper. } Throughout the paper let us use $0/0=:0$. Furthermore, if on an open subset $U$ of some $\mathbb{R}^n$ a differentiable function $q:U\mapsto \mathbb{R}$ is given with $q(x)=0$ for some $x\in U$ then we set $\nabla \ln q(x)=0$. 
\medskip

In Section 2 we first properly define the characterizing terms $p_0,B,\gamma,M$ of the class of $N$-particle Boltzmann type processes we are dealing with. We then carry out preliminary finite dimensional partial integrations in order to prepare the infinite dimensional ones of the following section. 

In Section 3 we collect all items of the calculus we use in the paper. In particular we adopt a certain object from the Malliavin calculus on Riemannian manifolds, the flow equation of Theorem \ref{Theorem3.6}. Derivatives along  solutions to this equation result in certain directional derivatives based on which the gradient and, after integration by parts, the divergence are defined. 

In Section 4 we first establish families of trajectories which are specified by the order in which the times $\{\sigma_k(i,j):1\le i<j\le N\}\cup\{ \tau_l (i'):1\le i'\le N\}$ occur. These families will induce the sets $\mathcal{C}=\{ \mathcal {C}_i \}$ and $\mathcal{D}=\{\mathcal{D}_i\}$. Then the calculus developed in the previous sections is applied in order to establish an appropriately adjusted integration by parts implying the  existence of the smooth densities $p_{v,i}$ and $p_{x,i}$ of the main result. 

\section{Finite dimensional integration by parts}\label{sec:2}
\setcounter{equation}{0} 

An $N$-particle Boltzmann type process $Y_t$, $t\ge 0$, has the following sources of randomness, the initial configuration, the random impact variables $\ve\equiv \ve(i,j)$ inducing the post-collision velocity distributions, the random control variables $\gamma\equiv \gamma(i,j)$, and the random velocity redistributions after hitting the physical boundary. 

For the random initial configuration $Y_0=(x_0,v_0)$ we assume that its distribution is given by a probability density 
\begin{eqnarray*}
\qquad p_0\ \mbox{\rm on}\ (D^N)_\beta\times V^N\ \mbox{\rm with logarithmic gradient}\ \nabla\ln p_0\in L^2\left((D^N)_\beta\times V^N;p_0\cdot\l \right)  
\end{eqnarray*} 
where $\l$ denotes the Lebesgue measure on $(D^N)_\beta \times V^N$. Moreover, let $\ve_k(i,j)$ be the random impact variable determining the post-impact velocity at time $\sigma_k(i,j)$ and let $\chi$ denote the indicator function. According to \cite{CPW98} (2.20) we have 
\begin{eqnarray*}
&&\hspace{-.5cm}\frac{\displaystyle{P\left(\ve_k(i,j)\in de\left|\, v_i(\sigma_k (i,j) -)=v\, ,\, v_j(\sigma_k(i,j)-)=v'\right.\vphantom{l^1}\right)}}{de} \\ 
&&\hspace{.5cm}=\frac{B\left(v,v',e\right)}{\int_{S^{d-1}_+ (v-v')}B(v,v',e)\chi_{\{(v^\ast,v'^\ast)\in V\times V\}}\, de}=:B_{v,v'}(e) 
\end{eqnarray*}
where $de$ refers to the Riemann-Lebesgue measure on $S^{d-1}_+$ and $B$ is the collision kernel precisely defined in Definition \ref{Definition2.1} below. The last equation implies 
\begin{eqnarray*}
B\left(v,v',e\right)=0 \quad\mbox{\rm and }\quad B_{v,v'}(e)=0 \ \mbox{\rm  a.e.~on}\ \{e\in S^{d-1}_+ (v-v'): (v^\ast, v'^\ast) \not\in V\times V\}
\end{eqnarray*}
and is the motivation for part (v) of the subsequent Definition \ref{Definition2.1}. 
\begin{definition}\label{Definition2.1}{\rm 
Let $B(v_i,v_j,e)$ with domain $v_i,v_j\in V$, $e\in S^{d-1}$ denote the {\it collision kernel} applied to the particles $i$ and $j$, $1\le i<j\le N$, whenever the following conditions (i)-(vi) are satisfied. 
\begin{itemize}
\item[(i)] $B$ is a non-negative, bounded, and continuous function on its domain 
given above 
\item[(ii)] such that $B(v_i,v_j,\cdot)$ is symmetric in $v_i$ and $v_j$, $i\neq j$, and  
\item[(iii)] for any fixed $v,v'\in V$ the function $S^{d-1}_+\equiv S^{d-1}_+ (v-v')\ni e\mapsto B_{v,v'}(e)$ is continuously differentiable and compactly supported such that 
\begin{eqnarray*}
\inf_{v,v'\in V,\, e\in\, {\rm supp }B_{v,v'}(\cdot)}|\langle v-v',e\rangle|=:\kappa>0\, .
\end{eqnarray*}
\item[(iv)] 
$B(v^\ast,v'{}^\ast,e)=B(v,v', e)$ for all $v,v'\in V$ and $e\in S^{d-1}_+$ for 
which $(v^\ast,v'{}^\ast)\in V\times V$ and 
\item[(v)] we have 
\begin{eqnarray}\label{2.1}
&&\hspace{-.5cm}\frac{\displaystyle{P\left(v_i(\sigma_k(i,j))\in v^\ast(de)\left| 
\, v_i(\sigma_k(i,j)-)=v\, ,\, v_j(\sigma_k(i,j)-)=v'\right.\vphantom{l^1}\right)}} 
{de}\nonumber \\ 
&&\hspace{.5cm}=\frac{\displaystyle{P\left(v_j(\sigma_k(i,j))\in v'^\ast(de)\left| 
\, v_i(\sigma_k(i,j)-)=v\, ,\, v_j(\sigma_k(i,j)-)=v'\right.\vphantom{l^1}\right)}} 
{de}\nonumber \\ 
&&\hspace{.5cm}=\frac{\displaystyle{P\left(\ve_k(i,j)\in de\left|\, v_i(\sigma_k 
(i,j)-)=v\, ,\, v_j(\sigma_k(i,j)-)=v'\right.\vphantom{l^1}\right)}}{de}\nonumber \\ 
&&\hspace{.5cm}=B_{v,v'}(e)\, ,\quad v,v'\in V,\ e\in S^{d-1}_+(v-v')\, .\vphantom{\frac11}
\end{eqnarray}
\item[(vi)] For the logarithmic spherical gradient of $B_{v,v'}(\cdot)$ on $S^{d-1}_+$ we suppose that 
\begin{eqnarray*}
\left\|\nabla_e B_{v,v'}(\cdot)\right\|_{L^2\left(S^{d-1}_+;B_{v,v'}(\cdot)\, de\right)}\quad\mbox{\rm is bounded on $(v,v')\in V\times V$. }
\end{eqnarray*}
\end{itemize}
}
\end{definition} 

The next definition is dedicated to the probabilistic description of the 
precise times where two particles $i$ and $j$ change their velocities due 
to collision. For this and for later use we post the following hypotheses 
on the position of the particles, $((x_1,v_1),\ldots ,(x_N,v_N))$. 
\begin{itemize} 
\item[$p\, $(i)] It holds that $x_i(0)\in D$ and $|x_i(0)-x_{j}(0)|>\beta$ for all $1\le i<j\le N$.
\item[$p\, $(ii)] Let $1\le i<j\le N$ and $0<s<t$ such that $|x_i(s)-x_{j}(s)| =|x_i(t)-x_{j} (t)|= \beta$ and $|x_i(u)-x_{j}(u)|<\beta$ for all $u\in (s,t)$. For a third particle $i'\in \{1,\ldots ,N\}\setminus\{i,j\}$ we have $|x_i(s)-x_{i'}(s)| \neq \beta$ and $|x_j(s)-x_{i'}(s)|\neq\beta$. If 
\begin{eqnarray*}
|x_i(s)-x_{i'}(s)|>\beta\quad\mbox{\rm and}\quad |x_j(s)-x_{i'}(s)|>\beta  
\end{eqnarray*}
and 
\begin{eqnarray*}
|x_i(u)-x_{i'}(u)|<\beta\quad\mbox{\rm or}\quad |x_j(u)-x_{i'}(u)|<\beta\quad  
\mbox{\rm for some}\ u\in (s,t)
\end{eqnarray*}
then the particle $i'$ is not involved in any interaction with one of the particles $i$ or $j$ during the time interval $[s,t]$. However the particles $i$ and $j$ collide. 

\item[$p\, $(iii)] If $|x_i(t)-x_{j}(t)|=\beta$ for some $1\le i<j\le N$ and some $t>0$ then there is $\ve>0$ such that $|x_i(u)-x_{j}(u)|<\beta$ for all $u\in (t- \ve, t)$ or all $u\in (t,t+\ve)$. 
\end{itemize}

{\sc Remark. } Hypotheses $p\, $(i) -- $p\, $(iii) have been introduced for a proper construction of a probability space corresponding to the particle system. We mention that hypothesis $p\, $(ii) excludes a probability zero set. To avoid that three or more particles are involved in a collision, hypothesis $p\, $(ii) allows that particles continue without any interaction although their midpoints have come closer to each other than the particle diameter. This phenomenon already appears in established models, see for example \cite{CPW98}, (2.17) and (2.20). Furthermore, to simplify the argumentation in Section \ref{sec:4}, hypothesis $p\, $(iii) also excludes a probability zero set. 
\begin{definition}\label{Definition2.2}{\rm 
Let $1\le i<j\le N$. (a) If 
\begin{itemize}
\item[(i)] $|x_i(s)-x_j(s)|=\beta$, $|x_i(s)+(t-s)v_i(s)-(x_j(s)+(t-s)v_j(s))|= \beta$, and $|x_i(s)+(u-s)v_i(s)-(x_j(s)+(u-s)v_j(s))|<\beta$ for some $0<s<t$ as well as all $s<u<t$ 
\end{itemize}
and the particles $i$ and $j$ {\it collide} in the sense of hypothesis $p\, $(ii) then we say that $s$ is the {\it entrance time} of the collision. 
\\ 
(b) Let $s>0$ be the $k$-th of such entrance times in chronological order of collisions of the particles $i$ and $j$. We call $\sigma_k(i,j)\in (s,t)$ the {\it $k$-th collision time} of the particles $i$ and $j$ if the following holds. 
\begin{itemize}
\item[(i)] $(v_i(u),v_j(u))=(v_i(\sigma_k(i,j)-),v_j(\sigma_k(i,j)-)\neq (v_i( \sigma_k(i,j)),v_j(\sigma_k(i,j))$ for all $s\le u<\sigma_k(i,j)$.  
\item[(ii)] $\sigma_k(i,j)=s+\frac12\, \gamma_k(i,j)\cdot(t-s)$ where 
\item[(iii)] $\gamma_k(i,j)\in (\ve_\gamma,1-\ve_\gamma)$ for some $\ve_\gamma>0$ is a random variable independent of $Y_u$, $0\le u\le s$. 
\item[(iv)] The distribution of $\gamma_k(i,j)$ admits a probability density 
$g_\gamma$ independent of $i,j,k$, with supp$\, g_\gamma=[\ve_\gamma,1-\ve_\gamma]$,  $\, g_\gamma>0$ on $(\ve_\gamma,1-\ve_\gamma)$, and having a logarithmic derivative $(\ln g_\gamma)'\in L^2([0,1]; g_\gamma\cdot dx)< \infty$ where $dx$ refers to the Lebesgue measure on $[0,1]$. 
\end{itemize}
}
\end{definition}

{\sc Remark. } From $\beta=|x_i(s)-x_j(s)|=|x_i(s)+v_i(s) (t-s)-x_j(s) -v_j(s) (t-s)|$ we obtain $0=|v_i(s)-v_j(s)|^2(t-s)^2+2\langle x_i(s)-x_j(s),v_i(s)-v_j(s) \rangle\, (t-s)$ and thus 
\begin{eqnarray*}
&&\hspace{-.5cm}t-s=-2\, \frac{\langle x_i(s)-x_j(s),v_i(s)-v_j(s)\rangle}{|v_i(s)-v_j(s)|^2}\, .
\end{eqnarray*}
Therefore, 
\begin{eqnarray*}
&&\hspace{-.5cm}x_i(\sigma(i,j))-x_j(\sigma(i,j))= x_i(s) -x_j(s) + {\textstyle \frac12} \gamma(i,j)(v_i(s) -v_j(s)) (t-s) \\ 
&&\hspace{.5cm}= x_i(s) - x_j(s) - \gamma(i,j)\cdot\frac{v_i(s)-v_j(s)} {|v_i(s)- v_j(s) |}\cdot \left\langle x_i(s) -x_j(s), \frac{v_i(s)-v_j(s)} {|v_i(s)-v_j(s) |} \right\rangle\, .
\end{eqnarray*}
Given $x_i(\sigma(i,j))$, this implies a distribution of $x_j(\sigma(i,j))$ which involves the distributions of $x_i(s) - x_j(s)\ $, $v_i(s) - v_j(s)\ $, $\langle x_i(s) - x_j(s)\, ,\, v_i(s) - v_j(s)\rangle\ $, and $\gamma(i,j)$. Its density $h(x_i (\sigma(i,j)),\cdot)$ with respect to the Lebesgue measure on $D$ is of compact support. It is called {\it smearing function} in the literature, see for example \cite{CPW98}. However note that by our construction it may depend on $i$ and $j$ and the entrance time $s$ of the collision. 
\bigskip

The following definition describes the particle reflection at the boundary $\partial D$. It is adapted from \cite{CPW98}. 
\bigskip

%\nid
{\sc Notation. } Let $n(x)$ denote the outer normal vector on $\partial D$ at $x\in\partial D$. Abbreviate $n_l(i)\equiv n(x_i(\tau_l(i)))$. Let $\l_V$ denote the Lebesgue measure on $V$. 
\begin{definition}\label{Definition2.3}{\rm \hspace{0.2cm} 
An $N$-particle Boltzmann type process $Y_t=((x_1(t),v_1(t)),\ldots ,(x_N(t), $ $v_N(t)))$, $t\ge 0$, is said to have {\it diffusive boundary conditions} 
\begin{eqnarray}\label{2.2}
\frac{P(v_i(\left.\tau_l(i))\in dv\right|x_i(\tau_l (i)))}{dv}=M(x_i(\tau_l (i)); v)\, \langle v,n_l(i)\rangle\quad P\mbox{\rm -a.e.}\, , 
\end{eqnarray}
$1\le i\le N$, $l\in {\mathbb N}$, if the following holds. 
\begin{itemize}
\item[(i)] The {\it redistribution kernel} $M:\partial D\times V\mapsto {\mathbb R}_+$ is continuously differentiable on $\partial D \times V$ where 
\item[(ii)] $M(x;v)>0$ on $\{v\in V: \langle v/|v|,n(x)\rangle < -\ve\}$ and $\int_{ \langle v/|v|,n(x)\rangle<-\ve} M(x;v) \langle v,n(x)\rangle\, dv$ $=-1$ for all $x\in\partial D$ and some sufficiently small $\ve>0$, 
\item[(iii)] $M(x;v)=0$ on $\{v\in V: \langle v/|v|,n(x)\rangle \ge -\ve\}$, and  
\item[(iv)] $\|\nabla \ln [M(x;\, \cdot\, )\langle\, \cdot\, ,n(x)\rangle] \|_{L^2 (V;M(x;\cdot)\cdot\l_V )}$ is bounded on $x\in\partial D$.
\end{itemize} 
}
\end{definition}

Next we deal with a calculation which is important for the calculus of Section \ref{sec:3}. Since it just uses ordinary vector calculus we formulate it in form of an exercise. Here the symbols $\nabla_e$ and $\mbox{\rm div}_e$ denote the gradient and divergence with respect to the variable $e$, here as an element of $\mathbb {R}^d$. The multiplication sign is used for the usual product of scalars with matrices as well as matrices with matrices. For this, consider also the coordinate representation of $\nabla_e$ as a line vector. Below from (2.3) on, we may drop the multiplication signs from the notation.
\medskip 

\nid 
{\sc Exercise. } (a) Show that $\ \displaystyle{\frac{\partial v^\ast(e)} {\partial e}}= -\left[e \cdot (v-v')^T\right]-\langle e,v-v'\rangle\, I_d\quad$ where $I_d$ denotes the $d$-dimensional unit matrix. 
\medskip

\nid
(b) Verify $\ \displaystyle{\left(\frac{\partial v^\ast}{\partial e}(e) \right)^{-1}}= \frac{1} {2\langle e,v-v'\rangle^2}\cdot \left[e\cdot (v-v')^T \right]-\frac{1} {\langle e,v-v'\rangle}\, I_d\ $. 
\medskip

\nid
(c) Verify also $\ \displaystyle{\nabla_e\cdot\displaystyle{ \left(\frac{ \partial v^\ast}{\partial e}(e)\right)^{-1}}=\frac{d}2\cdot\frac{(v-v')^T}{\langle e,v-v' \rangle^2} }\ $. 
\medskip

\nid
(d) Show that 
\begin{eqnarray*}
&&\hspace{-.0cm}\mbox{\rm div}_e\left(\left(\frac{ \partial v^\ast(e)}{\partial e}\right)^{-1}\cdot \psi\left(v^\ast(e)\right)\cdot B_{v,v'}(e)\right) \\ 
&&\hspace{.5cm}=(\mbox{\rm div}\, \psi)\left( v^\ast(e)\right)\cdot B_{v,v'}(e)+\left\langle\psi\left(v^\ast(e)\right),\mathbf{d}(v,v',e) \right\rangle\vphantom{\left(\frac{\partial v^\ast(e)}{\partial e}\right)^{-1}}\cdot B_{v,v'}(e)\, ,\quad \psi\in C_b^\infty(\mathbb{R}^d;\mathbb{R}^d), 
\end{eqnarray*}
where 
\begin{eqnarray*}
\mathbf{d}(v,v',e)=\frac{(v-v')\cdot (d+\langle e,\nabla_e \ln B_{v,v'}(e)\rangle)-2 \langle e,v-v'\rangle\cdot\nabla_e \ln B_{v,v'}(e)}{2\langle e,v-v'\rangle^2}\, .
\end{eqnarray*}
%

%\smallskip

\medskip

%\nid
{\sc Notation. }(1) To ease the notation let us abbreviate below $v_i(\sigma_k) =v_i(\sigma_k(i,j))$, $v_j (\sigma_k) =v_j(\sigma_k(i,j))$, and $(v_i,v_j) (\sigma_k-) =\left(v_i(\sigma_k (i,j)-), v_j(\sigma_k(i,j)-)\right)$. We may also indicate dependence on $\ve_k (i,j)$ or its outcome $e\in S^{d-1}_+$, for example by $v_i(\sigma_k)\equiv v_i (\sigma_k)(\ve_k(i,j))$ or short $v_i(\sigma_k)\equiv v_i(\ve_k(i,j))$. 

\medskip

\nid
{\sc (2) } Let us suppose that ${\mathbb R}^{N\cdot d}$ is spanned by the ONB $e_k$, $k=1,\ldots ,N\cdot d$. Let $E_i$ denote the subspace of ${\mathbb R}^{N \cdot d}$ which is spanned by $e_{i\cdot d+1},\ldots ,e_{(i+1)\cdot d}$. Let $\nabla_i$ be the gradient with respect to $E_i$. In addition, let $\mbox{\rm div}_i$ denote the divergence with respect to $E_i$. This notation is actually not necessary for the subsequent calculations (\ref{2.3}) and (\ref{2.4}). It is introduced in order to ease the application of (\ref{2.3}) and (\ref{2.4}) in the proof of Proposition \ref{Proposition3.12}. 

%\medskip

%In order to take advantage of the exercise in (\ref{2.3}) below, introduce $S^{d -1}_{+ ;\ve} (v-v'):=\bigcup_{e\in S^{d-1}_+(v,v')}\{ae:a\in (-\ve,\ve)\}$ and let $v^\ast (ae) :=v^\ast(e)$ as well as $B_+(v-v')(ae):=B_+(v-v')(e)$, $e\in S^{d-1}_+ (v,v')$, $a\in (-\ve,\ve)$. The partial integration step of (\ref{2.3}) can then be carried out on $S^{d -1}_{+;\ve} (v-v')$ instead of $S^{d -1}_{+} (v-v')$ using the exercise. The partial integration on $S^{d -1}_{+} (v-v')$ is complete by dividing by $2\ve$ and letting $\ve\to 0$.

\medskip

Using part (d) of the above exercise, part (iii) of Definition \ref{Definition2.1}, and relation (\ref{2.1}) we obtain the following finite dimensional integration by parts. In this form it prepares the infinite dimensional integration by parts of Section \ref{sec:3}. Let $i\in\{1,\ldots,N\}$, $\vp_i\in C_b^\infty(\mathbb{R}^d)$, and $\psi_i\in C_b^\infty(\mathbb{R}^d;\mathbb{R}^d)$ such that supp$\, \vp_i \subset V$ or supp$\, \psi_i\subset V$. Denoting by $E$ the expectation with respect to $P$ it holds that 
\begin{eqnarray}\label{2.3}
&&\hspace{-.0cm}E\left[\left\langle\nabla_i\varphi_i(v_i(\ve_k(i,j)))\, ,\, \psi_i 
\left(v_i(\ve_k(i,j))\right)\vphantom{\dot{f}}\right\rangle\right]\vphantom{\int_{ 
v\in V}}\nonumber \\ 
&&\hspace{.5cm}=E\left[E\left[\left\langle\nabla_i\varphi_i(v_i(\ve_k(i,j)))\, ,\, 
\psi_i\left(v_i(\ve_k(i,j))\right)\vphantom{\dot{f}}\right\rangle\left|\vphantom 
{\dot{f}}\, (v_i,v_j)(\sigma_k -)\right.\right]\right]\vphantom{\int_{v\in V}} 
\nonumber \\ 
&&\hspace{.5cm}=\int_{(v,v')\in V\times V}\int_{e\in S^{d-1}_+(v-v')}\left\langle  
\nabla_i\varphi_i(v^\ast(e))\, ,\, \psi_i\left(v^\ast(e)\right)\vphantom{\dot{f}}
\right\rangle\, B_{v,v'}(e)\times\nonumber \\ 
&&\hspace{7cm}\times \, de\, P((v_i,v_j)(\sigma_k-)\in d(v,v'))\vphantom {\int}\nonumber \\ 
&&\hspace{.5cm}=\int_{(v,v')\in V\times V}\int_{e\in S^{d-1}_+(v-v')}\left\langle  \hspace{-.1cm}\left(\frac{\partial v^\ast(e)}{\partial e}\right)^T\, \nabla_i \varphi_i (v^\ast (e))\, ,\, \left(\frac{\partial v^\ast(e)}{\partial e}\right)^{-1}\, \psi_i\left( v^\ast(e)\right)\hspace{-.1cm}\right\rangle \times \nonumber \\ 
&&\hspace{7cm}\times\, B_{v,v'}(e)\, de\, P((v_i,v_j)(\sigma_k-)\in d(v,v'))\vphantom{\int_{v\in V}}\nonumber \\ 
&&\hspace{.5cm}=\int_{(v,v')\in V\times V}\int_{e\in S^{d-1}_+(v-v')}\left\langle  \nabla_e\varphi_i(v^\ast(e))\, ,\, \left(\frac{\partial v^\ast(e)}{\partial e}\right)^{-1}\, \psi_i\left(v^\ast(e) \right)\right\rangle\, B_{v,v'}(e)\times \nonumber \\ 
&&\hspace{7cm}\times \, de\, P((v_i,v_j)(\sigma_k-)\in d(v,v'))\vphantom {\int_{v\in V}}\nonumber \\
&&\hspace{.5cm}=-\int_{(v,v')\in V\times V}\int_{e\in S^{d-1}_+(v-v')}\varphi_i (v^\ast(e))\cdot \mbox{\rm div}_e\left(\left(\frac{\partial v^\ast(e)}{\partial e}\right)^{-1}\, \psi_i\left(v^\ast(e)\right)\, B_{v,v'}(e)\right) \times 
\nonumber \\ 
&&\hspace{7cm}\times\, de\, P((v_i,v_j)(\sigma_k-)\in d(v,v'))\vphantom{\int_{v\in V}}\nonumber \\
&&\hspace{.5cm}=-\int_{(v,v')\in V\times V}\int_{e\in S^{d-1}_+(v-v')}\varphi_i 
(v^\ast(e))\cdot\mbox{\rm div}_i\psi_i(v^\ast(e))\, B_{v,v'}(e)\times \nonumber \\ 
&&\hspace{7cm}\times\, de\, P((v_i,v_j)(\sigma_k-)\in d(v,v'))\vphantom{\int_{v\in V}}\nonumber \\
&&\hspace{1.0cm}-\int_{(v,v')\in V\times V}\int_{e\in S^{d-1}_+(v-v')}\varphi_i 
(v^\ast(e))\left\langle\psi_i\left(v^\ast(e)\right),\mathbf{d}(v,v',e)\right\rangle \times\nonumber \\ 
&&\hspace{7cm}\times\,  B_{v,v'}(e)\, de\, P((v_i,v_j)(\sigma_k-)\in d(v,v'))\vphantom{\int_{v\in V}}\nonumber \\ 
&&\hspace{.5cm}=E\left[\varphi_i(v_i(\ve_k(i,j)))\vphantom{\dot{f}}\right.\times \nonumber \\ 
&&\hspace{1.0cm}\times\left. \left(\vphantom{\dot{f}}-\mbox{\rm div}_i\psi_i(v_i (\ve_k(i,j)))-\left\langle\psi_i(v_i(\ve_k(i,j))),\mathbf{d}\left((v_i,v_j) (\sigma_k-), \ve_k(i,j)\vphantom{l^1}\right)\right\rangle\vphantom{\dot{f}}\right) \right]\, . 
\end{eqnarray}
Furthermore, using (\ref{2.2}), another standard calculation shows that for $i\in\{1, \ldots, N\}$ and $\vp_i\in C_b^\infty(\mathbb{R}^d)$ as well as $\psi_i\in C_b^\infty (\mathbb{R}^d;\mathbb{R}^d)$ with supp$\, \vp_i\subset V$ or supp$\, \psi_i\subset V$ we have 
\begin{eqnarray*}
&&\hspace{-.5cm}E\left[\left.\left\langle\nabla_i\varphi_i(v_i(\tau_l))\, ,\, \psi_i (v_i(\tau_l))\vphantom{l^1}\right\rangle\vphantom{\dot{f}}\right|x_i(\tau_l (i))\right] \\ 
&&\hspace{.5cm}=\int_V \left\langle\nabla\varphi_i(v)\, ,\, \psi_i(v)\vphantom {l^1}\right\rangle \, M(x_i(\tau_l(i));v)\langle v,n_l(i)\rangle \, dv \\ 
&&\hspace{.5cm}=-\int_V\varphi_i(v)\cdot\mbox{\rm div}_i\psi_i(v)\, M(x_i 
(\tau_l(i));v)\langle v,n_l(i)\rangle\, dv\nonumber \\ 
&&\hspace{1.0cm}-\int_V\varphi_i(v)\cdot\left\langle\psi_i(v)\, ,\, \nabla_i\left[ 
M(x_i(\tau_l(i));v)\langle v,n_l(i)\rangle\vphantom{l^1}\right]\right\rangle\, dv\, . 
\end{eqnarray*}
It follows that 
\begin{eqnarray}\label{2.4} 
&&\hspace{-.5cm}E\left[\left\langle\nabla_i\varphi_i(v_i(\tau_l))\, ,\, \psi_i( v_i(\tau_l))\vphantom{l^1}\right\rangle\vphantom{\dot{f}}\right]=E\left[\varphi_i (v_i(\tau_l))\cdot\left(\vphantom{\dot{f}}-\mbox{\rm div}_i\psi_i(v_i(\tau_l))\right. 
\right.\nonumber \\ 
&&\hspace{.5cm}\left.\left. -\left\langle\psi_i(v_i(\tau_l))\, ,\, \nabla_i\left[\ln\left(M(x_i(\tau_l(i));v_i(\tau_l) )\, \langle v_i(\tau_l) ,n_l(i)\rangle \vphantom {l^1} \right)\right]\right\rangle\vphantom{\dot{f}}\right)\right]\, . 
\end{eqnarray}

\section{Infinite dimensional integration by parts and related calculus}\label{sec:3}
\setcounter{equation}{0}

In this section we introduce a reduced version of trajectories of an $N$-particle Boltzmann type process. These {\it reduced} trajectories carry the same information as the original ones. 

It is the objective to define derivatives of functions of reduced trajectories of an $N$-particle Boltzmann type process by means of derivatives along appropriate flows. Such a flow takes values in the set of all reduced trajectories of the $N$-particle Boltzmann-type process. For this, Definition 3.1 provides an appropriate sequence representation of the directions which are tangential to such a flow, see also Theorem 3.5. Furthermore, Definition 3.2, introduces a sequence representation of the reduced trajectories of the $N$-particle Boltzmann type process, see also Lemma 3.3. The most challenging technical issue in order to construct derivatives along a flow as above, is to incorporate reflections of trajectories on the boundary of the physical space. This leads to Definition 3.4. As Theorem 3.5 shows, these efforts pay off and result in Proposition 3.9 (a) in the proper definition of a pre-version of a directional derivative. 
All other derivatives introduced in Section 3 are based on the one of Proposition 3.9 (a).
\bigskip

{\sc Notation. (1)} For an $N$-particle Boltzmann type process $Y_t=((x_1(t),v_1 (t)), \ldots ,$ $(x_N(t), v_N(t)))$, $t\ge 0$, let us consider the set $D_Y$ of all functions $Y(\omega): [0,\infty)\mapsto (D\times V)^N$, $\omega\in\Omega$, constructed in Section \ref{sec:1} satisfying hypotheses $p\, $(i)-$p\, $(iii) of Section \ref{sec:2}. We will call these functions $Y(\omega)$, $\omega\in\Omega$, the {\it trajectories of} $Y$. This means that we suppose that all trajectories $Y(\omega)$ of $Y$ are cadlag. For the probability space $(\Omega, \mathcal{F},P)$ we assume that there is an identification of $\omega \in\Omega$ with the trajectory $Y(\omega)$ and that $\mathcal{F}$ is generated by the cylinder sets over the trajectories. We mention that restrictions, for example by the cut-offs of Definition \ref{Definition2.1} (iii) or Definition \ref{Definition2.3} (iii), may  result in $P$-zero sets. 
\medskip

\nid 
{\sc (2)} Let $D_{[0,\infty)\mapsto {\mathbb R}^{N\cdot d}}$ denote the set of all ${\mathbb R}^{N\cdot d}$-valued cadlag functions on $[0,\infty)$. Furthermore, let $K_{[0,\infty) \mapsto {\mathbb R}^{N\cdot d}}$ denote the set of all functions in $D_{[0,\infty) \mapsto {\mathbb R}^{N\cdot d}}$ which are piecewise constant such that jump points do not converge on $[0,\infty)$. Let ${\mathbb D}:=\left({\mathbb R}^{N\cdot 2d};K_{[0, \infty)\mapsto {\mathbb R}^{N\cdot d}}\right)$. 
\medskip

\nid
{\sc (3)} For an $N$-particle Boltzmann type process $Y_t=((x_1(t),v_1 (t)), \ldots ,$ $(x_N(t), v_N(t)))$, $t\ge 0$, let us introduce the {\it reduced process} 
\begin{eqnarray*}
\mathbb{Y}\equiv \left(\mathbb{Y}(0); (v_1(t),\ldots ,v_N(t))_{t\ge 0}\vphantom{l^1}\right)\quad\mbox{\rm where}\quad\mathbb{Y}(0):=Y_0\, . 
\end{eqnarray*}
For the trajectory  $Y(\omega)\equiv y(\omega)(t)=((x_1 (\omega) (t),v_1(\omega) (t)), \ldots ,(x_N(\omega) (t),v_N(\omega)(t)))$, $t\ge 0$, $\omega\in\Omega$, define the {\it reduced trajectory} constructed from $Y(\omega)$ by 
\begin{eqnarray*}
{\mathbb Y}(\omega):=\left(y(\omega)(0);(v_1(\omega)(t),\ldots ,v_N(\omega)(t))_{t\ge 0}\vphantom{l^1}\right)\, . 
\end{eqnarray*}

Denote by $D_{\mathbb Y}$ the set of all reduced trajectories, i.e. $D_{\mathbb Y}= \{{\mathbb Y}\in {\mathbb D}:Y\in D_Y\}$. We mention that $Y(\omega)$ can be recovered from ${\mathbb Y}(\omega)$ by $x_i(\omega)(t)=x_i (\omega)(0)+\int_0^t v_i(\omega)(s) \, ds$, $1\le i\le N$, $t\ge 0$. Observe also that the definition of the term {\it reduced trajectory} ${\mathbb Y}(\omega)$ is compatible with the set of test functions of the form $(\ref{3.1})$ introduced below. 
\medskip 

\nid
{\sc (4)} Since the map $\Omega\ni\omega\mapsto {\mathbb Y}(\omega)\in D_{\mathbb Y}$ is a bijection, let us also identify $\omega\in\Omega$ with the reduced trajectory ${\mathbb Y}(\omega)$. Denote by $\mu_{\mathbb Y}$ the image measure of $P$ and let $\mathcal {F}_{\mathbb Y}$ be the image of $\mathcal{F}$ under this identification. If there is no ambiguity, we may drop the symbol $\omega$ from the notation. To simplify the notation, we shall write $E[f(\mathbb{Y}(\omega))]$ or $E[f(\mathbb{Y})]$ instead of $\int f(\mathbb{Y})\, d\mu_{\mathbb Y}$. 

\bigskip

\nid
{\sc (5)} Let ${\bf e}_i\equiv(({\bf e}_i)_1,\ldots ,({\bf e}_i)_{N\cdot d})\in {\mathbb R}^{N \cdot d}$ be the vector with $({\bf e}_i)_n=1$ if $n\in \{i\cdot d+1,\ldots ,(i+1) \cdot d\}$ and $({\bf e}_i)_n=0$ otherwise. For $b\equiv (b_1, \ldots ,b_d)\in {\mathbb R}^d$ let $b\circ {\bf e}_i\equiv ((b\circ{\bf e}_i)_1,\ldots ,(b\circ{\bf e}_i)_{N \cdot d})$ given by $(b\circ {\bf e}_i)_n=b_r({\bf e}_i)_n$ if $n=i\cdot d+r$ and $r\in \{1,\ldots ,d\}$, and $(b\circ {\bf e}_i)_n=0$ otherwise. For $a^{(1)}, a^{(2)} \in {\mathbb R}^d$, $a:=(a^{(1)},a^{(2)})\in {\mathbb R}^{2d}$, and $1\le i<j\le N$ introduce $a\circ ({\bf e}_i+{\bf e}_j):=a^{(1)}\circ {\bf e}_i+a^{(2)}\circ {\bf e}_j$. For $c_1,\ldots ,c_N\in {\mathbb R}^d$, $c:=(c_1, \ldots ,c_N)\in {\mathbb R}^{N\cdot d}$, and $1\le i<j\le N$ let $\overline{c(i)} :=c-c_i\circ {\bf e}_i$ and $\overline{c(i,j)}:=c-(c_i,c_j)\circ ({\bf e}_i+{\bf e}_j)$. For an arbitrary set $S$ let $\ell (S)$ denote the space of all sequences over $S$. 

\begin{definition}\label{Definition3.2}{\rm 
(a) For $\omega\in\Omega$ and ${\mathbb Y}(\omega)\in D_{\mathbb Y}$, let $\rho_0:=0$ and $\rho_m\equiv\rho_m({\mathbb Y}(\omega))$, $m\in {\mathbb N}$, be the a.e. strictly increasing sequence of random times such that 
\begin{eqnarray*}
\left\{\rho_m:m\in {\mathbb N}\right\}=\left\{\sigma_k(i,j):1\le i<j\le N,\ k\in 
{\mathbb N}\right\}\cup\left\{\tau_l(i'):1\le i'\le N,\ l\in {\mathbb N}\right\}\, . 
\end{eqnarray*}
(b) Introduce the sequence spaces 
\begin{eqnarray*}
\mathcal{A}:=\left(\ell ({\mathbb R}^d)\right)^{\frac{(N-1)N}{2}}\times\left(\ell ({\mathbb R}^d)\right)^N\quad\mbox{\rm and}\quad\mathcal{G}:=\left(\ell ({\mathbb R})\right)^{\frac{(N-1)N}{2}}\, . 
\end{eqnarray*}
(c) Given $\omega\in\Omega$ as well as the associated ${\mathbb Y}(\omega)\in D_{ \mathbb Y}$, and $h(0)=(h_1(0),\ldots ,h_N(0))\in {\mathbb R}^{N\cdot d}$, and moreover
\begin{eqnarray*}
\alpha=\left(\left((a_{i,j;k})_{k\in {\mathbb N}}\right)_{\{1\le i<j\le N\}}\, ; \, \left((b_{i';l})_{l\in {\mathbb N}}\right)_{\{1\le i'\le N\}}\right)\in\mathcal{A} 
\end{eqnarray*}
define the function 
\begin{eqnarray*} 
(h(\alpha,{\mathbb Y}(\omega))(t))_{t\ge 0}\equiv\left(\left( {h}_1(\alpha,{\mathbb Y}(\omega))(t),\ldots ,{h}_N(\alpha,{\mathbb Y}(\omega))(t)\right) \right)_{t\ge 0}\in K_{[0,\infty)\mapsto {\mathbb R}^{N\cdot d}} 
\end{eqnarray*}
by 
\begin{itemize}
\item[(i)] $h(\alpha,{\mathbb Y}(\omega))(t):=h(0)$, $0\le t<\rho_1$,  
\item[(ii)] for $m\in\mathbb{N}$, 
\begin{eqnarray*}
h(\alpha,{\mathbb Y}(\omega))(t):=\left\{
\begin{array}{l}
\overline{h(\alpha,{\mathbb Y}(\omega))(\rho_m-)(i,j)}+(a_{i,j;k},a_{j,i;k})\circ 
({\bf e}_i+{\bf e}_j)\vphantom{\displaystyle\int} \\  
\hphantom{yy} {\rm if}\ \rho_m\le t<\rho_{m+1}\ {\rm and}\ \rho_m=\sigma_k 
(i,j)\ \mbox{\rm for some}\ 1\le i<j\le N\vphantom{\displaystyle\int} \\ 
\overline{h(\alpha,{\mathbb Y}(\omega))(\rho_m-)(i')}+b_{i';l}\circ {\bf e}_{i'} 
\vphantom{\displaystyle\int}  \\ 
\hphantom{yy} {\rm if}\ \rho_m\le t<\rho_{m+1}\ {\rm and}\ \rho_m=\tau_l(i') 
\ \mbox{\rm for some}\ 1\le i'\le N\vphantom{\displaystyle\int}
\end{array}
\right. \hspace{-.5cm}
\end{eqnarray*}
where 
\begin{eqnarray*}
a_{j,i;k}:=h_i(\alpha,{\mathbb Y})(\rho_m-)+h_j(\alpha,{\mathbb Y})(\rho_m-)-a_{i,j;k} 
\end{eqnarray*}
if $\rho_m=\sigma_k(i,j)$ for some $1\le i<j\le N$ and $k\in\mathbb{N}$.
\end{itemize} 
}
\end{definition}
\begin{definition}\label{Definition3.3}{\rm 
(a) Introduce the {\it set $\hat{V}$ of equivalence classes} $\hat{v}\equiv 
\{v,-v\}$, $v\in V$, i.e.  
\begin{eqnarray*}
\hat{V}:=\left\{\hat{v}\equiv\{v,-v\}:v\in V\right\}\quad\mbox{\rm and the 
{\it absolute value}}\quad\left|\hat{v}\right|:=|v|\, . 
\end{eqnarray*} 
Moreover, introduce the sequence spaces  
\begin{eqnarray*} 
\mathcal{A}_V:=\left(\ell (V)\right)^{\frac{(N-1)N}{2}}\times \left(\ell (\hat{V})\right)^N\quad\mbox{\rm and}\quad\mathcal{G}_X:=\bigcup_{\ve_\gamma\in (0,\frac12)}\left(\ell ((\ve_\gamma, 1-\ve_\gamma)\right)^{\frac{(N-1)N}{2}}\, . 
\end{eqnarray*} 
(b) Fix $y(0)=(x(0),h(0))$ where $x(0)=(x_1(0),\ldots ,x_N(0))\in (D^N)_\beta$ as well as $h(0)=(h_1(0),\ldots ,$ $h_N(0)) \in V^N$. Choose furthermore 
\begin{eqnarray*}
\eta=\left(\left((a_{i,j;k})_{k\in {\mathbb N}}\right)_{\{1\le i<j\le N\}}\, ;\,  
\left((\widehat{b_{i';l}})_{l\in {\mathbb N}}\right)_{\{1\le i'\le N\}}\right)\in 
\mathcal{A}_V   
\end{eqnarray*}
as well as 
\begin{eqnarray*}
\g=\left(\gamma_k(i,j)\right)_{k \in {\mathbb N},\, 1\le i<j\le N}\in \mathcal{G}_X\, .
\end{eqnarray*}
Let $\hat{h}(\eta,y)$ be defined by the following. 
\begin{itemize}
\item[(i)] $\hat{h}(\eta,y)(0):=h(0)$, $0\le t<\rho_1$.  
\item[(ii)] Given $\hat{h}(\eta,y)(s)=\left(\hat{h}_1(\eta,y)(s), 
\ldots ,\hat{h}_N(\eta,y)(s)\right)\in V^N$, $s\in [0,t]$, let 
\begin{eqnarray*}
x_i(t)=x_i(0)+\int_0^t\hat{h}_i(\eta,y)(s)\, ds\, ,\quad 1\le i\le N,\ t 
\ge 0, 
\end{eqnarray*}
and assume for the following that $x(t)=\left(x_1(t),\ldots ,x_N(t)\right) \in \overline{D}^N$. 
\item[(iii)] Let $\left\{\rho_m:m\in {\mathbb N}\right\}$ denote the set of times introduced in Definition \ref{Definition3.2} (a). Given $x(t)=\left(x_1(t),\ldots ,x_N(t)\right)\in\overline{D}^N$, $t \in [0,s]$ such that $x_i(s)\in\partial D$ for some $0\le i\le N$, assume that $s=\tau_l(i)\equiv\rho_m$ for some $l,m\in \mathbb{N}$. 

Given $x(t)=\left(x_1(t),\ldots ,x_N(t)\right)\in\overline {D}^N$, $t \in [0,s]$, such that $|x_i(s)-x_j(s)|=\beta$, $|x_i(u)-x_j(u)|>\beta$ for all $u\in(s- \varepsilon, s)$ and some $\varepsilon>0$ as well as some $0\le i<j\le N$, and that the particles $i$ and $j$ collide in the sense of hypothesis $p(ii)$, assume that $s+\gamma_k(i,j) =\sigma_k(i,j) \equiv \rho_m$ for some $k,m\in \mathbb{N}$.
\smallskip

Recall that $n(x_i(\rho_m))$ denotes the outer normal vector on $\partial D$ at $x_i(\rho_m)$ if $x_i(\rho_m) \in \partial D$. Note that for the equivalence class $\widehat {b_{i';l}}\in\hat{V}$, the term 
\begin{eqnarray*}
\wt{ b_{i';l}}:=-\mbox{\rm sign}\left\langle b_{i';l}, n(x_i(\rho_m))\vphantom{l^1} \right\rangle\cdot b_{i';l}
\end{eqnarray*}
is independent of the choice of the representative $b_{i';l}\in\widehat{b_{i'; 
l}}$ or $-b_{i';l}\in\widehat{b_{i';l}}$. 
\item[(iv)] For $m\in\mathbb{N}$ let 
\begin{eqnarray*}
\hat{h}(\eta,y)(t):=\left\{
\begin{array}{l}
\overline{\hat{h}(\eta,y)(\rho_m-)(i,j)}+(a_{i,j;k},a_{j,i;k})\circ ({\bf e}_i+{\bf e}_j)\vphantom{\displaystyle\int} \\  
\hphantom{yy} {\rm if}\ \rho_m\le t<\rho_{m+1}\ {\rm and}\ \rho_m=\sigma_k 
(i,j)\ \mbox{\rm for some}\ 1\le i<j\le N\vphantom{\displaystyle\int} \\ 
\overline{\hat{h}(\eta,y)(\rho_m-)(i')}+\wt{ b_{i';l}}\circ {\bf e}_{i'}\vphantom{\displaystyle\int}  \\ 
\hphantom{yy} {\rm if}\ \rho_m\le t<\rho_{m+1}\ {\rm and}\ \rho_m=\tau_l(i')\ \mbox{\rm for some}\ 1\le i'\le N\vphantom{\displaystyle\int}
\end{array}
\right.  
\end{eqnarray*}
where 
\begin{eqnarray*}
a_{j,i;k}:=\hat{h}_i(\eta,y)(\rho_m-)+\hat{h}_j(\eta,y)(\rho_m-)-a_{i,j;k} 
\end{eqnarray*}
if $\rho_m=\sigma_k(i,j)$ for some $1\le i<j\le N$.
\end{itemize}
}
\end{definition} 
%

%\nid
{\sc Remark. } Implicitly this definition introduces a forward in time construction of a reduced trajectory from $(x(0),v(0))\in (D^N)_\beta\times V^N$, $\g=\left( \gamma_k(i,j)\right)_{k \in {\mathbb N},\, 1\le i<j\le N}\in \mathcal {G}_X$, and $\eta\in\mathcal{A}_V$, since the times $\rho_m$, $m\in\mathbb{N}$, are generated along this forward in time construction. The subsequent lemma follows.
\begin{lemma}\label{Lemma3.4}
Let $y(0)=(x(0),v(0))\in (D^N)_\beta\times V^N$, $\g=\left(\gamma_k(i,j) \right)_{ k\in {\mathbb N},\, 1\le i<j\le N}\in \mathcal {G}_X$, and $\eta\in\mathcal{A}_V$. There is a unique 
\begin{eqnarray*}
{\mathbb Y}\equiv\left({\mathbb Y}(0);(v_1(t),\ldots ,v_N(t))_{t\ge 0}\vphantom{l^1} 
\right)=:{\mathbb F}(y(0);\g;\eta)\in D_{\mathbb Y} 
\end{eqnarray*} 
such that with 
\begin{eqnarray*}
x_i(t):=x_i(0)+\int_0^t v_i(s)\, ds\, ,\quad 1\le i\le N,\ t\ge 0, 
\end{eqnarray*} 
we have the following. 
\begin{itemize}
\item[(i)] ${\mathbb Y}(0)=y(0)$.  
\item[(ii)] If $|x_i(s)-x_j(s)|=\beta$, $|x_i(s)+(t-s)v_i -(x_j(s)+(t-s)v_j)| =\beta$, $|x_i(s)+(u-s)v_i-(x_j(s)+(u-s)v_j)|<\beta$ for $s<u<t$, and the particles $i$ and $j$ collide in the sense of hypothesis $p(ii)$ such that $[s,t]$ is the first interval after $\sigma_{k-1}(i,j)$ with this property, then $\sigma_k(i,j) =s+ \frac12\, \gamma_k(i,j)\cdot (t-s)$, $k\in {\mathbb N}$. Here, $\sigma_0(i,j) :=0$. 
\item[(iii)] $(v_1(t),\ldots ,v_N(t))=\hat{h}(\eta,{\mathbb Y})(t)$, $t\ge 0$. 
\end{itemize}
The map ${\mathbb F}:(D^N)_\beta\times V^N\times\mathcal{G}_X\times\mathcal{A}_V \mapsto D_{\mathbb Y}$ is a bijection. 
\end{lemma}
\medskip

%\nid
{\sc Notation. (1) } This lemma gives rise to add dependencies to the notation. Given ${\mathbb Y}\in D_{ \mathbb Y}$ we may write 
\begin{eqnarray*}
\g\equiv\g({\mathbb Y})\in\mathcal{G}_X\, ,\quad\eta\equiv\eta({\mathbb Y})\in 
\mathcal{A}_V\, ,\quad\mbox{\rm and}\quad y(0)={\mathbb Y}(0)\in (D^N)_\beta\times V^N. 
\end{eqnarray*}
In particular, the coefficients of part (iv) of Definition \ref{Definition3.3} can be written as $a_{i,j;k}\equiv a_{i,j;k}(\mathbb{Y})$ and $\wt{ b_{i';l}}\equiv \wt{ b_{i';l}}(\mathbb{Y})$. Furthermore, in order to mark the jump times we may write 
\begin{eqnarray*}
\sigma_k(i,j)\equiv\sigma_k(i,j)({\mathbb Y})\, ,\quad\tau_l(i')\equiv\tau_l(i')({\mathbb Y})\, ,\quad\mbox{\rm and}\quad \rho_m\equiv\rho_m(\mathbb Y). 
\end{eqnarray*}
{\sc (2) } Note that with the norm 
\begin{eqnarray*}
&&\hspace{-.5cm}\|(y(0);\g,\alpha)\|^2_{{\mathbb H}}:=|y(0)|^2+\sum_{1\le i\le N}\sum_{l\in {\mathbb N}}\left|b_{i;l}\right|^2+\sum_{1\le i<j\le N}\sum_{k\in {\mathbb N}}\left(\left|a_{i,j;k}\right|^2+\left(\gamma_k(i,j)\right)^2\vphantom{l^1} \right)\, , \\  
&&\hspace{.5cm}\mbox{\rm where}\quad y(0)\in {\mathbb R}^{2N\cdot d},\ \g=\left( 
\gamma_k(i,j) \right)_{k\in{\mathbb N},\, 1\le i<j\le N}\in\mathcal{G},\vphantom 
{\sum_{1\le i<j\le N}} \\  
&&\hspace{.5cm}\mbox{\rm and}\quad\alpha=\left(\left((a_{i,j;k})_{k\in{\mathbb N}} 
\right)_{\{1\le i<j\le N\}}\, ;\, \left((b_{i';l})_{l\in {\mathbb N}}\right)_{\{1 
\le i'\le N\}}\right)\in\mathcal{A},\vphantom{\sum_{1\le i<j\le N}} 
\end{eqnarray*}
the set 
\begin{eqnarray*}
{\mathbb H}:=\left\{(y(0);\g,\alpha)\in{\mathbb R}^{2N\cdot d}\, \times\mathcal{G}\times\mathcal{A}:\|(y(0);\g,\alpha)\|_{\mathbb H}<\infty\right\}
\end{eqnarray*}
becomes a separable Hilbert space. 
\begin{definition}\label{Definition3.5}{\rm 
(a) Let $y(0)\in (D^N)_\beta\times V^N$, $\g\in\mathcal{G}_X$, and $\eta\in\mathcal {A}_V$ as in part (b) of Definition \ref{Definition3.3} and consider $\mathbb{Y} \in D_{\mathbb{Y}}$ constructed as in Lemma \ref{Lemma3.4} from $y(0)$, $\g$, and $\eta$. Furthermore, let 
\begin{eqnarray*}
\t \eta\equiv \t \eta (\mathbb{Y}):=\left(\left((a_{i,j;k}(\mathbb{Y}))_{k\in {\mathbb N}}\right)_{\{1\le i<j\le N\}}\, ;\, \left((\wt {b_{i';l}}(\mathbb{Y}))_{l\in {\mathbb N}}\right)_{\{1\le i'\le N\}}\right)\, .
\end{eqnarray*}
In addition, for 
\begin{eqnarray*}
\alpha=\left(\left((a_{i,j;k}^\ast)_{k\in {\mathbb N}}\right)_{\{1\le i<j\le N\}}\, ;\, \left((b_{i';l}^\ast)_{l\in {\mathbb N}}\right)_{\{1\le i'\le N\}}\right)\in\mathcal{A}\, ,
\end{eqnarray*}
let 
\begin{eqnarray*}
\widehat{b_{i';l}}\plus b_{i';l}^\ast:=\left\{\wt{b_{i';l}}(\mathbb{Y})+b_{i';l}^\ast\, ,\, -(\wt{b_{i';l}}(\mathbb{Y})+b_{i';l}^\ast)\right\}\in\hat{V}
\end{eqnarray*}
provided that $\wt{b_{i';l}}(\mathbb{Y})+b_{i';l}^\ast\in V$. Note the non-symmetry of the operation $``\plus"$ which, in this sense, we use for an element $\hat{b} \in \hat{V}$ on the left-hand side of this operation and an element $b^\ast\in \mathbb{R}^d$ on the right-hand side of this operation. Correspondingly, for $s\in \mathbb{R}$ introduce 
\begin{eqnarray*}
\eta\+ s\alpha := \left(\left((a_{i,j;k} +sa_{i,j;k}^\ast)_{k\in {\mathbb N}}\right)_{\{1\le i<j\le N\}}\, ;\, \left((\widehat{b_{i';l}}\plus sb_{i';l}^\ast)_{l\in {\mathbb N}}\right)_{\{1\le i'\le N\}}\right)\in\mathcal{A}_V\, ,
\end{eqnarray*}
provided that $a_{i,j;k}+sa_{i,j;k}^\ast\in V$ and $\wt{b_{i';l}}(\mathbb{Y}) +sb_{i';l}^\ast\in V$, $1\le i<j\le N$, $1\le i'\le N$, $k,l\in \mathbb{N}$.
}
\end{definition}

{\sc Remark. } As a consequence of Lemma \ref{Lemma3.4} the map 
\begin{eqnarray*} 
D_{\mathbb Y}\ni {\mathbb Y}\mapsto\left( \mathbb{Y}(0);\g(\mathbb{Y});\t \eta (\mathbb{Y})\right)\in \left\{\left( \mathbb{X}(0);\g(\mathbb{X});\t \eta(\mathbb{X})\right):{\mathbb X}\in D_{\mathbb Y}\right\}
\end{eqnarray*} 
is bijective. Its inverse we shall denote by $\t {\mathbb F}$. 
\bigskip

\nid
{\bf Definition 3.4.} (continuation)
(b) For $\ve>0$ and $(-\ve,\ve)\ni u\mapsto {\mathbb X}_{s+u}\in D_{\mathbb Y}$ we say that ${\mathbb X}_\cdot$ is {\it differentiable} in $s\in {\mathbb R}$ if we have the following. There exists $H_s=((x_s(0),h_s(0));\c_s;\alpha_s)\in {\mathbb H}$ with $\alpha_s=\left(( ((a_s)_{ i,j; k}^\ast)_{k \in {\mathbb N}})_{\{1\le i<j\le N\}}\, ;\, (((b_s)_{i';l }^\ast )_{l \in {\mathbb N}})_{\{1\le i'\le N\}} \right) \in\mathcal{A}$ and $\c_s = \left(((c_s)_k(i,j) \right.$ $\left. )_{k \in {\mathbb N}}\right)_{\{1 \le i<j\le N\}}\in \mathcal{G}$ such that in the notation of part (a) we have 
\begin{eqnarray*}
&&\hspace{-.5cm}\frac{d}{ds}{\mathbb X}_{s}(0) =(x_s(0),h_s(0))\, ,\quad\frac{d}{ds}\gamma_k(i,j)({\mathbb X}_{s}) = (c_s)_k(i,j)\, ,
\end{eqnarray*}
and 
\begin{eqnarray*}
&&\hspace{-.5cm}\frac{d}{ds}\t \eta({\mathbb X}_{s})\equiv\lim_{u\to 0}\frac{\t \eta({\mathbb X}_{s+u})-\t \eta({\mathbb X}_{s})}{u}=\alpha_s
\end{eqnarray*}
componentwise, i.e., 
\begin{eqnarray*}
&&\hspace{-.5cm}\frac{d}{ds}a_{i,j;k}({\mathbb X}_{s}) = (a_s)^\ast_{i,j;k}\, ,\quad\frac{d}{ds}\wt{b_{i';l}}({\mathbb X}_{s}) = (b_s)^\ast_{i';l}\, ,
\end{eqnarray*}
$1\le i<j\le N$, $1\le i'\le N$, $k,l\in\mathbb{N}$. In this case we shall write 
\begin{eqnarray*}
&&\hspace{-.5cm}\lim_{u\to 0}\frac{\t {\mathbb F}^{-1}\left({\mathbb X}_{s+u}\right)-\t {\mathbb F}^{-1}\left({\mathbb X}_s\right)}{u} \\ 
&&\hspace{.5cm}=\lim_{u\to 0}\frac{\left(\mathbb{X}_{s+u}(0);\g(\mathbb{X}_{s+u});\t \eta(\mathbb{X}_{s+u})\vphantom{f^f}\right)-\left(\mathbb{X}_s(0);\g(\mathbb{X}_s);\t \eta(\mathbb{X}_s)\vphantom{f^f}\right)}{u}=H_s 
\end{eqnarray*}
or 
\begin{eqnarray*}
\left\{
\begin{array}{rcl}
\displaystyle{\vphantom{\sum^1}\dot{{\mathbb X}}_s\equiv\frac{\partial{\mathbb X}_s} {\partial s}} &=& \left((x_s(0),h_s(0));h(\alpha_s,{\mathbb X}_s)\vphantom {\displaystyle {l^1}}\right)  \\ 
\displaystyle{\vphantom{\sum^1}\dot{\g}({\mathbb X}_s)\equiv\frac{\partial}{\partial s}\g({\mathbb X}_s)} &=& \c_s
\end{array}
\right.\, . 
\end{eqnarray*}
\medskip 

\nid 
(c) Given ${\mathbb Y}\in D_{\mathbb Y}$, let $T_{{\mathbb Y}}$ be the set of all {\it admissible} directions from ${\mathbb Y}$. That is the set of all 
\begin{eqnarray*}
\h=\left(a(0);\c\, ;h(\alpha,{\mathbb Y})\vphantom{\displaystyle{l^1}}\right)\, ,
\end{eqnarray*}
where 
\begin{eqnarray*}
a(0)=(x(0),h(0)), \quad h(\alpha,{\mathbb Y})(0)=h(0)\in\mathbb{ R}^{N \cdot d}, \quad \mbox{\rm and}\quad H:=(a(0);\c\, ;\alpha)\in {\mathbb H}  
\end{eqnarray*}
such that the following holds. There is a family $\{{\mathbb Y}_s\in D_{\mathbb Y}:|s|< \ve\}$ for some $\ve>0$ such that ${\mathbb Y}=\mathbb{Y}_0$ and 
\begin{eqnarray*} 
(a(0);\c\, ;\alpha)=\lim_{s\to 0}\frac{({\mathbb Y}_s(0);\g ({\mathbb Y}_s);\t \eta({\mathbb Y_s}))-({\mathbb Y}(0);\g({\mathbb Y});\t \eta(\mathbb {Y}))}{s}\, .  
\end{eqnarray*}
\medskip

\nid
{\sc Notation. (1) } Referring to a particular $\omega\in\Omega$, i.e. to a particular reduced trajectory $\mathbb{Y}(\omega)\in D_{\mathbb{Y}}$, let $T_{\mathbb{Y} (\omega)}$ denote the set of all admissible directions from $\mathbb{Y}(\omega)$. 
\medskip 

\nid 
{\sc (2) } Let $\mathbb{H}^f$ denote the space of all $(a(0);\c\, ; \alpha) \in\mathbb{H}$ where $a(0)=(x(0), h(0))$ and $x(0),h(0)\in {\mathbb R}^{N \cdot d}$ such that $\alpha $ and $\c$ have just finitely many elements different from zero. 

\medskip

\begin{theorem}\label{Theorem3.6} 
Let $\omega\in\Omega$, ${\mathbb Y}(\omega)\in D_{\mathbb Y}$, and $(a(0);\c\, ; \alpha) \in\mathbb{H}^f$ in the sense of the last Notation (2). Furthermore, let $\h=\left(a(0);\c\, ;h(\alpha,\mathbb{Y} (\omega) ) \vphantom{ \displaystyle{l^1}}\right)$ where $h(\alpha, {\mathbb Y} (\omega))(0) =h(0)$. Set $y(0):=y(\omega)(0)$, $\g:=\g({\mathbb Y(\omega )})$, and $\eta:=\eta ({\mathbb Y}(\omega))$. Then, for $P$-a.e. $\omega\in\Omega$,
\begin{eqnarray*} 
\h=\left(a(0);\c\, ;h(\alpha,\mathbb{Y} (\omega)) \vphantom{ \displaystyle{l^1}}\right)\in T_{\mathbb {Y}(\omega)} 
\end{eqnarray*}
and the {\rm flow equation} 
\begin{eqnarray*} 
\left\{
\begin{array}{rcll}
\displaystyle{\frac{\partial\sigma^\alpha_s(\omega)}{\partial s}}&=&\left(a(0); h( \alpha,\sigma^\alpha_s(\omega))\vphantom{\displaystyle{l^1}}\right) & \\ \displaystyle{\vphantom{\sum^1}\dot{\g}(\sigma^\alpha_s(\omega))}&=&\c &   \\ 
\displaystyle{\vphantom{\sum^1}\sigma^\alpha_0(\omega)}&=& {\mathbb Y}(\omega) &
\end{array}
\right. 
\end{eqnarray*}
has a unique solution $(-\ve,\ve)\ni s\mapsto\sigma^\alpha_s(\omega)\in D_{\mathbb Y}$ for some $\ve\equiv\varepsilon (\omega,\h)>0$. This solution is given by 
\begin{eqnarray*} 
\sigma^\alpha_s(\omega)={\mathbb F}\left(sa(0)+y(0);\g+s\c\, ;\eta\+s\alpha\right)\, ,\quad s\in (-\ve,\ve). 
\end{eqnarray*}
\end{theorem}
Proof. {\it Step 1 } As in Definition \ref{Definition3.5} (a), let us use the notation 
\begin{eqnarray*}
\alpha=\left(\left((a_{i,j;k}^\ast)_{k\in {\mathbb N}}\right)_{\{1\le i<j\le N\}}\, ;\, \left((b_{i';l}^\ast)_{l\in {\mathbb N}}\right)_{\{1\le i'\le N\}}\right)\in\mathcal{A}\, .
\end{eqnarray*}
Consider ${\mathbb Y}(\omega)=(y(0);(v_1(t), \ldots ,v_N(t) )_{t\ge 0})\in D_{\mathbb Y}$ and fix $i\in\{1,\ldots ,N\}$. Let $\sigma^\alpha_s (\omega)$ be given as above. By Definition \ref{Definition2.3} (ii), (iii) we have $P$-a.e. $\langle n(x_i (\tau_l(i) (\mathbb{Y}))),\wt {b_{i;l}}(\mathbb{Y})/|\wt {b_{i;l}} (\mathbb{Y})| \rangle <-\ve$ for all $i\in\{1,\ldots ,N\}$ and $l\in \mathbb{N}$, for some $\ve>0$. Furthermore, $\tau_l(i) (\sigma^\alpha_s(\omega)) \stack{s\to 0}{\lra}\tau_l (i) (\mathbb{Y})$. Thus $\wt{b_{i;l}}(\sigma^\alpha_s (\omega))=\wt {b_{i;l}} (\mathbb{Y})+sb_{i;l }^\ast\in V$ as well as $\langle n(x_i (\tau_l(i) (\sigma^\alpha_s (\omega)))),(\wt {b_{i;l}}(\mathbb{Y})+sb_{i;l }^\ast)/|\wt {b_{i;l}} (\mathbb{Y})+sb_{i;l }^\ast| \rangle <-\ve$ for some $\ve\equiv\ve(i;l) >0$ and all $s\in (-\ve ,\ve)$. 
 
By $v_{min}<|v_i|<v_{max}$ we may suppose $|v_i(t)|\in [v_{min} + \ve, v_{max} -\ve]$, $0\le t\le T$ for any $T>0$ and some $\ve\equiv\ve(T)>0$. Recall now that $(a(0);\c\, ; \alpha) \in\mathbb{H}^f$. As a consequence, $\eta\+ s\alpha$ is well-defined for some $\ve>0$ and $s\in (-\ve ,\ve)$ by Definition \ref{Definition3.5} (a). 

Moreover, by $(a(0);\c\, ; \alpha) \in\mathbb{H}$ there is a bound on $|c_k(i,j)|$ uniformly on $1\le i<j\le N$ and $k\in \mathbb{N}$. According to Definition \ref{Definition2.2} (b) (iii), there exists $\ve>0$ such that for all $s\in (-\ve,\ve)$ and all $1\le i<j\le N$ as well as $k_0\in \mathbb{N}$ and $k\le k_0$ it holds that $s c_k(i,j) +\gamma_k(i,j)\in (\ve_\gamma, 1-\ve_\gamma)$. Thus, there is $\, \ve>0$ such that for $s\in (-\ve ,\ve)$ a unique reduced trajectory $\sigma^\alpha_s (\omega)\in D_{\mathbb Y}$ can be constructed from $y(0)+sa(0)$, $\g+s\c\, $, and $\eta\+s\alpha$ as introduced by Definition \ref{Definition3.3} and Lemma \ref{Lemma3.4}. In other words, there is $\, \ve>0$ such that for $s\in (-\ve ,\ve)$
\begin{eqnarray*} 
\sigma^\alpha_s (\omega):={\mathbb F}\left(y(0)+sa(0);\g+s\c\, ;\eta\+s\alpha\right)
\end{eqnarray*}
and, according to the remark after Definition \ref{Definition3.5} (a), 
\begin{eqnarray*} 
\sigma^\alpha_s(\omega)\equiv\t{\mathbb F}\left(\sigma^\alpha_s(\omega)(0);\g (
\sigma^\alpha_s(\omega));\t \eta(\sigma^\alpha_s(\omega))\right)=\t{\mathbb F}
\left(y(0)+sa(0);\g+s\c\, ;\wt{\eta\+ s\alpha}\right)\, . 
\end{eqnarray*}

\nid
{\it Step 2 } We show that the construction of Step 1 complies with Definition \ref{Definition3.5} (b) and (c). For $\, \ve>0$ as in the end of Step 1 and $s\in (-\ve ,\ve)$ we observe 
\begin{eqnarray*} 
&&\hspace{-.5cm}L(s):=\frac{(\sigma^\alpha_s(\omega)(0);\g (\sigma^\alpha_s (\omega)); \t \eta(\sigma^\alpha_s(\omega)))-({\mathbb Y}(0); \g({\mathbb Y});\t \eta(\mathbb {Y}))}{s}-(a(0);\c\, ;\alpha) \\ 
&&\hspace{.5cm}=\frac{(0;0;\t \eta(\sigma^\alpha_s (\omega)))-(0;0;\t \eta( \mathbb {Y}))}{s}-(0;0;\alpha) \\ 
&&\hspace{.5cm}=\left(0;0;\left(0;\left(\frac{\wt{b_{i;l}}(\sigma^\alpha_s (\omega))-\wt{b_{i;l}}(\mathbb {Y})}{s} -b^\ast_{i;l}\right)_{1\le i\le N,\ l\in {\mathbb N}} \right)\right)
\end{eqnarray*}
since $a_{i,j;k}+sa_{i,j;k}^\ast\in V$ by Step 1 and therefore $\frac{a_{i,j; k} +sa_{i,j;k}^\ast -a_{i,j;k}}{s} -a_{i,j; k}^\ast=0$, $1\le i<j\le N$, $k\in {\mathbb N}$. 

By Step 1 it holds that $\wt{b_{i;l}}(\sigma^\alpha_s (\omega))=\wt {b_{i;l}} (\mathbb {Y}) +sb_{i;l }^\ast$ for some $\ve \equiv\ve(i;l)>0$ and $s\in (-\ve ,\ve)$. We obtain 
\begin{eqnarray*} 
&&\hspace{-.5cm}\lim_{s\to 0}L(s)=0 
\end{eqnarray*}
which, in particular says $\h=\left(a(0);\c\, ;h(\alpha,\mathbb{Y} (\omega)) \vphantom{ \displaystyle{l^1}}\right)\in T_{\mathbb {Y}(\omega)}$. The theorem follows. 
\qed 
\medskip

According to Definition \ref{Definition3.5} (c) and the calculation in Step 2 of the proof of Theorem \ref{Theorem3.6} we have the following. 
\begin{corollary}\label{Corollary3.7} 
Let $\mathbb{Y}\in D_\mathbb{Y}$. The map $H=(a(0);\c\, ;\alpha)\mapsto\h =\left( a(0);\c\, ;h(\alpha,{\mathbb Y}) \vphantom{\displaystyle{l^1}} \right)$ is an injection $\mathbb{H}^f\mapsto T_{\mathbb {Y}}$.
\end{corollary}
{\sc Notation } Denote $T^f_{\mathbb {Y}}:=\{\left( a(0);\c\, ;h(\alpha,{\mathbb Y}) \vphantom{\displaystyle{l^1}} \right):(a(0);\c\, ;\alpha)\in\mathbb{H}^f\}$. 
\begin{definition}\label{Definition3.8}{\rm (a) For ${\mathbb Y}\in D_{\mathbb Y}$  introduce by 
\begin{eqnarray*}
\|\h(H,{\mathbb Y})\|_{T_{{\mathbb Y}}}:=\|H\|_{\mathbb H}\, ,\quad H\in {\mathbb H}^f,
\end{eqnarray*}
an isometry between $\mathbb {H}^f$ and $T^f_{{ \mathbb Y}}$. The corresponding completion $T^c_{{\mathbb Y}}$ of $T^f_{{ \mathbb Y}}$ as well as $T_{{\mathbb Y}}$ defines a norm $\| \cdot \|_{T^c_{{\mathbb Y}}}$ on $T^c_{{\mathbb Y}}$ and a related inner product $\langle\cdot, \cdot \rangle_{T^c_{{\mathbb Y}}}$. 
\bigskip

\nid 
{\sc Notation. (1) } For the elements of $T^c_{{\mathbb Y}}$, we take over the notation of Definition \ref{Definition3.5} (c). 
\medskip

\nid
{\sc (2) } Given ${\mathbb Y}\in D_\mathbb{Y}$, by the previous corollary and part (a) of the current definition there is a one to one correspondence between $H =(a(0); \c \, ;\alpha) \in {\mathbb H}$ and $\h=\left(a(0);\c\, ;h (\alpha, {\mathbb Y}) \vphantom{ \displaystyle {l^1}} \right) \in T^c_{{\mathbb Y}}$ which we will indicate by 
\begin{eqnarray*}
\h\equiv\h(H,{\mathbb Y})\, . 
\end{eqnarray*}

\nid
{\bf Definition 3.7.} (continuation)  (b) Define ${\mathbf T}_{\mathbf Y}:=\{T^c_{ \mathbb Y}: {\mathbb Y}\in D_{ \mathbb Y}\}$. For a measurable map $D_{\mathbb Y} \ni {\mathbb Y} \mapsto H({\mathbb Y}) \in {\mathbb H}$ we say that the induced map 
\begin{eqnarray*}
\bh\equiv \bh_{H({\mathbb Y})}:D_{\mathbb Y}\ni {\mathbb Y}\mapsto \h(H({\mathbb Y}),{\mathbb Y})\in T^c_{\mathbb Y} 
\end{eqnarray*}
{\it takes values in} ${\bf T}_{\mathbf Y}$ and indicate it with ``${;{\bf T}_{\mathbf Y}}$'' in the notation of the function space the map $\bh\equiv \bh_{H ({\mathbb Y})}$ belongs to.
}
\end{definition} 

{\sc Notation. (1) } Let $L^2(D_{\mathbb Y},\mu_{\mathbb Y};{\mathbf T}_{\mathbb Y})$ denote the set of all measurable maps (equivalence classes) $D_{\mathbb Y}\ni {\mathbb Y}\mapsto \Phi({\mathbb Y})\in T^c_{{\mathbb Y}}$ such that ${\mathbb Y} \mapsto\|\Phi( {\mathbb Y}) \|_{T^c_{{ \mathbb Y}}}$ is a measurable real function on $(D_{\mathbb Y},\mathcal {F}_{\mathbb Y})$ and 
\begin{eqnarray*} 
\|\Phi\|^2_{L^2}\equiv \|\Phi\|^2_{L^2(D_{\mathbb Y},\, \mu_{\mathbb Y};\, {\mathbf T}_{\mathbb Y})} :=E\left[\|\Phi({\mathbb Y})\|^2_{T^c_{{\mathbb Y}}}\right] <\infty.
\end{eqnarray*}
The space $\left(L^2(D_{\mathbb Y},\mu_{\mathbb Y};{\mathbf T}_{\mathbb Y}),\| \cdot \|^2_{L^2} \right)$ is a Hilbert space. We denote the inner product by $\langle \cdot ,\cdot \rangle_{L^2}\equiv\langle\cdot ,\cdot\rangle_{L^2(D_{\mathbb Y},\, \mu_{\mathbb Y};\, {\mathbf T}_{\mathbb Y})}$. 
\bigskip

\nid
% {\bf Definition 3.7.} (continuation)
% (c) Introduce by 
%
% \begin{eqnarray*}
% \|\bh_{H(\mathbb{Y})}\|_{\mathbf{T}_\mathbb{Y}}:=\|\h(H(\mathbb{Y}),\mathbb{Y})\|_{L^2} 
% \end{eqnarray*}
% 
% a norm on $\mathbf{T}_{\mathbb{Y}}:=\{\bh_{H(\mathbb{Y})}\in\t {\mathbf{T}}_{ \mathbb {Y}}:\h(H(\mathbb{Y}),\mathbb{Y})\in L^2(D_{\mathbb Y},\mu_{\mathbb Y};\t {\bf T}_{\mathbb Y}) \}$ and a corresponding inner product $\langle \cdot, \cdot \rangle_{ \mathbf{T }_{\mathbb{Y}}}$. With this, the notation can be simplified to $L^2 (D_{\mathbb Y},\mu_{\mathbb Y};\t {\bf T}_{\mathbb Y})\equiv L^2(D_{\mathbb Y},\mu_{\mathbb Y}; {\bf T}_{\mathbb Y})$. 
% \bigskip

\nid 
{\sc (2) } For an open subset $S$ of some Euclidean space $\mathbb{R}^n$, let $C_b^\infty(S)$ denote the space of all bounded real infinitely differentiable functions on $S$. In addition, let $C_0^\infty(S)$ denote the space of all compactly supported real infinitely differentiable functions on $S$.

Let $f_0\in C_b^\infty((D^N)_\beta\times V^N)$, $f^\gamma_{i,j;k}\in C_b^\infty (0,1)$, and $f_{i,j;k},f_{i';l}\in C_b^\infty(V)$, $1\le i<j\le N$, $1\le i'\le N$, $k,l \in {\mathbb N}$. Moreover, let $f_{i';l}(v)=f_{i';l}(-v)$, $v\in V$, $1\le i'\le N$, $l\in {\mathbb N}$. For reduced trajectories ${\mathbb Y}(\omega)$, finite sets  $C_F$ of indices $(i,j;k)$, finite sets $R_F$ of indices $(i';l)$, and test functions
\begin{eqnarray}\label{3.1}
&&\hspace{-1cm}F({\mathbb Y}(\omega))\equiv F({\mathbb Y}):=f_0(y(0))\cdot\prod_{ (i,j;k)\in C_F}\left(f^\gamma_{i,j;k}\left(\gamma_k(i,j)\vphantom{l^1}\right)\cdot f_{i,j;k}\left(v_i (\sigma_k (i,j))\vphantom{l^1}\right)\vphantom{\dot{f}}\right)\times \nonumber \\ 
&&\hspace{.0cm}\times\prod_{(i';l)\in R_F}f_{i';l}\left( v_{i'}(\tau_l (i'))\vphantom{l^1}\right)
\end{eqnarray}
we may introduce the set $\hat{C}_b^\infty (D_{\mathbb Y})$ of all such test functions $F$ given by (\ref{3.1}). Furthermore, let $\hat{ C}_0^\infty(D_{\mathbb Y})$ denote the subset of all $F\in\hat{C}_b^\infty (D_{\mathbb Y})$ with $f_0\in C_0^\infty ((D^N)_\beta\times V^N)$, $f^\gamma_{i,j;k} \in C_0^\infty (0,1)$, and $f_{i,j;k},f_{i';l}\in C_0^\infty (V)$ for all $(i,j;k)\in C_F$ and all $(i';l)\in R_F$. 
\bigskip

\nid
{\sc (3) } Let $\t C_b^\infty(D_{\mathbb Y})$ denote the space of all finite linear combinations of elements of $\hat{C}_b^\infty(D_{\mathbb Y})$. Furthermore, let $\t C_0^\infty (D_{\mathbb Y})$ denote the space of all finite linear combinations $F=\sum_{r=1}^n\hat{F}_r$ of elements $\hat{F}_r\in \hat{C}_0^\infty(D_{\mathbb Y})$ such that for all $p,q\in\{1,\ldots ,n\}$
\begin{eqnarray*}
C_{\hat{F}_p}=C_{\hat{F}_q}\quad\mbox{\rm and}\quad R_{\hat{F}_p}=R_{\hat{F}_q}\, .
\end{eqnarray*}

\nid
{\sc (4) } As just introduced, for $F=\sum_{r=1}^n\hat{F}_r\in \t C_b^\infty (D_{ \mathbb Y})$ there is a unique set $I(F)$ of indices $(i,j;k)$ as well as $(i';l)$ which are necessary to define $F$. If $F\in\t C_0^\infty (D_{\mathbb Y})$ then these are precisely the indices belonging to $C_{\hat{F}_p}$ or $R_{\hat{F}_p}$, i.e. $I(F)=C_{\hat{F}_p}\cup R_{\hat{F}_p}$ for any $p\in\{1,\ldots ,n\}$. If $F \in \t C_b^\infty (D_{\mathbb Y})\setminus \t C_0^\infty (D_{\mathbb Y})$ then $I(F) \subseteq C_{\hat{F}_p}\cup R_{\hat{F}_p}$ for all $p$.

Denote by $\mathbb {H}_F$ the linear subspace of $\mathbb{H}$ consisting of all $(a(0);\c;\alpha)$ for which the components corresponding to the indices not belonging to $I(F)$ are zero. Keep in mind that $F\in\t C_b^\infty (D_{\mathbb Y})$ implies $\mathbb {H}_F\subset \mathbb {H}^f$.

\begin{lemma}\label{Lemma3.9}
The sets $\t C_0^\infty(D_{\mathbb Y})$ and $\t C_b^\infty (D_{{\mathbb Y}})$ are dense in $L^2(D_{\mathbb Y},\mu_{\mathbb Y})$. 
\end{lemma}
Proof. {\it Step 1 } Introducing the abbreviation $\left\{\gamma_k(i,j) \vphantom {l^1}\right\}\equiv\left(\gamma_k(i,j)\vphantom{l^1}\right)_{{k\in {\mathbb N}}\atop { 1\le i<j\le N}} $ and in the same manner $\left\{v_i(\sigma_k (i,j)) \vphantom{l^1} \right\}$ as well as $\left\{\hat{v}_i(\tau_l(i)) \vphantom{l^1} \right\}$ let 
\begin{eqnarray*} 
S_{\mathbb Y}\equiv S_{\mathbb Y}(\omega) :=\left((x(0),v(0)),\left\{\gamma_k(i,j) \vphantom {l^1}\right\},\left\{v_i(\sigma_k(i,j))\vphantom{l^1}\right\},\left\{ \hat{v}_i(\tau_l(i))\vphantom{l^1}\right\}\vphantom{\dot{f}}\right)\, .
\end{eqnarray*}
Recall that the map $\mathbb{F}^{-1}$, 
\begin{eqnarray*} 
D_{\mathbb Y}\ni {\mathbb Y}\mapsto S_{\mathbb Y}\in (D^N)_\beta\times V^N\times \mathcal{G}_X\times\mathcal{A}_V=:\mathcal{S}_{\mathbb Y}  
\end{eqnarray*}
is a bijection by Lemma \ref{Lemma3.4}. Let $(\mathcal{S}_{\mathbb Y},\mathcal{F}_{ \cal S},\mu)$ be the image of $(D_{\mathbb Y},\mathcal{F}_{\mathbb Y},\mu_{\mathbb Y})$ under this map. We note that the $\sigma$-algebra $\mathcal{F}_{\cal S}$ coincides with the $\sigma$-algebra on $\mathcal{S}_{\mathbb Y}$ generated by the random sequence $S_{\mathbb Y}$. Introduce $C_b^\infty(\hat{V}) :=\left\{f\in C_b^\infty(V): f(v)=f(-v) \right\}$, 
\begin{eqnarray*}
&&\hspace{-.5cm}\mathcal{C}_{\mathcal{B}}^\infty:=C_b^\infty \left((D^N)_\beta \times V^N\right) \times \left(\left(C_b^\infty((0,1))\vphantom{l^1}\right)^{\frac{(N-1)N}{2}}\right)^{\mathbb N}\times \\ 
&&\hspace{1cm}\times\left(\left(C_b^\infty(V)\vphantom{l^1}\right)^{\frac{(N-1)N} {2}}\right)^{\mathbb N}\times\left(\left(C_b^\infty(\hat{V})\vphantom{l^1}\right)^N \right)^{\mathbb N}\, , 
\end{eqnarray*}
and 
\begin{eqnarray*}
\hat{\mathcal{C}}_b^\infty:=\left\{\Phi\in\mathcal{C}_{\mathcal{B}}^\infty:\Phi\ \mbox{\rm consists of finitely many components}\right\}    
\end{eqnarray*}
where we identify the components of $\Phi\in\hat{\mathcal{C}}_b^\infty$ as factors, the product of which we identify with $\Phi$. Define also 
\begin{eqnarray*}
\mathcal{C}_b^\infty:=\left\{\sum_{r=1}^n\Phi_r:\Phi_r\in\hat{\mathcal{C}}_b^\infty, \ n\in {\mathbb N}\right\}\, . 
\end{eqnarray*}
Comparing the construction of ${\mathbb Y}(\omega)$ from the components of $\mathcal {S}_{\mathbb Y}$ in Definition \ref{Definition3.3} (b) and Lemma \ref{Lemma3.4} with the definition of $S_{\mathbb Y}$ by components of $\mathcal{S}_{\mathbb Y}$ we even get $\mathcal{F}_{\cal S}=\sigma \left( {\mathcal{C }_b^\infty}\right)$. 

The aim of this step is to demonstrate that $\t C_b^\infty(D_{\mathbb Y})$ is dense in $L^2(D_{\mathbb Y},\mu_{\mathbb Y})$. Sufficient for this is to show that $\mathcal {C}_b^\infty$ is dense in $L^2\left( \mathcal{S}_{\mathbb Y},\mu \right)$. For the latter, we adapt the proof of \cite{DS01}, Lemma 2.7, where we apply the Monotone Class theorem in the form of \cite{Pr05}, Theorem 8 of Chapter I. 

Let $H$ denote the set of all bounded functions in the closure of $\mathcal{C }_b^\infty$ with respect to $L^2\left(\mathcal{S}_{\mathbb Y},\mu\right)$. Then $H$ is a vector space with $\1\in H$ where $\1=1$ on $\mathcal{S}_{\mathbb Y}$. Furthermore, we observe that for any increasing sequence of non-negative functions in $H$, $0\le f_1\le f_2\le\ldots\ $ such that the limit $f:=\lim_{n\to\infty}f_n$ is bounded we have $f\in H$. 

In addition, $\mathcal{C}_b^\infty$ is closed under multiplication. These properties of $\mathcal{C}_b^\infty$ and $H$ imply now that $H$ contains all bounded $\sigma \left({\mathcal{C}_b^\infty}\right)=\mathcal{F}_{\cal S}$-measurable functions. Thus $H$ is dense in $L^2\left(\mathcal{S}_{\mathbb Y}, \mu \right)$. It follows that $\mathcal {C}_b^\infty$ is dense in $L^2 \left( \mathcal {S}_{\mathbb Y},\mu\right)$ which entails that $\t C_b^\infty (D_{\mathbb Y})$ is dense in $L^2(D_{\mathbb Y}, \mu_{\mathbb Y})$. 
\medskip 

\noindent
{\it Step 2 } Let $C_0^\infty(\hat{V}):=\left\{f\in C_0^\infty(V): f(v)=f(-v) \right\}$
and 
\begin{eqnarray*}
&&\hspace{-.5cm}\mathcal{C}_{\mathcal{O}}^\infty:=C_0^\infty \left((D^N)_\beta\times V^N\right)\times\left(\left(C_0^\infty((0,1))\vphantom{l^1}\right)^{\frac{(N-1)N}{2}}\right)^{\mathbb N}\times 
 \\ 
&&\hspace{1cm}\times\left(\left(C_0^\infty(V)\vphantom{l^1}\right)^{\frac{(N-1)N}{2}}\right)^{\mathbb N}\times\left(\left(C_0^\infty(\hat{V})\vphantom{l^1}\right)^N\right)^{\mathbb N}\, , 
\end{eqnarray*}
and 
\begin{eqnarray*}
\hat{\mathcal{C}}_0^\infty:=\left\{\hat{\Phi}\in\mathcal{C}_{\mathcal{O}}^\infty: \hat{\Phi}\ \mbox {\rm consists of finitely many factors}\right\} 
\end{eqnarray*}
as well as, with the same identifications as in Step 1, 
\begin{eqnarray*}
\mathcal{C}_0^\infty:=\left\{\sum_{r=1}^n\hat{\Phi}_r:\hat{\Phi}_r\in\hat{\mathcal{C}}_0^\infty, \ n\in {\mathbb N}\right\} 
\end{eqnarray*}
where we suppose that for all $p,q\in\{1,\ldots ,n\}$
\begin{eqnarray}\label{3.2*}
&&\hat{\Phi}_p=\vp^{(p)}_0\prod \vp^{\gamma,(p)}_{i,j;k}\cdot\prod \vp^{(p)}_{i,j;k}\cdot \prod \vp^{(p)}_{i;l} \nonumber \\ 
&&\hspace{-1.1cm}\quad\mbox{\rm implies}\quad\hat{\Phi}_q =\vp^{(q)}_0\prod \vp^{\gamma,(q)}_{i,j;k}\cdot\prod \vp^{(q)}_{i,j;k}\cdot \prod \vp^{(q)}_{i;l} \\ 
&&\hspace{-0.6cm}\ \mbox{\rm for the same indices $(i,j;k)$ as well as $(i;l)$} \vphantom{\prod \vp^{(q)}_{i,j;k}} \nonumber
\end{eqnarray}
and compactly supported factors $\vp^{(p)}_0$, $\vp^{\gamma,(p) }_{i,j;k}$ $\vp^{(p)}_{i,j;k}$, $\vp^{(p)}_{i;l}$ as well as $\vp^{(q)}_0$, $\vp^{\gamma, (q)}_{i,j;k}$ $\vp^{(q)}_{i,j;k}$, $\vp^{(q)}_{i;l}$, cf. (3) of the last block of notations. 
\medskip

The set $\mathcal{C}_0^\infty$ is a dense subset of $\mathcal{C}_b^\infty$ with respect to the topology in $L^2\left( \mathcal{S}_{\mathbb Y},\mu\right)$ since any $\Phi=\sum_{r=1}^n{\Phi}_r\in\mathcal{C}_b^\infty$ can be adjusted to property (\ref{3.2*}) by appending factors equal to one to the summands ${\Phi}_r$. By Step 1, $\mathcal{C}_0^\infty$ is dense in $L^2\left(\mathcal{S}_{\mathbb Y},\mu\right)$ which implies that $\t C_0^\infty (D_{\mathbb Y})$ is dense in $L^2(D_{\mathbb Y},\mu_{\mathbb Y})$. 
\qed 
\begin{proposition}\label{Proposition3.10}  
Let $a(0) \in {\mathbb R}^{2N \cdot d}$, $\c\in\mathcal{G}$, $\alpha \in\mathcal {A}$, such that $\h=\left(a(0);\c\, ;h(\alpha,{\mathbb Y}(\omega))\vphantom {\displaystyle {l^1}} \right)\in T^f_{{\mathbb Y}(\omega)}$. \\ 
(a) For any $F\in \t C_b^\infty(D_{\mathbb Y})$ and ${\mathbb Y}(\omega)\in D_{\mathbb Y}$ the limit 
\begin{eqnarray*} 
\partial_{\sh} F({\mathbb Y}(\omega)):=\lim_{s\to 0}\frac1s(F(\sigma_s^\alpha (\omega)) -F({\mathbb Y}(\omega))) 
\end{eqnarray*}
exists. \\ 
(b) For any $F\in\t C_b^\infty(D_{\mathbb Y})$ and ${\mathbb Y}(\omega)\in D_{\mathbb Y}$ there is some $\Nab F({\mathbb Y}(\omega))\in T^f_{{\mathbb Y}(\omega)}$ such that  
\begin{eqnarray*} 
\partial_{\sh}F({\mathbb Y}(\omega))=\langle\Nab F({\mathbb Y}(\omega)),\h\rangle_{T_{ {\mathbb Y}(\omega)}}\, . 
\end{eqnarray*}
In particular, if $F\in\t C_b^\infty(D_{\mathbb Y})$ is given by (\ref{3.1}) then 
\begin{eqnarray*} 
\Nab F({\mathbb Y}(\omega))=\left(b({\mathbb Y}(\omega))(0);\d({\mathbb Y}(\omega)), h(\beta,{\mathbb Y}(\omega))\vphantom{l^1}\right)\in T^f_{ {\mathbb Y}(\omega)}  
\end{eqnarray*}
with 
\begin{eqnarray*} 
b({\mathbb Y}(\omega))(0)=F({\mathbb Y}(\omega))\frac{\nabla f_0}{f_0}(y(\omega)(0)) \in {\mathbb R}^{2N\cdot d}\, ,  
\end{eqnarray*}
\begin{eqnarray*} 
\d({\mathbb Y}(\omega))=\left(F({\mathbb Y}(\omega))\frac{(f^\gamma_{i,j;k})'}{f^\gamma_{i,j;k}}\left(\gamma_k(i,j)\vphantom{l^1}\right)\right)_{ %\genfrac{}{}{0pt} {}{k\in {\mathbb N}}{1\le i<j\le N}}\in\mathcal{G}\, , 
{k\in {\mathbb N}}\atop {1\le i<j\le N}}\in\mathcal{G}\, ,  
\end{eqnarray*}
and 
\begin{eqnarray*} 
&&\hspace{-.5cm}\beta=\left(\left(F({\mathbb Y}(\omega))\frac{\nabla_i f_{i,j;k}} 
 {f_{i,j;k}}\left(v_i(\sigma_k(i,j))\vphantom{l^1}\right)\right)_{ %\genfrac{}{}{0pt}{}{k\in {\mathbb N} }{ 1\le i<j\le N}}\, ;\right. \\  
{k\in {\mathbb N} }\atop { 1\le i<j\le N}}\, ;\right. \\   
&&\hspace{3.5cm}\left.\left(F({\mathbb Y}(\omega))\frac{\nabla_{i'} f_{i';l}} {f_{i';l}}\left(v_{i'}(\tau_l(i'))\vphantom{l^1}\right)\right)_{ 
%\genfrac{}{}{0pt} {}{l\in {\mathbb N}}{ 1\le i'\le N}}\right)\in\mathcal{A}\, . 
{l\in {\mathbb N}}\atop{ 1\le i'\le N}}\right)\in\mathcal{A}\, .  
\end{eqnarray*}
\end{proposition} 
Proof. (b) By Theorem \ref{Theorem3.6} and ordinary differential calculus we obtain  
\begin{eqnarray*} 
&&\hspace{-.5cm}\partial_{\sh}F({\mathbb Y}(\omega))=\lim_{s\to 0}\frac1s\left( F(\sigma_s^\alpha(\omega))-F({\mathbb Y}(\omega))\vphantom{\displaystyle{l^1}} \right) \\ 
&&\hspace{.5cm}=F({\mathbb Y}(\omega))\left(\left\langle\frac{\nabla f_0}{f_0} (y(\omega)(0))\, ,\, a(0)\right\rangle\vphantom{\sum_{1\le i\le N}}\right. \\ 
&&\hspace{1cm}+\sum_{1\le i<j\le N}\ \sum_k\frac{(f^\gamma_{i,j;k})'} {f^\gamma_{i,j;k}}\left(\gamma_k(i,j)\vphantom{l^1}\right)\cdot c_k(i,j) \\ 
&&\hspace{1cm}+\sum_{1\le i<j\le N}\ \sum_k\left\langle\frac{\nabla_if_{i,j;k}} {f_{i,j;k}}\left(v_i(\cdot)\vphantom{l^1}\right)\, ,\, h_i(\alpha, {\mathbb Y} (\omega))(\cdot)\right\rangle\left(\sigma_k(i,j)\vphantom{l^1}\right) \\ 
&&\hspace{1cm}\left.+\sum_{1\le i\le N}\sum_l\left\langle\frac{\nabla_i f_{ i;l}} {f_{i;l}}\left(v_i(\cdot)\vphantom{l^1}\right)\, ,\, h_i(\alpha,{\mathbb Y} (\omega))(\cdot)\right\rangle\left(\tau_l(i)\vphantom{l^1}\right)\right) \\ 
&&\hspace{0.5cm}=:\langle\Nab F({\mathbb Y}(\omega)),\h\rangle_{T_{{\mathbb Y} (\omega)}}\vphantom{\sum^c}\, , 
\end{eqnarray*}
where, for the last line, we note that $\d({\mathbb Y}(\omega))$ and $\beta$ have only finitely many non-zero components which implies that $\Nab F({\mathbb Y}(\omega))=\left( b({\mathbb Y} (\omega))(0);\d({\mathbb Y}(\omega)), h(\beta, {\mathbb Y}(\omega))\vphantom{l^1} \right)$ belongs to $T^f_{ {\mathbb Y}(\omega)}$. 
\qed
\bigskip 

For the subsequent lemma keep in mind (\ref{3.1}) and the paragraph below it. Recall also the notation $\mathbb{H}_F$ if $F\in\t C_0^\infty(D_{\mathbb Y})$, cf. part (4) of the last block of notations. 
\begin{lemma}\label{Lemma3.11}
Using the notation $\bh_{H}:=\{\h(H,{\mathbb Y}):{\mathbb Y}\in D_{\mathbb Y}\}$, $H\in \mathbb{H}$, the sets 
\begin{eqnarray*}
\t C_0^\infty(D_{\mathbb Y};{\bf T}_{\mathbb Y}):=\left\{\sum_{j=1}^k{F}_j\cdot\bh_{H_j}: {F}_j\in\t C_0^\infty(D_{\mathbb Y}),\ H_j\in \mathbb{H}_{F_j},\ k\in {\mathbb N}\right\} 
\end{eqnarray*}
and 
\begin{eqnarray*}
\t C_b^\infty(D_{\mathbb Y};{\bf T}_{\mathbb Y}):=\left\{\sum_{j=1}^kF_j\cdot\bh_{H_j}:F_j\in\t C_b^\infty(D_{\mathbb Y}),\ H_j\in \mathbb{H}^f,\ k\in {\mathbb N}\right\}
\end{eqnarray*}
are dense in $L^2(D_{\mathbb Y},\mu_{\mathbb Y};{\bf T}_{\mathbb Y})$.
\end{lemma}
Proof. Assume that there is a $\Psi\in L^2(D_{\mathbb Y},\mu_{\mathbb Y};{\bf T}_{\mathbb Y})$ 
such that 
\begin{eqnarray*} 
\left\langle F\cdot\bh_H,\Psi\right\rangle_{L^2}=0\quad\mbox{\rm for all}\quad F\in \t C_0^\infty(D_{\mathbb Y})\quad \mbox{\rm and all}\quad \bh_H\ \mbox{\rm where}\ H\in\mathbb{H}_F.
\end{eqnarray*}
Then we have 
\begin{eqnarray*} 
E\left[ F({\mathbb Y}(\omega))\cdot\langle\h(H,{\mathbb Y}(\omega)),\Psi({\mathbb Y}(\omega)) \rangle_{T^c_{{\mathbb Y}(\omega)}}\right]=0
\end{eqnarray*}
for any fixed $H\in\mathbb{H}^f$ and all $F\in\t C_0^\infty(D_{\mathbb Y})$ such that $H\in\mathbb{H}_F$. Write $\1_{H,0}\in\t C_0^\infty(D_{\mathbb Y})$ for a function with $0\le \1_{H,0}\le 1$ and $H\in\mathbb{H}_{\1_{H,0}}$ such that $\|\1_{H,0}-1 \|_{L^2}$ is sufficiently small. 

Since any $F\in\t C_0^\infty(D_{\mathbb Y})$ can be approximated in $L^2(D_{\mathbb Y},\mu_{\mathbb Y})$ by elements $F\cdot\1_{H,0}\in \t C_0^\infty (D_{\mathbb Y})$ and it holds that $H\in\mathbb{H}_{F\cdot\1_{H,0}}$, by Lemma \ref{Lemma3.9} the last equality implies
\begin{eqnarray*} 
\langle\h(H,{\mathbb Y}(\omega)),\Psi({\mathbb Y}(\omega)) \rangle_{T^c_{{\mathbb Y}(\omega)}}=0\quad\mu_{\mathbb Y}\mbox{\rm -a.e. for all} \ H\in\mathbb{H}^f. 
\end{eqnarray*}

Note that the dense set $\mathbb{H}^f\subset \mathbb{H}$ is separable. Since $\mathbb {H}\ni H\mapsto \langle \h(H,{\mathbb Y}(\omega)),\Psi({\mathbb Y}(\omega)) \rangle_{T^c_{{\mathbb Y} (\omega)}}$ is for fixed ${\mathbb Y} (\omega)$ continuous, it is a standard conclusion that 
\begin{eqnarray*} 
\langle\h(H,{\mathbb Y}(\omega)),\Psi({\mathbb Y}(\omega)) \rangle_{T^c_{{\mathbb Y}(\omega)}}=0\quad\mbox{\rm for all} \ H\in\mathbb{H}\quad\mu_{\mathbb Y}\mbox{\rm -a.e.} 
\end{eqnarray*}
In other words, $\left\langle F\cdot\bh_H,\Psi\right\rangle_{L^2}=0$ for all $F\in \t C_0^\infty(D_{\mathbb Y})$ and $H\in\mathbb{H}_F$ yields $\Psi=0$ $\mu_{\mathbb Y}$-a.e. This implies that $\t C_0^\infty(D_{\mathbb Y};{\bf T}_{\mathbb Y})$ is dense in $L^2(D_{\mathbb Y},\mu_{\mathbb Y};{\bf T}_{\mathbb Y})$ which entails that $\t C_b^\infty (D_{\mathbb Y};{\bf T}_{\mathbb Y})$ is dense in $L^2(D_{\mathbb Y},\mu_{\mathbb Y};{\bf T}_{ \mathbb Y})$. 
\qed
\bigskip 

%\nid
{\sc Notation. (1) } For $F\in\t C_b^\infty(D_{\mathbb Y})$ denote by $c\{F\}$ the subset of the {\it coordinates} 
\begin{eqnarray*}
\left\{{\mathbb Y}(0),\left\{ \gamma_k (i,j)\vphantom {l^1}\right\}, \left\{v_i (\sigma_k (i,j))\vphantom{l^1}\right\},\left\{ v_i(\tau_l(i)) \vphantom{l^1} \right\}\right\}
\end{eqnarray*}
which are necessary to define $F$ in (\ref{3.1}). We observe that $c\{F\}$ is related to the same set of indices $(i,j;k)$ and $(i';l)$ as $\mathbb{H}_F$. 

More generally, for $F\in L^p(D_{\mathbb Y}; \mu_\mathbb {Y})$ or $F\in L^p(D_{\mathbb Y};\mu_\mathbb {Y};{\bf T}_{\mathbb {Y}})$ denote by $c\{F\}$ the set of all coordinates $\{{\mathbb Y}(0),\left\{ \gamma_k (i,j)\vphantom {l^1} \right\},\left\{ v_i(\sigma_k (i,j)) \vphantom{l^1} \right\}, \left\{ v_i(\tau_l(i) )\vphantom {l^1} \right\}\}$, except for those $F$ does not depend on, $1\le p\le \infty$. 
\medskip

\nid 
{\sc (2) }  Correspondingly, for $F\in L^p(D_{\mathbb Y}; \mu_\mathbb {Y})$ or $F\in L^p(D_{\mathbb Y};\mu_\mathbb {Y};{\bf T}_{\mathbb {Y}})$ introduce $\mathbb{H}_F$ as follows. That is the linear subspace of $\mathbb{H}$ consisting of all $(a(0);\c; \alpha)$ for which the components with the indices $(i',j'; k')$ as well as $(i'';l')$ are zero if and only if the coordinates with these indices are not included in $c\{F\}$.

\begin{proposition}\label{Proposition3.12} 
Let $F\in\t C_b^\infty(D_{\mathbb Y})$ as well as $G\in\t C_0^\infty (D_{\mathbb Y})$ and $c\{F\}\subseteq c\{G\}$ or let $F\in\t C_0^\infty(D_{\mathbb Y})$ as well as $G\in\t C_b^\infty (D_{\mathbb Y})$ and $c\{F\}=c\{G\}$. 

Furthermore, let $\h\equiv\h(H,{\mathbb Y}(\omega))\in T^f_{{\mathbb Y}(\omega)}$ for some $H=(a(0);\c\, ;\alpha)\in {\mathbb H}_G$ and all ${\mathbb Y}(\omega)\in D_{\mathbb Y}$. In addition, let $\Phi({\mathbb Y} (\omega)) := G({\mathbb Y}(\omega))\cdot\h$. Then we have 
\begin{eqnarray*} 
&&\hspace{-.5cm}E\left[\left\langle\Nab F({\mathbb Y}(\omega)),\Phi({\mathbb Y}(\omega)) \vphantom{l^1}\right\rangle_{T^c_{{\mathbb Y}(\omega)}}\right]=E\left[ \left\langle\Nab F({\mathbb Y}(\omega)),\h\vphantom{l^1}\right\rangle_{T^c_{ {\mathbb Y}(\omega)}}\cdot G({\mathbb Y}(\omega))\right]\vphantom{\dot{f}} \\ 
&&\hspace{.5cm}=E\left[\partial_{\sh}F({\mathbb Y}(\omega))\cdot G({\mathbb Y} (\omega)) \vphantom{l^1}\right]\vphantom{\dot{f}} \\ 
&&\hspace{.5cm}=-E\left[F({\mathbb Y}(\omega))\left(\partial_{\sh} G({\mathbb Y}(\omega)) +z^{\sh}({\mathbb Y}(\omega))\cdot G({\mathbb Y}(\omega))\vphantom{l^1} \right)\right]\vphantom{\dot{f}}
\end{eqnarray*}
where 
\begin{eqnarray*} 
&&\hspace{-.5cm}z^{\sh}({\mathbb Y}(\omega)):=\left\langle a(0)\, ,\nabla\ln p_0 
(y(\omega)(0))\vphantom{l^1}\right\rangle +\sum_{i<j}\ \sum_k c_k 
(i,j)(\ln g_\gamma)'(\gamma_k(i,j)) \\ 
&&\hspace{1cm}+\sum_{i<j}\ \sum_k\left\langle\vphantom{l^1}a_{i,j;k}\, ,\, \mathbf{d}\left((v_i,v_j) (\sigma_k-),\ve_k(i,j)\vphantom{l^1}\right)\right\rangle \vphantom{\dot{f}}\nonumber \\ 
&&\hspace{1cm}+\sum_{i}\sum_l\left\langle b_{i;l}\, ,\nabla_i\left[\ln \left(M(x_i (\tau_l(i)); v_i(\tau_l))\, \langle v_i(\tau_l),n_l(i)\rangle\vphantom {l^1} \right)\right] \right\rangle 
\end{eqnarray*}
and the indices over which we sum up correspond to $c\{G\}$. We have $z^{\sh}\in L^2(D_{\mathbb Y},\mu_{\mathbb Y})$. 
\end{proposition}
Proof. It follows from the assumptions on $p_0$, $g_\gamma$, $B$, and $M$ formulated in Section \ref{sec:2} that $z^{\sh}\in L^2(D_{\mathbb Y},\mu_{\mathbb Y})$. In particular recall Definitions \ref{Definition2.1} (vi), \ref{Definition2.2} (iv), and \ref{Definition2.3} (iv), calculation (\ref{2.3}), and the calculation prior to (\ref{2.4}). 

For the sake of simplicity, let $F\in\hat{C}_b^\infty(D_{\mathbb Y})$ and $G\in \hat{C}_0^\infty(D_{\mathbb Y})$ or $F\in\hat{C}_0^\infty(D_{\mathbb Y})$ and $G\in \hat{C}_b^\infty(D_{\mathbb Y})$. The more general claim follows then by linear combination. 
\medskip

Let us recall that we assume the $\sigma$-algebra $\mathcal{F}$ to be generated by the cylinder sets over the trajectories of $Y$. This implies $\mathcal{F} =\sigma\{ {\mathbb Y}(0),\left\{\gamma_k(i,j) \vphantom {l^1}\right\}, \left\{v_i (\sigma_k (i,j))\vphantom{l^1}\right\},\left\{ v_i(\tau_l(i)) \vphantom{l^1}\right \}\}$ in the notation of the proof of Lemma \ref{Lemma3.9}. Furthermore, we stress that the sums in the subsequent calculations are arranged that the set of indices over which we sum up corresponds to $c\{G\}$. 
By means of the finite dimensional partial integrations of Section \ref{sec:2}, 
(\ref{2.3}) as well as (\ref{2.4}), Theorem \ref{Theorem3.6}, and ordinary 
integral calculus we get 
\begin{eqnarray*} 
&&\hspace{-.0cm}E\left[\left\langle\Nab F({\mathbb Y}(\omega)),\Phi({\mathbb Y} (\omega)) \vphantom{l^1}\right\rangle_{T_{{\mathbb Y}(\omega)}}\right]=E\left[ \left \langle\Nab F({\mathbb Y}(\omega)),\h\vphantom{l^1}\right\rangle_{T_{ {\mathbb Y}(\omega)}}\cdot G({\mathbb Y}(\omega))\right] \\ %\vphantom{\sum_{i}} \\ 
&&\hspace{1.0cm}=E\left[F({\mathbb Y}(\omega))G({\mathbb Y}(\omega))\left( \left\langle \frac{\nabla f_0}{f_0}(y(\omega)(0))\, ,\, a(0)\right\rangle\vphantom{\sum_{ i}}\right.\right. \\ 
&&\hspace{1.5cm}+\sum_{i<j}\ \sum_k\frac{(f^\gamma_{i,j;k})'} {f^\gamma_{i,j;k}} \left(\gamma_k(i,j)\vphantom{l^1}\right)\cdot c_k(i,j) \\ 
&&\hspace{1.5cm}+\sum_{i<j}\ \sum_k\left\langle\frac{\nabla_i f_{i,j;k}}{f_{i,j;k}}\left(v_i(\sigma_k)\vphantom{l^1}\right)\, ,\, h_i(\alpha, {\mathbb Y}(\omega)) (\sigma_k)\right\rangle\nonumber \\ 
&&\hspace{1cm}\left.\left.+\sum_{i}\sum_l\left\langle\frac{\nabla_i f_{i;l}} {f_{i;l}} \left(v_i(\tau_l)\vphantom{l^1}\right)\, ,\, h_i(\alpha,{\mathbb Y} (\omega)) (\tau_l)\right\rangle\right)\right] \\ 
&&\hspace{1.0cm}=E\left[E\left[F({\mathbb Y}(\omega))G({\mathbb Y}(\omega))\left \langle \frac{\nabla f_0}{f_0}(y(\omega)(0))\, ,\, a(0)\right\rangle\left|\vphantom{\frac{\nabla f_0}{f_0}}\, {\mathbb Y}(0)\, \right]\right.\right] \\ 
&&\hspace{1.5cm}+\sum_{i<j}\ \sum_kE\left[E\left[F({\mathbb Y}(\omega)) G({\mathbb Y}(\omega))\frac{(f^\gamma_{i,j;k})'}{f^\gamma_{i,j;k}}\left(\gamma_k (i,j) \vphantom{l^1}\right)\cdot c_k(i,j)\left|\vphantom{\frac{(f^\gamma_{i,j;k})'} {f^\gamma_{i,j;k}}}\, \gamma_k(i,j)\, \right]\right. \right] \\ 
&&\hspace{1.5cm}+\sum_{i<j}\ \sum_kE\left[E\left[F({\mathbb Y}(\omega))G({\mathbb Y}(\omega))\vphantom{\frac{\nabla f_0}{f_0}}\times\right.\right.\\ 
&&\hspace{5.0cm}\left.\left.\times\left\langle\frac{\nabla_if_{i,j;k}}{f_{i,j;k}} \left( v_i(\sigma_k)\vphantom{l^1}\right)\, ,\, h_i(\alpha,{\mathbb Y}(\omega)) (\sigma_k)\right\rangle\left|\vphantom{\frac{\nabla f_0}{f_0}}\, v_i(\sigma_k)\, \right]\right.\right]\nonumber \\ 
&&\hspace{1.5cm}+\sum_{i}\sum_lE\left[E\left[F({\mathbb Y}(\omega))G ({\mathbb Y}(\omega))\vphantom{\frac{\nabla f_0}{f_0}}\times\right.\right.\\ 
&&\hspace{5.7cm}\left.\left.\times \left\langle\frac{\nabla_if_{i;l}}{f_{i;l}}\left(v_i(\tau_l) \vphantom{l^1}\right)\, ,\, h_i(\alpha,{\mathbb Y}(\omega))(\tau_l)\right\rangle \left|\vphantom{\frac{\nabla f_0}{f_0}}\, v_i(\tau_l)\, \right]\right.\right] 
\end{eqnarray*}
and hence 
\begin{eqnarray*} 
&&\hspace{-.0cm}E\left[\left\langle\Nab F({\mathbb Y}(\omega)),\Phi({\mathbb Y} (\omega)) \vphantom{l^1}\right\rangle_{T_{{\mathbb Y}(\omega)}}\right]\vphantom {\sum_{1\le i\le N}} \\ 
&&\hspace{1.0cm}=-E\left[F({\mathbb Y}(\omega))G({\mathbb Y}(\omega))\times\vphantom{ \sum_{1\le i\le N}}\right. \nonumber \\ 
&&\hspace{2.0cm}\times\left.\left(\left \langle \frac{\nabla g_0}{g_0}(y(\omega)(0))\, ,\, a(0)\right\rangle+\left \langle\nabla\ln p_0(y(\omega)(0))\, ,\, a(0)\vphantom{l^1}\right\rangle \right)\vphantom{\sum_{1\le i\le N}}\right] \\ 
&&\hspace{1.5cm}-E\left[F({\mathbb Y}(\omega))G({\mathbb Y}(\omega))\times\vphantom{ \sum_{i<j}}\right. \nonumber \\ 
&&\hspace{2.0cm}\times\left.\sum_{i<j}\ \sum_k c_k(i,j)\left(\frac{(g^\gamma_{i,j;k})'} {g^\gamma_{i,j;k}}\, (\gamma_k(i,j))+(\ln g_\gamma)'(\gamma_k(i,j))\right)\vphantom{\sum_{i<j}}\right] \\ 
&&\hspace{1.5cm}-E\left[F({\mathbb Y}(\omega))G({\mathbb Y}(\omega))\times\vphantom{ \sum_{i<j}}\right. \nonumber \\ 
&&\hspace{2.0cm}\times\left(\sum_{i<j}\ \sum_k\left\langle\frac{\nabla_i g_{i,j;k}} {g_{i,j;k}}\left(v_i(\sigma_k)\vphantom{l^1}\right)\, ,\, h_i(\alpha,{\mathbb Y}(\omega))(\sigma_k)\right\rangle\right. \\ 
&&\hspace{2.0cm}+\left.\left.\sum_{i<j}\ \sum_k\left\langle\vphantom {l^1}h_i(\alpha,{\mathbb Y}(\omega))(\sigma_k),\mathbf{d}\left((v_i,v_j) (\sigma_k-), \ve_k(i,j)\vphantom{l^1}\right)\right\rangle\vphantom{\dot{f}}\right)\right] \nonumber \\ 
&&\hspace{1.5cm}-E\left[F({\mathbb Y}(\omega))G({\mathbb Y}(\omega))\times\vphantom{ \sum_{i<j}}\right. \nonumber \\ 
&&\hspace{2.0cm}\left.\times\left(\sum_{i}\sum_l\left\langle\frac {\nabla_ig_{i;l}} {g_{i;l}}\left(v_i(\tau_l)\vphantom{l^1}\right)\, ,\, h_i(\alpha,{\mathbb Y}(\omega)) (\tau_l(i))\right\rangle\right.\right. \\ 
&&\hspace{2.0cm}\left.\left.+\sum_{i}\sum_l\left\langle h_i(\alpha,{\mathbb Y} (\omega))(\tau_l)\, , \nabla_i\left[\ln\left(M(x_i(\tau_l);v_i(\tau_l))\, \langle v_i(\tau_l) ,n_l(i)\rangle\vphantom{l^1}\right)\right]\right\rangle \right)\right]\, . 
\end{eqnarray*}
The claim is now a consequence of the definition of $h(\alpha,{\mathbb Y}(\omega))$ 
in Definition \ref{Definition3.2} (c). 
\qed 
\bigskip

%\nid
{\sc Notation. } For $\Phi= \sum_{j=1}^kG_j\cdot\bh_j\in \t C_b^\infty(D_{\mathbb Y};{\bf T}_{\mathbb Y})$ with $G_j\in\t C_b^\infty (D_{\mathbb Y})$ and $\bh_j\equiv \{\h_j(H_j,{\mathbb Y}): {\mathbb Y}\in D_{\mathbb Y}\}$, $H_j\in \mathbb{H}_{G_j}$, $j\in\{1,\ldots ,k\}$, $k\in {\mathbb N}$, set 
\begin{eqnarray}\label{8}
\delta\Phi({\mathbb Y}(\omega)):=-\sum_{j=1}^k\left(\partial_{\sh_j} G_j({\mathbb Y}
(\omega))+z^{\sh_j}({\mathbb Y}(\omega))\cdot G_j({\mathbb Y}(\omega))\right)
\end{eqnarray}
where $z^{\sh_j}$ is given in Proposition \ref{Proposition3.12}. Furthermore, let us recall the definition of $\Nab F$, $F\in\t C_b^\infty (D_{\mathbb Y})$, from Proposition \ref{Proposition3.10}. 
\begin{proposition}\label{Proposition3.13}
The operator $\left(\Nab,\t C_b^\infty(D_{\mathbb Y})\mapsto L^2(D_{\mathbb Y},\mu_{\mathbb Y};{\bf T}_{\mathbb Y})\right)$ is closable on the space $L^2(D_{ \mathbb Y},\mu_{\mathbb Y})$ and the operator $\left(\delta, \t C_0^\infty (D_{\mathbb Y};{\bf T}_{\mathbb Y})\mapsto L^2 (D_{\mathbb Y},\mu_{\mathbb Y})\right)$ is closable on $L^2(D_{\mathbb Y},\mu_{\mathbb Y};{\bf T}_{ \mathbb Y})$. 
\end{proposition}
Proof. {\it Step 1 } Let $F=\sum_{r=1}^n\hat{F}_r\in\t C_b^\infty(D_{\mathbb Y})$, where $\hat{F}_r\in \hat{C}_b^\infty(D_{\mathbb Y})$, and 
\begin{eqnarray*}
\Phi=\sum_{j'=1}^{k'} {\Phi}_{j'}\cdot\bh_{H_{j'}}\in\t C_0^\infty (D_{\mathbb Y};{\bf T}_{\mathbb Y}). \quad\mbox{\rm Note}\ \ {\Phi}_{j'}\in\t C_0^\infty (D_{\mathbb Y}),\ H_{j'} \in \mathbb{H}_{{\Phi}_{j'}},\ {k'}\in {\mathbb N} \, , 
\end{eqnarray*}
by definition of $\t C_0^\infty (D_{\mathbb Y};{\bf T}_{\mathbb Y})$. Recalling the Notation (3) below (\ref{3.1}), we observe that the indices for which the factors in (\ref{3.1}) of the $\widehat{(\Phi_{j'})}_r$ are not constant (equal to one) correspond to $c \{\Phi_{j'}\}$. Therefore $E\left[\left. F ({\mathbb Y} (\omega)) \right| c\{\Phi_{j'}\} \right]$ is obtained from $F$ by keeping the factors in (\ref{3.1}) of the $\hat{F}_r$ which correspond to $c \{\Phi_{j'} \}$ and replacing the others with the resulting constants. Thus, 
\begin{eqnarray}\label{3.3}
&&\hspace{-.5cm}E\left[\left\langle\Nab F({\mathbb Y}(\omega)),\Phi({\mathbb Y}(\omega)) \vphantom{l^1} \right\rangle_{T_{{\mathbb Y}(\omega)}}\right]=\sum_{j'=1}^{k'}E\left[\left\langle\Nab F({\mathbb Y}(\omega)), \bh_{H_{j'}}\vphantom{l^1} \right\rangle_{T_{{\mathbb Y}(\omega)}}\Phi_{j'}({\mathbb Y}(\omega))\right]\nonumber \\ 
&&\hspace{.5cm}=\sum_{j'=1}^{k'}E\left[\left\langle\Nab E\left[\left. F ({\mathbb Y}(\omega)) \right|c\{\Phi_{j'}\}\ \right], \bh_{H_{j'}}\vphantom{l^1} \right\rangle_{T_{{\mathbb Y}(\omega)}}\Phi_{j'}({\mathbb Y}(\omega))\right]\, . 
\end{eqnarray}
By virtue of Proposition \ref{Proposition3.12}, it follows now that 
\begin{eqnarray*}
&&\hspace{-.5cm}E\left[\left\langle\Nab F({\mathbb Y}(\omega)),\Phi({\mathbb Y}(\omega)) \vphantom{l^1} \right\rangle_{T_{{\mathbb Y}(\omega)}}\right] \\ &&\hspace{.5cm}=\sum_{j'=1}^{k'}E\left[E\left[\left. F({\mathbb Y}(\omega)) \right| c\{\Phi_{j'}\} \right]\cdot\delta\left( \Phi_{j'}({\mathbb Y}(\omega)) \bh_{H_{j'}}\right)\vphantom{l^1}\right]\, .
\end{eqnarray*}
Now assume that, instead of $F$ we have a sequence $F_n\in \t C_b^\infty(D_{\mathbb Y})$, $n\in {\mathbb N}$, with $F_n\stack{n\to\infty}{\lra}0$ in $L^2(D_{\mathbb Y},\mu_{\mathbb Y})$ and $\Nab F_n\stack{n\to\infty}{\lra}f$ in $L^2 (D_{\mathbb Y},\mu_{\mathbb Y};{\bf T}_{\mathbb Y})$. Then 
\begin{eqnarray*}
E\left[\left\langle f,\Phi({\mathbb Y}(\omega))\vphantom{l^1}\right\rangle_{T_{{\mathbb Y} (\omega)}}\right]=\lim_{n\to\infty}\sum_{j'=1}^{k'}E\left[E\left[\left. F_n({\mathbb Y}(\omega)) \right| c\{\Phi_{j'}\} \right]\cdot\delta\left( \Phi_{j'}({\mathbb Y}(\omega)) \bh_{H_{j'}}\right)\vphantom{l^1}\right]=0. 
\end{eqnarray*}
Since $\Phi$ is chosen arbitrarily from $\t C_0^\infty(D_{\mathbb Y};{\bf T}_{\mathbb Y})$ Lemma \ref{Lemma3.11} implies $f=0$, which means that  $(\Nab,\t C_b^\infty(D_{\mathbb Y}))$ is closable on $L^2(D_{\mathbb Y},\mu_{\mathbb Y})$. 
\medskip

\nid
{\it Step 2 } Let $F\in\t C_b^\infty(D_{\mathbb Y})$ and $\Phi_n\in\t C_0^\infty (D_{\mathbb Y}; {\bf T}_{\mathbb Y})$, $n\in {\mathbb N}$, be a sequence with 
\begin{eqnarray*}
\Phi_n\equiv\sum_{j'=1}^{{k_n}'} (\Phi_n)_{j'}\cdot\bh_{H^n_{j'}} \stack {n\to\infty} {\lra}0\quad\mbox{\rm in }\quad L^2(D_{\mathbb Y},\mu_{\mathbb Y};{\bf T}_{\mathbb Y})
\end{eqnarray*}
and $\delta\Phi_n\stack {n\to\infty}{\lra}\vp$ in $L^2(D_{\mathbb Y},\mu_{\mathbb Y})$. By the arguments of part (a) we observe 
\begin{eqnarray*}
\left\|\Nab E\left[\left. F({\mathbb Y}(\omega))\right|c\{(\Phi_n)_{j'}\} \right]\right\|_{L^2}\le \left\|\Nab F({\mathbb Y}(\omega))\right\|_{L^2}\, ,\quad j'=1,\ldots ,{k_n}', 
\end{eqnarray*}
and from Proposition \ref{Proposition3.12} it follows that 
\begin{eqnarray*}
&&\hspace{-.5cm}E\left[\left\langle\Nab E\left[\left. F({\mathbb Y}(\omega))\right| c\{(\Phi_n)_{j'}\} \right] ,(\Phi_n)_{j'}({\mathbb Y}(\omega))\cdot \bh_{H^n_{j'}}\vphantom{l^1} \right\rangle_{T_{{\mathbb Y}(\omega)}}\right] \\ 
&&\hspace{.5cm}=E\left[E\left[\left. F({\mathbb Y}(\omega)) \right| c\{(\Phi_n)_{j'}\} \right]\cdot \delta\left( (\Phi_n)_{j'}({\mathbb Y}(\omega))\cdot\bh_{H^n_{j'}}\right) \vphantom{l^1}\right]\, ,\quad j'=1,\ldots ,{k_n}'. 
\end{eqnarray*}
Since $\delta((\Phi_n)_{j'}\cdot\bh_{H^n_{j'}})$ is $c\{(\Phi_n)_{j'}\}$-measurable by Proposition \ref{Proposition3.12} and (\ref{8}), the last two relations imply 
\begin{eqnarray*}
E[F({\mathbb Y}(\omega))\cdot\vp]=0  
\end{eqnarray*}
for all $F\in\t C_b^\infty(D_{\mathbb Y})$. Lemma \ref{Lemma3.9} yields $\vp=0$. Therefore, $(\delta,\t C_0^\infty(D_{\mathbb Y};{\bf T}_{\mathbb Y}))$ is closalbe on $L^2(D_{\mathbb Y},\mu_{\mathbb Y};{\bf T}_{\mathbb Y})$.
\qed 
\begin{definition}\label{Definition3.14}{\rm
(a) The closure $(\Nab,D(\Nab))$ of 
\smallskip

\centerline{$\left(\Nab,\t C_b^\infty(D_{\mathbb Y})\mapsto L^2(D_{\mathbb Y},\mu_{\mathbb Y};{\bf T}_{\mathbb Y})\right)$ on $L^2(D_{\mathbb Y},\mu_{\mathbb Y})$} 

\nid
is called {\it gradient}. Furthermore, we denote the closure of 
\smallskip

\centerline{$\left(\Nab,\t C_0^\infty(D_{\mathbb Y})\mapsto L^2(D_{\mathbb Y},\mu_{\mathbb Y};{\bf T}_{\mathbb Y})\right)$ on $L^2(D_{\mathbb Y},\mu_{\mathbb Y})$} 

\nid
by $(\Nab,D_0 (\Nab))$. \\ 
(b) The closure $(-\delta,D_0(\delta))$ of 
\smallskip

\centerline{
$\left(-\delta,\t C_0^\infty(D_{\mathbb Y}; {\bf T}_{\mathbb Y})\mapsto L^2(D_{\mathbb Y},\mu_{\mathbb Y})\right)$ on $L^2(D_{\mathbb Y},\mu_{\mathbb Y};{\bf T}_{\mathbb Y})$ } 

\nid
is called {\it divergence}. 
}
\end{definition}
\begin{corollary}\label{Corollary3.15}
(a) (Duality) Let $F\in D(\Nab)$ and $\Phi\in D_0(\delta)$ such that $c\{F\} \subseteq c\{\Phi\}$. We have 
\begin{eqnarray*}
E\left[\left\langle\Nab F({\mathbb Y}(\omega)),\Phi({\mathbb Y}(\omega))\vphantom{l^1} \right\rangle_{T_{{\mathbb Y}(\omega)}}\right]=E\left[F({\mathbb Y}(\omega)) \cdot \delta \Phi({\mathbb Y}(\omega))\right]\, . 
\end{eqnarray*}
(b) Let $F\in D(\Nab)\cap L^\infty(D_{\mathbb Y},\mu_{\mathbb Y})$ and $G\in D_0(\Nab)$ such that $c\{F\}\subseteq c\{G\}$. Furthermore, let $\h\equiv\h(H, {\mathbb Y}(\omega))\in T_{{\mathbb Y}(\omega)}$ for some $H\in {\mathbb H}_G$ and all ${\mathbb Y}(\omega) \in D_{\mathbb Y}$. Then we have 
\begin{eqnarray*} 
&&\hspace{-.5cm}E\left[\left\langle\Nab F({\mathbb Y}(\omega)),G({\mathbb Y}(\omega)) \cdot\h)\vphantom{l^1}\right\rangle_{T_{{\mathbb Y}(\omega)}}\right] \\ 
&&\hspace{.5cm}=-E\left[F({\mathbb Y}(\omega))\left(\left\langle\Nab G({\mathbb Y}(\omega)),\h\vphantom{l^1}\right\rangle_{T_{{\mathbb Y}(\omega)}}+z^{\sh}({\mathbb Y}(\omega))\cdot G({\mathbb Y}(\omega))\vphantom{l^1}\right)\right]\vphantom{\dot{f}}
\end{eqnarray*}
where $z^{\sh}$ is given in Proposition \ref{Proposition3.12}. \\ 
(c) (Chain Rule) Let $\ve>0$ and $(-\ve,\ve)\ni u\mapsto {\mathbb X}_{s+u}\in D_{\mathbb 
Y}$ such that ${\mathbb X}_\cdot$ is differentiable in $s\in {\mathbb R}$ and let 
$F\in D(\Nab)$. Then, using the notation of Definition \ref{Definition3.5} (b), 
\begin{eqnarray*}
\frac{d}{ds}F({\mathbb X}_s)=\left\langle\Nab F\, ,\, \left(\dot{{\mathbb X}}_s, \dot {\g}_s\right)\right\rangle_{T_{{\mathbb X}_s}}\equiv\left\langle\Nab F\, ,\, \left((x_s(0) ,h_s(0));\dot{\g}_s\, ;h(\alpha_s,{\mathbb X}_s)\vphantom{\displaystyle {l^1}}\right) \right\rangle_{T_{{\mathbb X}_s}}\ a.e. 
\end{eqnarray*}
\end{corollary}
Proof. Part (a) is an immediate consequence of Proposition \ref{Proposition3.12} and Proposition \ref{Proposition3.13}. Let us turn to part (b). Because of $G\in D_0(\Nab)$ there exists $G_n\in\t C_0^\infty (D_{\mathbb Y})$, $n\in\mathbb{N}$, such that $G_n\stack{n\to \infty}{\lra}G$ in $L^2(D_{\mathbb Y},\mu_{\mathbb Y})$ as well as $\Nab G_n\stack{n\to \infty}{\lra}\Nab G$ in $L^2(D_{\mathbb Y},\mu_{\mathbb Y})$. Let $H_n$ be the orthogonal projection of $H\in\mathbb{H}$ to $\mathbb{H}_{G_n}$ and set $\h_n\equiv\h(H_n, {\mathbb Y} (\omega))$. By Proposition \ref{Proposition3.12} and a calculation similar to (\ref{3.3}) it holds that 
\begin{eqnarray*} 
&&\hspace{-.5cm}E\left[\left\langle\Nab F({\mathbb Y}(\omega)),G_n({\mathbb Y}(\omega)) \cdot\h_n)\vphantom{l^1}\right\rangle_{T_{{\mathbb Y}(\omega)}}\right] \\ 
&&\hspace{.5cm}=E\left[\left\langle\Nab E\left[\left.F({\mathbb Y}(\omega))\right|c\{G_n\} \right],G_n({\mathbb Y}(\omega)) \cdot\h_n)\vphantom{l^1}\right\rangle_{T_{{\mathbb Y}(\omega)}}\right] \\ 
&&\hspace{.5cm}=-E\left[E\left[\left.F({\mathbb Y}(\omega))\right|c\{G_n\} \right] \left( \left\langle\Nab G_n({\mathbb Y}(\omega)),\h_n\vphantom{l^1}\right\rangle_{T_{{\mathbb Y}(\omega)}}+z^{\sh_n}({\mathbb Y}(\omega))\cdot G_n({\mathbb Y}(\omega))\vphantom{l^1}\right)\right]\vphantom{\dot{f}}\, .
\end{eqnarray*}

Recall the first paragraph of the proof of Proposition \ref{Proposition3.12}. Since $E\left[\left.F\right|c\{G_n\} \right]$ converges (on some subsequence of indices $n$) boundedly a.e. to $F$ and $z^{\sh_n}({\mathbb Y} (\omega))\cdot G_n({\mathbb Y}(\omega))$ converges in $L^1(D_{\mathbb Y},\mu_{ \mathbb Y})$ to $z^{\sh}({\mathbb Y}(\omega))\cdot G({\mathbb Y}(\omega))$ as $n\to \infty$, it holds that  
\begin{eqnarray*}
&&\hspace{-.5cm}E\left[E\left[\left.F({\mathbb Y}(\omega))\right|c\{G_n\} \right]\cdot z^{\sh_n}({\mathbb Y}(\omega))\cdot G_n({\mathbb Y}(\omega ))\right] \\ 
&&\hspace{-.5cm}\stack{n\to\infty} {\lra}E\left[F({\mathbb Y}(\omega))\cdot z^{\sh}({\mathbb Y}(\omega))\cdot G({\mathbb Y}(\omega)) \right]\, .\vphantom{\dot{f}} 
\end{eqnarray*}
We get (b). 

Concerning part (c), for $F\in \t C_b^\infty(D_{\mathbb{Y}})$ the claim is a consequence of Definition \ref{Definition3.5} (b) and ordinary differential calculus. For $F\in D(\Nab)$ it follows then from Proposition \ref{Proposition3.13}. 
\qed 
\bigskip

%\nid
{\sc Notation. }For $1\le p\le \infty$ and $F\in L^2(D_{\mathbb Y}, \mu_{\mathbb Y})$ introduce $L_F^p (D_{\mathbb Y},\mu_{\mathbb Y}):=\{G\in L^p(D_{\mathbb Y},\mu_{\mathbb Y}):c\{G\}  =c\{F\}\} $. 
\begin{proposition}\label{Proposition3.16} 
Let $\h\equiv\h(H,{\mathbb Y}(\omega))\in T_{{\mathbb Y}(\omega)}$ for some $H\in {\mathbb H}$ and all ${\mathbb Y}(\omega)\in D_{\mathbb Y}$ and 
\begin{eqnarray*}
&&\hspace{-.5cm}D(\partial_{\sh}):=\left\{\vphantom{\dot{f}}F\in L^2(D_{\mathbb Y}, \mu_{\mathbb Y}):H\in \mathbb{H}_F\ \mbox{\rm and there exists } \mathbf{C}_F>0\ \mbox{\rm such that for all }\right. \\  
&&\hspace{3.3cm}\left. G\in D_0(\Nab)\cap L_F^\infty(D_{\mathbb Y},\mu_{\mathbb Y}) \ \mbox{\rm we have }\right. \vphantom{\left(\dot{f}\right)}\\  
&&\hspace{1.3cm}\left.\left|\vphantom{\dot{f}}E\left[FG\cdot z^{\sh}\vphantom{l^1}\right]+E\left[F\cdot\langle\Nab G,\h\rangle_{{T}_{{\mathbb Y}(\omega)}}\vphantom{l^1}\right]\right|\le \mathbf{C}_F\, \|G\|_{L^2(D_{\mathbb Y},\mu_{\mathbb Y})}\right\}\, .  
\end{eqnarray*}
(a) For $F\in D(\partial_{\sh})$ there exists a unique bounded linear functional $\mathcal{L}_F$ on $L_F^2(D_{\mathbb Y},\mu_{\mathbb Y})$ such that for all $G\in D_0(\Nab)\cap L_F^\infty(D_{\mathbb Y},\mu_{\mathbb Y})$ we have 
\begin{eqnarray*}
&&\hspace{-.5cm}-E\left[FG\cdot z^{\sh}\vphantom{l^1}\right]-E\left[F\cdot \langle \Nab G,\h\rangle_{{T}_{{\mathbb Y}(\omega)}}\vphantom{l^1}\right]= \mathcal{L}_F(G)\, .  
\end{eqnarray*}
In particular there exists a representing element $\t \partial_{\sh}F$ of $\mathcal {L}_F$, that is the unique element $\t \partial_{\sh}F\in L_F^2(D_{\mathbb Y}, \mu_{\mathbb Y})$ satisfying 
\begin{eqnarray*}
\langle\t G,\t \partial_{\sh}F\rangle_{L^2(D_{\mathbb Y}, \mu_{\mathbb Y})}=\mathcal{L}_F(\t G)\quad\mbox{\rm for all }\quad\t G\in L_F^2(D_{\mathbb Y},\mu_{\mathbb Y})\, . 
\end{eqnarray*}
(b) We have $D_{\sh}(\Nab):=\{F\in D(\Nab):H\in\mathbb{H}_F\}\subseteq D(\partial_{\sh})$ and 
\begin{eqnarray*}
\t \partial_{\sh}F=\langle\Nab F,\h\rangle_{{T}_{{\mathbb Y}(\omega)}}\, ,\quad F\in D_{\sh}(\Nab)\quad \mu_{\mathbb Y}\mbox{\rm -a.e.}
\end{eqnarray*}
(c) $\left(\t \partial_{\sh},D(\partial_{\sh})\right)$ is a closed operator extending $\left(\partial_{\sh}\, ,\{F\in\t C_b^\infty(D_{\mathbb Y}):H\in\mathbb {H}_F\}\right)$ of Proposition \ref{Proposition3.10} (a). 
\end{proposition} 
Proof. Noting that the set of all $G\in D_0(\Nab)\cap L_F^\infty(D_{\mathbb Y},\mu_{ \mathbb Y})$ is dense in $L^2_F(D_{\mathbb Y}, \mu_{\mathbb Y})$ by the  arguments of the proof of Lemma \ref{Lemma3.9}, part (a) follows immediately from  the Riesz representation theorem. 
\medskip 

\nid 
Next we focus on part (b). For $F\in D_{\sh}(\Nab)$ we obtain from Corollary \ref{Corollary3.15} (b)
\begin{eqnarray*}
&&\hspace{-.5cm}\left|E\left[FG\cdot z^{\sh}\vphantom{l^1}\right]+E\left[F\cdot \langle\Nab G,\h\rangle_{{T}_{{\mathbb Y}(\omega)}}\vphantom{l^1}\right]\vphantom{\dot{f}}\right|=\left|\left\langle G,\langle\Nab F,\h\rangle_{{T}_{{\mathbb Y}(\omega)}}\right\rangle_{L^2(D_{\mathbb Y},\mu_{\mathbb Y})}\right| \\ 
&&\hspace{.5cm}\le\left\|\langle\Nab F,\h\rangle_{{T}_{{\mathbb Y}(\omega)}}\right\|_{ L^2(D_{\mathbb Y},\mu_{\mathbb Y})}\, \|G\|_{L^2(D_{\mathbb Y},\mu_{\mathbb Y})} 
\end{eqnarray*}
for all $G\in D_0(\Nab)\cap L_F^\infty(D_{\mathbb Y},\mu_{\mathbb Y})$. We get $D_{\sh} (\Nab)\subseteq D(\partial_{ \sh})$. Furthermore, 
\begin{eqnarray*}
&&\hspace{-.5cm}\langle G,\t \partial_{\sh}F\rangle_{L^2(D_{\mathbb Y},\mu_{\mathbb Y})} =-E\left[FG\cdot z^{\sh}\vphantom{l^1}\right]-E\left[F\cdot\langle\Nab G,\h \rangle_{{T}_{{\mathbb Y}(\omega)}}\vphantom{l^1}\right] \\ 
&&\hspace{.5cm}=\left\langle G,\langle\Nab F,\h\rangle_{{T}_{{\mathbb Y}(\omega)}}\right\rangle_{L^2(D_{\mathbb Y},\mu_{\mathbb Y})} 
\end{eqnarray*}
for $F\in D_{\sh}(\Nab)$ and $G\in D_0(\Nab)\cap L_F^\infty(D_{\mathbb Y},\mu_{ \mathbb Y})$ follows from part (a) and again from Corollary \ref{Corollary3.15} (b). Again the fact that the set of all $G\in D_0(\Nab)\cap L_F^\infty (D_{\mathbb Y},\mu_{ \mathbb Y})$ is dense in $L^2_F(D_{\mathbb Y}, \mu_{\mathbb Y})$ implies now 
\begin{eqnarray*}
\t \partial_{\sh}F=\langle\Nab F,\h\rangle_{{T}_{{\mathbb Y}(\omega)}}\, ,\quad F\in D_{\sh}(\Nab)\quad \mu_{\mathbb Y}\mbox{\rm -a.e.} 
\end{eqnarray*}
Let us now prove part (c). Let $\h\equiv\h(H,{\mathbb Y}(\omega))\in T_{{\mathbb Y}(\omega)}$ for some $H\in {\mathbb H}$ and all ${\mathbb Y}(\omega)\in D_{\mathbb Y}$. Assume 
\begin{eqnarray*}
D(\partial_{\sh})\ni F_n\stack{n\to\infty}{\lra}F\quad{\rm in}\ L^2(D_{\mathbb Y}, 
\mu_{\mathbb Y}) 
\end{eqnarray*}
and 
\begin{eqnarray*}
\t \partial_{\sh}F_n\stack{n\to\infty}{\lra}\vp\quad{\rm in}\ L^2(D_{\mathbb Y}, 
\mu_{\mathbb Y})  
\end{eqnarray*}
for some $\vp\in L^2(D_{\mathbb Y},\mu_{\mathbb Y})$. Without loss of generality, by these assumptions we may suppose $c\{F_n\}=c\{F\}$ and $\mathbb{H}_{F_n} =\mathbb{H}_F$. We have to demonstrate that $F\in D(\partial_{\sh})$ and $\vp=\t \partial_{\sh}F$. From part (a) we obtain for all $G\in D_0(\Nab)\cap L^\infty_F (D_{\mathbb Y},\mu_{\mathbb Y})$ 
\begin{eqnarray*}
&&\hspace{-.5cm}E\left[\t \partial_{\sh}F_n\cdot G\vphantom{l^1}\right]=-E\left[ F_nG\cdot z^{\sh}\vphantom{l^1}\right]-E\left[F_n\cdot\langle\Nab G,\h\rangle_{{T}_{{\mathbb Y}(\omega)}}\vphantom{l^1}\right] \\ 
&&\hspace{-.0cm}\stack{n\to\infty}{\lra}-E\left[FG\cdot z^{\sh}\vphantom{l^1} \right]-E\left[F\cdot\langle\Nab G,\h\rangle_{{T}_{{\mathbb Y}(\omega)}}\vphantom{l^1} \right]\, .  
\end{eqnarray*}
Thus, for all $G\in D_0(\Nab)\cap L^\infty_F (D_{\mathbb Y},\mu_{\mathbb Y})$ it holds that 
\begin{eqnarray*}
E\left[\vp\cdot G\vphantom{l^1}\right]=-E\left[FG\cdot z^{\sh}\vphantom{l^1} \right]-E\left[F\cdot\langle\Nab G,\h\rangle_{{T}_{{\mathbb Y}(\omega)}}\vphantom{l^1} \right] \, , 
\end{eqnarray*}
which by part (a) entails that $F\in D(\partial_{\sh})$ and $\vp=\t \partial_{\sh} F$. In other words, $(\t \partial_{\sh},D(\Nab))$ is a closed operator extending $\left(\partial_{\sh}\, ,\{F\in\t C_b^\infty(D_{\mathbb Y}): H\in\mathbb{H}_F\} \right)$ of Proposition \ref{Proposition3.10}. 
\qed 
\medskip 

%\nid
{\sc Notation. } Proposition \ref{Proposition3.16} (c) gives rise to denote 
\begin{eqnarray*}
\left(\partial_{\sh},D(\partial_{\sh})\right)\equiv \left(\t \partial_{\sh}, 
D(\partial_{\sh})\right)\, . 
\end{eqnarray*}
\begin{definition}\label{Definition3.17}{\rm
(a) Using the notation $\bh_{H}:=\{\h(H,{\mathbb Y}):{\mathbb Y}\in D_{\mathbb Y}\}$, $H\in \mathbb{H}$, the operator $(\partial_{\bh_{H}},D(\partial_{\bh_{H}}))$ is called {\it directional derivative} in direction of $\bh_{H}\in {\bf T}_{\mathbb Y}$. 
\\ 
(b) In particular, $\partial_{\sh}F({\mathbb Y}(\omega))$ is called {\it directional derivative} in direction of $\h\equiv\h(H,{\mathbb Y})\in T_{{\mathbb Y}(\omega)}$ {\it at} ${\mathbb Y}(\omega) \in D_{\mathbb Y}$. For $F\in\t C_b^\infty(D_{\mathbb Y})$ such that $H\in \mathbb{H}_F$ it is defined everywhere on $D_{\mathbb Y}$ and for $F\in D(\partial_{\sh})$ such that $H\in \mathbb{H}_F$ it is defined $\mu_{\mathbb Y}$-a.e. on $D_{\mathbb Y}$. 
}
\end{definition}

{\sc Remark. } Recall the arguments of Step 1 of the proof of Proposition \ref{Proposition3.13}. For $G\in D_0(\Nab)\cap L^\infty (D_{ \mathbb Y}, \mu_{\mathbb Y})$ there is an approximating sequence $\t G_n\in \t C_0^\infty( D_\mathbb {Y})$, $n\in \mathbb{N}$, i.e. $\t G_n \stack{n\to \infty}{\lra} G$ and $\Nab \t G_n \stack{n\to \infty} {\lra} \Nab G$ in $L^2 (D_{\mathbb Y}, \mu_{ \mathbb Y})$. The functions $G_n:=E[G|c\{\t G_n\}]$ satisfy $G_n\in L^\infty (D_{ \mathbb Y}, \mu_{ \mathbb Y})$, $n\in \mathbb {N}$. Moreover, by orthogonal projection of $G$ to $L^2_{\t G_n} (D_{ \mathbb Y}, \mu_{\mathbb Y})$, we have $G_n\in D_0(\Nab)$, $n\in \mathbb {N}$, and $G_n \stack{n\to \infty}{\lra} G$ as well as $\Nab G_n\stack{n\to \infty} {\lra} \Nab G$ in $L^2 (D_{\mathbb Y}, \mu_{ \mathbb Y})$. 
\bigskip

Keeping the last Remark in mind, the proofs of the following two assertions are based on standard conclusions of the closedness of $(\Nab,D_0(\Nab))$ as well as $(\delta,D_0(\delta))$ and are therefore omitted. First we reestablish formula (\ref{8}) relative to the operator $(-\delta,D_0(\delta))$ of Definition \ref{Definition3.14} (b). 
\begin{corollary}\label{Corollary3.18} 
Let $k\in {\mathbb N}$ and $\bh_j\equiv\bh_{H_j}:=\{\h(H_j,{\mathbb Y}):{\mathbb Y}\in D_{\mathbb Y}\}$, $H_j\in \mathbb{H}$, $j\in\{1,\ldots ,k\}$. Then for all $j\in\{1,\ldots ,k\}$
\begin{itemize} 
\item[{}] $G_j\in D_0(\Nab)\cap L^\infty(D_{\mathbb Y},\mu_{\mathbb Y})$ and $H_j\in \mathbb{H}_{G_j}$ implies $G_j\cdot \bh_j\in D_0(\delta)$.
\end{itemize} 
For all $j$ let $z^{ \sh_j}$ be given by Proposition \ref{Proposition3.12}. For $\Phi=\sum_{j=1}^kG_j\cdot\bh_j$, 
\begin{eqnarray*}
\delta\Phi=-\sum_{j=1}^k\left(\partial_{\sh_j} G_j+z^{\sh_j}\cdot G_j\right)\, . 
\end{eqnarray*}
\end{corollary}

The next proposition contains two known rules from infinite dimensional stochastic 
calculus adjusted to the situation of the present paper.  
\begin{proposition}\label{Proposition3.19}
(a) (Chain rule for the gradient) Let $u=(u_1,\ldots ,u_N)$ be a random variable 
with values in $V^N$ and let $u_i=(u_{(i,1)},\ldots u_{(i,d)})$, $i\in\{1,\ldots 
,N\}$. Suppose $u_{(i,r)}\in D(\Nab)$ for all $i\in\{1,\ldots ,N\}$ and $r\in\{1, 
\ldots ,d\}$. Then for all $\phi\in C^1_b(V^N)$ we have $\phi(u)\in D(\Nab)$ and 
it holds that 
\begin{eqnarray*}
\Nab\phi(u)=\sum_{i=1}^N\sum_{r=1}^d\partial_{(i,r)}\phi(u)\Nab u_{(i,r)}\, .  
\end{eqnarray*}
(b) (Product rule for the divergence) Let $\Phi\in D_0(\delta)$ and $G\in D(\Nab)\cap L^\infty(D_{\mathbb Y},\mu_{\mathbb Y})$. Then we have $\Phi\cdot G\in D_0(\delta)$ and 
\begin{eqnarray*}
\delta(\Phi\cdot G)=G\cdot\delta(\Phi)-\langle\Nab G\, ,\Phi\rangle_{{T}_{\mathbb Y}} 
\, . 
\end{eqnarray*}
\end{proposition} 

\section{Regularity of densities}\label{sec:4}
\setcounter{equation}{0}

In this section suppose that the initial value $\mathbb{Y}(0)=(x(0),$ $v(0))$ possesses a smooth positive density $p_0$ with respect to the Lebesgue measure on $(D^N)_\beta\times V^N$. In addition, assume that the compactly on $(0,1)$ supported probability density $g_\gamma$ is smooth, cf. Definition \ref{Definition2.2}, that the redistribution kernel $M$ is smooth on $\partial D \times V$, cf. Definition \ref{Definition2.3}, and that the collision kernel $B$ is positive and smooth on $\{(v,v',e):v,v'\in V,\, e\in S^{d-1}_+ (v-v')$ such that $(v^\ast, v'^\ast) \in V\times V$ and $|\langle v-v',e\rangle|\ge \kappa\}$, cf. Definition \ref{Definition2.2} (iii) and the paragraph before Definition \ref{Definition2.2}. 

Fix $t>0$ for the whole section. Recalling $|v|<v_{max}<\infty$ observe that the number of collisions during the period $[0,t]$ is a.e. uniformly bounded by some $I_c\in\mathbb{N}$ because of Definition \ref{Definition2.2} (b) (ii)-(iv) and hypothesis $p\, $(ii) on the position of the particles posted below Definition \ref{Definition2.1}. Similarly, the number of reflections during $[0,t]$ is a.e. uniformly bounded by some $I_r\in \mathbb {N}$ by Definition \ref{Definition2.3} (ii) and (iii). 
\bigskip

%\nid
{\sc Notation. (1) } Let $\mathcal{I}$ denote the set of all sequences $I=(\rho_1, \ldots ,\rho_{I_c+I_r})$ of $I_c+I_r$ random times of collisions or reflections taken from $\{\sigma_k(i,j),\, \tau_l (i'):\ k\in\{1,\ldots ,I_c\},\ l\in\{1,\ldots ,I_r\},\ 1\le i<j\le N,\ 1\le i'\le N\}$ such that 
\begin{itemize}
\item[(i)] if $\sigma_k(i,j)=\rho_m$ for some $1<k\le I_c+I_r$ and $1<m\le I_c+I_r$ then $\sigma_{k-1}(i,j)=\rho_n$ for some $1\le n<m$ and 
\item[(ii)] if $\tau_l (i')=\rho_m$ for some $1<l\le I_c+I_r$ and $1<m\le I_c+I_r$ then $\tau_{l-1}(i')=\rho_n$ for some $1\le n<m$\, .
\end{itemize} 
Note here, that for defining the elements $I\in \mathcal{I}$ we do not consider the concrete outcomes of these random times. Recall also that by Definition \ref{Definition3.2} (a) and the construction of Section \ref{sec:2} the corresponding sequence of random times $\rho_1(\mathbb{Y}), \ldots ,\rho_{I_c+I_r} (\mathbb{Y})$ is a.e. strictly increasing where, with $\rho_0 =0$, we have $\rho_0 <t<\rho_{I_c +I_r+1} (\mathbb{Y})$.
\bigskip

%\nid 
{\sc (2) } Consider the sequence $I=(\rho_1, \ldots ,\rho_{I_c+I_r})\in\mathcal{I}$, and define the set 
\begin{eqnarray*}
D_{{\mathbb{Y}},I}:=\left\{\mathbb{Y}\in D_\mathbb{Y}: (\rho_1 (\mathbb{Y}), \ldots ,\rho_{I_c+I_r} (\mathbb{Y}))=I\right\}\, .
\end{eqnarray*}
\smallskip

%\nid
{\sc (3) } In order to prepare Lemmas \ref{Lemma4.4} and \ref{Lemma4.5} below, for $I\in \mathcal{I}$, let us introduce 
\begin{eqnarray*} 
J_c\equiv J_c(I)\in\{0,1\}^{(N(N-1)/2)\cdot (I_c+I_r)}\quad\mbox{\rm and}\quad J_r\equiv J_r(I)\in\{0,1\}^{N\cdot (I_c+I_r)} 
\end{eqnarray*}
in the following manner. For a concrete $I=(\rho_1 (\mathbb{Y}), \ldots ,\rho_{I_c +I_r} (\mathbb {Y}))\in \mathcal{I}$, set 
\begin{eqnarray*} 
&&\hspace{-.5cm}(J_c)_{i,j;k}=1\quad\begin{array}{l}\mbox{\rm if there is } m\in \{1,\ldots , I_c+I_r\}\ \mbox{\rm such that }\\ \rho_m=\sigma_k(i,j) \ \mbox{\rm where }k\in\{1,\ldots, I_c+I_r\},\ 1\le i<j\le N, \end{array}
\end{eqnarray*}
otherwise set $(J_c)_{i,j;k}=0$. Similarly put 
\begin{eqnarray*} 
&&\hspace{-.5cm}(J_r)_{i;l}=1\quad\begin{array}{l}\mbox{\rm if there is } m\in \{1,\ldots , I_c+I_r\}\ \mbox{\rm such that }\\ \rho_m=\tau_l(i)\ \mbox{\rm where }l\in\{1,\ldots, I_c+I_r\},\ 1\le i \le N, \end{array}
\end{eqnarray*}
otherwise put $(J_r)_{i;l}=0$. Let 
\begin{eqnarray*} 
|J_c|:=\sum_{k\in\{1,\ldots, I_c+I_r\},\ 1\le i<j\le N}(J_c)_{i,j;k}\quad\mbox{\rm and}\quad |J_r|:=\sum_{l\in\{1,\ldots, I_c+I_r\},\ 1\le i \le N}(J_r)_{i;l}\, .
\end{eqnarray*}
\smallskip

%\nid
{\sc (4) } In this sense we shall introduce the random vectors 
\begin{eqnarray*} 
&&\hspace{-.5cm}\mathcal{S}\equiv (\sigma_k(i,j))_{1\le i<j\le N;\ k\in\{1,\ldots, I_c+I_r\}} \quad\mbox{\rm and}\quad\mathcal{T}\equiv (\tau_l (i'))_{1\le i'\le N;\ l\in\{1,\ldots, I_c+I_r\}} 
\end{eqnarray*}
and write $(\mathcal{S},{J_c})$ and $(\mathcal{T},{J_r})$ for the corresponding labeled vectors. Observe that the times labeled with one are precisely the first $I_c+I_r=|J_c|+|J_r|$ times of collision or reflection. %In this way identify $I$ 
Introduce also the random vectors 
\begin{eqnarray*} 
&&\hspace{-.5cm}\t \Gamma\equiv (\gamma_k(i,j))_{1\le i<j\le N;\ k\in\{1,\ldots, I_c+I_r\}} 
\end{eqnarray*}
and 
\begin{eqnarray*} 
&&\hspace{-.5cm}\mathcal{V}_\mathcal{S}\equiv (v_i(\sigma_k(i,j)))_{1\le i<j\le N;\ k \in\{1,\ldots, I_c+I_r\}} \ \mbox{\rm and} \ \mathcal{V}_\mathcal{T}\equiv (v_{i'}(\tau_l (i')))_{1\le i'\le N;\ l\in\{1,\ldots, I_c+I_r\}} \, .
\end{eqnarray*}
\smallskip

%\nid
{\sc (5) } For $I 
% \sim \left\{(\mathcal{S},{J_c})\, ,\, (\mathcal{T},{J_r}) \right\} 
\in \mathcal {I}$ write as above $(\t \Gamma,{J_c})\equiv (\t \Gamma,{J_c}(I))$ for the vector obtained from $\t \Gamma$ by labeling all those components of $\t \Gamma$ with one or zero, for which the corresponding component of $J_c$ is one or zero. In this way, $(\t \Gamma,{J_c})$ is associated with the times in $(\mathcal{S},{J_c})$, where the label one appears for the collisions within the first $I_c+I_r=|J_c|+|J_r|$ collisions or reflections. Note that this is compatible with Definitions \ref{Definition2.2} and \ref{Definition3.3}. Furthermore write $(\mathcal{V}_\mathcal {S},{J_c})$ and $(\mathcal{V}_\mathcal{T}, {J_r})$ for the labeled vectors obtained accordingly from $\mathcal{V}_\mathcal {S}$ and $\mathcal {V}_\mathcal{T}$. 
\medskip

Consider also the non-labeled vectors $\hat\Gamma$, $\hat{V}_\mathcal {S}$, and $\hat {V}_\mathcal{T}$ introduced as follows. Whenever a component of $(\t \Gamma,{J_c})$, $(\mathcal{V}_\mathcal {S},{J_c})$, or $(\mathcal{V}_\mathcal {T},{J_r})$ is labeled with one then the corresponding component of $\hat\Gamma$, $\hat{V}_\mathcal {S}$, or $\hat {V}_\mathcal{T}$ is the outcome relative to the choice of $\mathbb{Y} \in D_{\mathbb {Y},I}$ of the respective component of $\t \Gamma$, $\mathcal {V}_\mathcal {S}$, or $\mathcal{V}_\mathcal {T}$. Furthermore, whenever a component of $(\t \Gamma,{J_c})$, $(\mathcal{V}_\mathcal {S},{J_c})$, or $(\mathcal {V}_\mathcal {T},{J_r})$ is labeled with zero then the corresponding component of $\hat\Gamma$, $\hat{V}_\mathcal {S}$, or $\hat{V}_\mathcal{T}$ is the random variable of the respective component of $\t \Gamma$, $\mathcal{V}_\mathcal {S}$, or $\mathcal{V}_\mathcal {T}$. 
\bigskip

%\nid
{\sc (6) } 
The coordinates of $\mathbb{Y}\in D_{\mathbb{Y},I}$ are in the sense of Definition \ref{Definition3.3} of the form 
\begin{eqnarray*} 
&&\hspace{-.5cm}\left(\mathbb{Y}(0);\hat\Gamma ;\hat{V}_\mathcal {S} ,\hat{V}_\mathcal {T} \right)\, \
\end{eqnarray*}
where we remark that all coordinates from time zero up to the $(I_c+I_r)$-th collision or reflection are non-random and all other coordinates are random. Let $\Gamma ,\,  {V}_\mathcal {S} ,\, {V}_\mathcal {T}$ be obtained from $\hat\Gamma ,\, \hat{V }_\mathcal {S} ,\, \hat{V}_\mathcal {T}$ by keeping the non-random components and replacing the random components with $(0,1)$, $V$, and $V$, respectively.
\bigskip

%\nid
{\sc (7) } Introduce
\begin{eqnarray*}
&&\hspace{-.5cm}W_{\mathbb{Y},I}:=\left\{\mathbb{Y}(0)\times\Gamma\times{V}_\mathcal {S}\times {V}_\mathcal {T} :\mathbb{Y}\in D_{\mathbb{Y},I}\right\}\, ,\quad I\in \mathcal{I} .
\end{eqnarray*} 
Noting that, for fixed $I\in \mathcal{I}$, for all $\mathbb{Y}\in D_{\mathbb {Y}, I}$ the vectors $(\t \Gamma,{J_c})$, $(\mathcal{V}_\mathcal {S},{J_c})$, and $(\mathcal {V}_\mathcal {T},{J_r})$ are equally labeled, the inclusion 
\begin{eqnarray*}
W_{\mathbb{Y},I}\subseteq (D^N)_\beta\times V^N\times(0,1)^{\frac12 N(N-1)\cdot (I_c+I_r)}\times V^{\frac12 N(N-1)\cdot (I_c+I_r)}\times V^{N\cdot (I_c+I_r)} 
\end{eqnarray*} 
is meaningful. Recall the following. The identification of $\omega\in\Omega$ with the trajectory $Y(\omega)$, i.e. the map $\Omega\ni\omega\mapsto {\mathbb Y} (\omega) \in D_{\mathbb Y}$ is a bijection. Correspondingly we have introduced $\mu_{\mathbb Y}$ to be the image measure of $P$ under this map. 
\begin{lemma}\label{Lemma4.4} 
The sets $W_{\mathbb {Y},I}$, $I\in\mathcal{I}$, are mutually disjoint open subsets of $(D^N)_\beta\times V^N\times(0,1)^{\frac12 N(N-1)\cdot (I_c+I_r)}\times V^{\frac12 N(N-1)\cdot (I_c+I_r)}\times V^{N\cdot (I_c+I_r)}$ such that for the corresponding sets $D_{\mathbb {Y},I}$, $I\in\mathcal {I}$, of reduced trajectories it holds that 
\begin{eqnarray*}
\left\{ \mathbb{Y}(0)\times\Gamma\times{V}_\mathcal {S}\times {V}_\mathcal {T} :\mathbb{Y}\in D_{\mathbb{Y}} \right\}\subseteq\bigcup_{I\in\mathcal{I}}\overline {W_{\mathbb {Y},I}} 
\end{eqnarray*}
and 
\begin{eqnarray*} 
\mu_{\mathbb {Y}}\left(\bigcup_{I\in\mathcal{I}}D_{\mathbb{Y},I}\right)=1\, . 
\end{eqnarray*}
\end{lemma}
Proof. For $I\in \mathcal {I}$ and $\mathbb{Y} \in D_{\mathbb {Y},I}$ shorten the vectors $(\t \Gamma,{J_c})\equiv (\t \Gamma,{J_c}(I))$, $(\mathcal{V}_\mathcal {S},J_c)\equiv (\mathcal{V}_\mathcal {S},J_c(I))$, and $(\mathcal {V}_\mathcal {T},J_c)\equiv (\mathcal {V}_\mathcal {T},J_c(I))$ by just including the concrete outcomes, with respect to $\mathbb {Y} \in D_{\mathbb {Y},I}$, of the components labeled with one. Denote the shortened vectors by $(\Gamma,\1)$, $({V}_\mathcal {S},\1)$, and $({V}_\mathcal{T},\1)$ and define 
\begin{eqnarray*}
&&\hspace{-.5cm}W_{\mathbb{Y},I;\1}:=\left\{\left(\, \mathbb{Y}(0)\, ;\, (\Gamma, \1)\, ;\, ({V}_\mathcal {S},\1)\, ,\, ({V}_\mathcal{T},\1)\vphantom{l^1}\right) :\mathbb{Y}\in D_{\mathbb{Y},I}\right\}\, .
\end{eqnarray*} 

Let us demonstrate that every element of $W_{\mathbb {Y},I;\1}$ is an inner point of $W_{\mathbb{Y},I;\1}$ with respect to the topology in $(D^N)_\beta \times V^N \times (0,1)^{|J_c|}\times V^{|J_c|}\times V^{|J_r|}$. Keeping Lemma \ref{Lemma3.4} in mind, the rest is then obvious. 
\medskip

\nid
{\it Step 1 } Fix $I\in\mathcal{I}$, let $\mathbb{Y}\in D_{\mathbb{Y},I}$, and consider the outcome 
\begin{eqnarray*} 
\left(\, \mathbb{Y}(0)\, ;\, (\Gamma, \1)\, ;\, ({V}_\mathcal {S},\1)\, ,\, ({V}_\mathcal{T},\1)\vphantom{l^1}\right)\in W_{\mathbb {Y},I;\1} 
\end{eqnarray*}
with respect to $\mathbb{Y}$. Consider also open balls $B_{r}({\mathbb {Y}})$ with respect to the maximum norm in $(D^N)_\beta\times V^N\times (0,1)^{|J_c|}\times V^{|J_c|} \times V^{|J_r|}$ with center $(\, \mathbb{Y}(0)\, ;\, (\Gamma, \1)\, ;\, ({V}_\mathcal {S},\1)\, ,\, ({V}_\mathcal{T},\1))$ and radius $r$. Choose $r_0$ such that all points belonging to $B_{r_0} ({\mathbb {Y}})$  correspond by the construction of Definition \ref{Definition3.3} to reduced trajectories $\t {\mathbb{Y}}$ for which we have  
\begin{itemize} 
\item[(i)]
\begin{eqnarray*} 
\rho_1(\t{\mathbb {Y}})=\left\{\begin{array}{rl} 
\sigma_1(i,j)(\t{\mathbb {Y}})\quad & \mbox{\rm if }\quad \rho_1({\mathbb {Y}})=\sigma_1(i,j)({\mathbb {Y}}) \\ 
\tau_1(i')(\t{\mathbb {Y}})\quad & \mbox{\rm if }\quad \rho_1({\mathbb {Y}})=\tau_1(i')({\mathbb {Y}}) 
\end{array}\right.\, , \quad \t{\mathbb {Y}}\in B_{r_0}({\mathbb {Y}}). 
\end{eqnarray*}
\end{itemize} 
For the possibility of the choice of such an $r_0$ keep in mind hypothesis $p$(iii) of Section \ref{sec:2} and the convexity of the domain $D$ as well as the smoothness of the boundary $\partial D$. Recall Definitions \ref{Definition2.2} as well as \ref{Definition2.3}. As a consequence of (i), $\rho_1(\t{\mathbb {Y}})$ depends continuously on the coordinates $\left(\, \mathbb{Y}(0)\, ;\, (\Gamma, \1) \, ;\, ({V}_\mathcal {S},\1)\, ,\, ({V}_\mathcal{T},\1)\vphantom{l^1}\right)$ of $\t {\mathbb {Y}}\in B_{r_0} ({\mathbb {Y}})$ if $\rho_1({\mathbb {Y}})=\tau_1(i')$. If $\rho_1({\mathbb {Y}}) =\sigma_1(i,j) ({\mathbb {Y}})$ then the entrance time $s\equiv s(\t {\mathbb {Y}})$ of the first collision depends continuously on the coordinates $\left(\, \mathbb{Y}(0)\, ;\, (\Gamma, \1) \, ;\, ({V}_\mathcal {S},\1)\, ,\, ({V}_\mathcal{T},\1)\vphantom{l^1}\right)$ of $\t {\mathbb {Y}}\in B_{r_0} ({\mathbb {Y}})$. 

In the next three steps let us demonstrate that the subsequent statement (ii) is a consequence of (i). 
\medskip 

\begin{itemize} 
\item[(ii)] 
\begin{eqnarray*} 
\left\{\left(x(\rho_1(\t{\mathbb {Y}})-),v(\rho_1(\t{\mathbb {Y}})-), \gamma_1(i,j) (\t{\mathbb {Y}})\right): \t{\mathbb {Y}}\in B_{r_0}({\mathbb {Y}})\right\}
\end{eqnarray*}
is an open subsets of $D^N \times V^N\times (0,1)$ if $\rho_1(\t{\mathbb {Y}})=\sigma_1(i,j) (\t{\mathbb {Y}})$ and 
\begin{eqnarray*} 
\left\{\left(x(\rho_1(\t{\mathbb {Y}})-),v(\rho_1(\t{\mathbb {Y}})-)\right): \t{\mathbb {Y}}\in B_{r_0}({\mathbb {Y}})\right\}
\end{eqnarray*}
is an open subset of $D^{i'-1}\times \partial D\times D^{N-i'} \times V^N$ if $\rho_1(\t{\mathbb {Y}})=\tau_1(i') (\t{\mathbb {Y}})$. 
\end{itemize} 

\nid
{\it Step 2 } In case that $\rho_1 =\tau_1 (i')$ for some $1\le i'\le N$ we first observe that for fixed $v(0)=v_0$ the set 
\begin{eqnarray*}
&&\hspace{-.5cm}\left\{(x(\rho_{1} (\t{\mathbb {Y}})-), v(\rho_1 (\t{\mathbb {Y}})-)):\t{\mathbb {Y}} \in B_{r_0}({\mathbb {Y}}),\, v(0)=v_0\right\} \\ 
&&\hspace{1cm}=\left\{(x(\rho_1(\t{\mathbb {Y}})),v(0)):\t{\mathbb {Y}} \in B_{r_0}({\mathbb {Y}})\, ,v(0)=v_0\right\}
\end{eqnarray*}
coincides with $\mathcal {D}(v(0))\times\{v(0)\}$ for some open set $\mathcal {D}(v(0))\subset D^{i'-1}\times \partial D\times D^{N-i'}$. If $v(0)$ is no longer fixed but belongs to $\{v(0)(\t {\mathbb {Y}}):\t {\mathbb {Y}}\in B_{r_0}({\mathbb {Y}})\}$, it will be demonstrated in Step 4 below that the set $\{(x(\rho_{1} (\t{\mathbb {Y}})-), v(\rho_1 (\t{\mathbb {Y}})-)):\t{\mathbb {Y}} \in B_{r_0}({\mathbb {Y}})\}$, is an open subset of $D^{i'-1}\times \partial D\times D^{N-i'} \times V^N$. 
\medskip

\nid
{\it Step 3 } Furthermore, in case that  $\rho_1 =\sigma_1(i,j)$ for some $1\le i<j\le N$ we first claim that for fixed $v(0)=v_0$ the set 
\begin{eqnarray*}
&&\hspace{-.5cm}\left\{(x(\rho_1 (\t {\mathbb {Y}})-),\gamma_1 (i,j)(\t {\mathbb {Y}}),v (\rho_1 (\t{\mathbb {Y}})-)): \t{\mathbb {Y}} \in B_{r_0}({\mathbb {Y}}),\, v(0)=v_0\right\} \\ 
&&\hspace{1cm}=\left\{(x (\rho_1(\t {\mathbb {Y}})),\gamma_1 (i,j)(\t {\mathbb {Y}}), v(0)):\t{\mathbb {Y}} \in B_{r_0} ({\mathbb {Y}})\, ,v(0)=v_0\right\}
\end{eqnarray*}
coincides with $\mathcal {D}(v(0))\times\{v(0)\}$ for some open set $\mathcal{D} (v(0))\subset D^{N}\times (0,1)$. 

To see this, fix for a moment also the position $x(0)$ at time 0, let $s\equiv s(\mathbb{Y})$ denote the first entrance time of the first collision, and let $t>s$ be the associated time in the sense of Definition \ref{Definition2.2} (a) (i). Then by Definition \ref{Definition2.2} (b) (ii),(iv) the following holds. The position 
\begin{eqnarray*}
x(\sigma_1 (i,j)) = x(0)+\left(s +\frac12 \gamma\cdot (t-s)\right)v(0)\, ,\quad\gamma\equiv\gamma_1(i,j)\in (\ve_\gamma, 1-\ve_\gamma), 
\end{eqnarray*} 
can be attained at time $\rho_1 = \sigma_1 (i,j)$ as long as for the chosen $\gamma\in (\ve_\gamma, 1-\ve_\gamma)$, there is a ${\mathbb {Y}}'\in B_{r_0} ({\mathbb {Y}})$ which coincides with $\mathbb {Y}$ on $[0,s]$ such that the outcome of $\gamma_1 (i,j)({\mathbb {Y}}')$ is $\gamma$.
As a consequence, for fixed $v(0)$ and $x(0)$, there exists an open interval $(l,u)\equiv \left( l(x(0),v(0)), u(x(0), v(0)) \vphantom{\int}\right)$ such that 
\begin{eqnarray*}
x(\sigma_1 (i,j)) = x(0)+(s + a)v(0)\, ,\quad a\in (l,u),
\end{eqnarray*} 
are the possible positions of such ${\mathbb {Y}}'$ at time $\rho_1({\mathbb {Y}}') = \sigma_1(i,j)({\mathbb {Y}}')$. Since for fixed $v(0)$, the bounds $l(x(0),v(0))$ and $u(x(0), v(0))$ depend continuously on $x(0)$ on $\{x(0)(\t {\mathbb {Y}}):\t {\mathbb {Y}}\in B_{r_0}({\mathbb {Y}})\}$, cf. hypothesis $p$(iii) of Section 2 and (i) of this proof, the last sentence yields the statement from the beginning of Step 3. 
\medskip

\nid
{\it Step 4 } The set $\mathcal {V}_0:=\{v(0)\equiv v(0)(\t {\mathbb {Y}}):\t {\mathbb {Y}}\in B_{r_0}({\mathbb {Y}})\}$, is an open subset of $V^N$. For $v(0)\in\mathcal {V}_0$, let $\mathcal{D} (v(0))$ be the set constructed in Steps 2 and 3 in the respective cases. 

By the construction of Steps 2 and 3 for any $w(0)\in\mathcal {V}_0$ and $z\in  \mathcal {D}(w(0))$ there is a neighborhood $N(z) \subset \mathcal {D}(w(0))$ of $z$ and a neighborhood $N(w(0))\subset\mathcal{V}_0$ of $w(0)$ with $N(z)\subset \mathcal {D}(v(0))$ for all $v(0)\in N(w(0))$. Thus, given $w(0)\in \mathcal{V}_0$ and $z\in\mathcal{D} (w(0))$,  
\begin{eqnarray*}
N(z)\times N(w(0))\ \mbox{is a neighborhood of } (z,w(0)) \, , 
\end{eqnarray*}
being a subset of $\bigcup_{v(0)\in{\mathcal {V}}_0}\left(\mathcal{D}(v(0))\times \{v(0)\}\right)$. Denoting 
\begin{eqnarray*}
z(\rho_1)\equiv \left(x(\sigma_1(i,j) (\t{\mathbb {Y}}))\, ,\, \gamma_1(i,j) (\t{\mathbb {Y}}) \right) 
\end{eqnarray*}
if $\rho_1=\sigma_1(i,j)$ and $z(\rho_1) \equiv x(\tau_1(i') (\t{\mathbb {Y}}))$ if $\rho_1=\tau_1(i')$ we obtain that 
\begin{eqnarray*} 
&&\hspace{-.5cm}\left\{\left(z(\rho_1(\t{\mathbb {Y}})),v(\rho_1(\t{\mathbb {Y}} )-) \right): \t{\mathbb {Y}}\in B_{r_0}({\mathbb {Y}})\right\} = \left\{\left( z(\rho_1(\t{\mathbb {Y}})),v(0)(\t {\mathbb {Y}}) \right): \t{\mathbb {Y}}\in B_{r_0}({\mathbb {Y}})\right\} \\
&&\hspace{.5cm}=\bigcup_{v(0)\in{\mathcal {V}}_0}\mathcal{D}(v(0))\times \{v(0)\}=\bigcup_{v(0)\in \mathcal {V}_0,z\in\mathcal{D}(v(0))}N(z)\times 
N(v(0)) 
\end{eqnarray*}
is an open subset of $D^{N}\times (0,1)\times V^N$ if $\rho_1(\t{\mathbb {Y}})= \sigma_1(i,j) (\t{\mathbb {Y}})$ or an open subset of $D^{i'-1}\times \partial D\times D^{N-i'} \times V^N$ if $\rho_1(\t{\mathbb {Y}})=\tau_1(i') (\t{\mathbb {Y}})$. In other words we have verified (ii). 
\medskip

\nid
{\it Step 5 } Now repeat these arguments starting out at time $\rho_1(\t{\mathbb {Y}})$ from $\{( x(\rho_1(\t{\mathbb {Y}})),v(\rho_1(\t{\mathbb {Y}})) ): \t{\mathbb 
{Y}}\in B_{r_{1}}({\mathbb {Y}})\}$ with some $r_1\in (0,r_0]$ instead of $r_0$ such that 
\begin{itemize} 
\item[(i1)]
\begin{eqnarray*} 
\rho_2(\t{\mathbb {Y}})=\left\{\begin{array}{rl} 
\sigma_k(i,j)(\t{\mathbb {Y}})\quad & \mbox{\rm if }\quad \rho_2({\mathbb {Y}})=\sigma_k(i,j)({\mathbb {Y}}) \\ 
\tau_l(i')(\t{\mathbb {Y}})\quad & \mbox{\rm if }\quad \rho_2({\mathbb {Y}})=\tau_l(i')({\mathbb {Y}}) 
\end{array}\right.\, , \quad k,l\in \{1,2\}\, ,\quad \t{\mathbb {Y}}\in B_{r_1}({\mathbb {Y}}). 
\end{eqnarray*}
\end{itemize} 
Note that for repeating the arguments of Steps 2-4 we need (ii). 

Continue at time $\rho_2( \t{\mathbb {Y}})$ from $\{( x(\rho_2(\t{\mathbb {Y}})), v(\rho_2(\t{\mathbb {Y}})) ): \t{\mathbb {Y}}\in B_{r_{2}} ({\mathbb {Y}})\}$ with some suitable $r_2\in (0,r_1]$, and so forth to verify that for all $n\in\{1,\ldots ,I_c +I_r+1\}$
\begin{eqnarray*} 
\left\{\left(z(\rho_{n}(\t{\mathbb {Y}})),v(\rho_n(\t{\mathbb {Y}})-) \right): \t{\mathbb {Y}}\in B_{r_{n-1}}({\mathbb {Y}})\right\}
\end{eqnarray*}
is an open subset of $D^{N}\times (0,1)\times V^N$ if $\rho_{n}(\t{\mathbb {Y}}) =\sigma_k (i,j) (\t{\mathbb {Y}})$ for some $k\in\{1,\ldots ,n\}$ or an open subset of $D^{i'-1}\times \partial D\times D^{N-i'} \times V^N$ if $\rho_{n} (\t {\mathbb {Y}})=\tau_l(i') (\t{\mathbb {Y}})$ for some $l\in \{1, \ldots , n\}$. Here, $r_{n-1}$ and $z(\rho_n)$ are introduced in the same way as $r_1$ and $z(\rho_1)$ above. 

The last paragraph implies that for all $n\in\{1,\ldots ,I_c +I_r+1\}$ 
\begin{eqnarray*} 
\left\{\left(\gamma_k(\sigma_k (i,j) (\t{\mathbb {Y}})),v(\sigma_k (i,j) (\t{\mathbb {Y}})-) \right): \t{\mathbb {Y}}\in B_{r_{n-1}}({\mathbb {Y}})\right\}
\end{eqnarray*}
is an open subset of $(0,1)\times V^N$ if $\rho_{n}(\t{\mathbb {Y}}) =\sigma_k (i,j) (\t{\mathbb {Y}})$ for some $k\in\{1,\ldots ,n\}$. Furthermore, 
\begin{eqnarray*} 
\left\{v(\tau_l(i')(\t{\mathbb {Y}})-):\t{\mathbb {Y}}\in B_{r_{n-1}}({\mathbb {Y}})\right\}
\end{eqnarray*}
is an open subset of $V^N$ if $\rho_{n}(\t{\mathbb {Y}}) = \tau_l (i') (\t{\mathbb {Y}})$ for some $l\in\{1,\ldots ,n\}$. Therefore, any element of $W_{\mathbb {Y},I; \1}$ is an inner point of $W_{\mathbb{Y},I;\1}$ with respect to the topology in $(D^N)_\beta\times V^N\times (0,1)^{|J_c|}\times V^{|J_c|}\times V^{|J_r|}$. \qed 
\medskip

%\nid
{\sc Notation. } For fixed $I\in\mathcal{I}$ define 
\begin{eqnarray*}
D_{\mathbb{Y},I,m}:=\left\{\mathbb{Y}\in D_{\mathbb{Y},I}: \rho_m<t<\rho_{m+1}\right\}\, ,\quad m\in\{0,1,\ldots ,I_c+I_r\}, 
\end{eqnarray*}
and, for fixed $I\in\mathcal{I}$ as well as $m\in\{0,1,\ldots ,I_c+I_r\}$, 
\begin{eqnarray*}
W_{\mathbb{Y},I,m}:=\left\{\mathbb{Y}(0)\times\Gamma\times{V}_\mathcal {S}\times {V}_\mathcal {T} :\mathbb{Y}\in D_{\mathbb{Y},I,m}\right\}\, .
\end{eqnarray*}
\begin{lemma}\label{Lemma4.5} 
For fixed $I\in\mathcal{I}$, the sets $W_{\mathbb {Y},I,m}$, $m\in \{0,1,\ldots , I_c +I_r\}$, are mutually disjoint open subsets of $(D^N)_\beta\times V^N\times (0,1)^{ \frac12 N(N-1)\cdot (I_c+I_r)}\times V^{\frac12 N(N-1)\cdot (I_c+I_r)} \times V^{N\cdot (I_c+I_r)}$. For the sets $W_{\mathbb {Y},I,m}$ and the corresponding sets $D_{\mathbb{Y},I,m}$, $m\in \{0,1,\ldots , I_c +I_r\}$, of reduced trajectories it holds that 
\begin{eqnarray*}
W_{\mathbb {Y},I}\subset\bigcup_{m\in\{1,\ldots , I_c+I_r\}}\overline{W_{\mathbb {Y},I,m}} 
\end{eqnarray*}
and 
\begin{eqnarray*} 
\mu_{\mathbb {Y}}\left(D_{\mathbb {Y},I}\setminus \bigcup_{m\in \{0,1,\ldots ,I_c+I_r\}}D_{\mathbb {Y},I,m}\right)=0\, . 
\end{eqnarray*}
\end{lemma}
Proof. Keeping in mind that we require existence of an initial density with respect to the Lebesgue measure, the second sentence is a standard consequence of the construction of Boltzmann type processes, cf. Sections \ref{sec:1} as well as \ref{sec:2}, and in particular Definitions \ref{Definition2.1}-\ref{Definition2.3}. 

Let $I\in\mathcal{I}\, $. To verify the first claim choose ${\mathbb {Y}} \in D_{\mathbb {Y},I,m}$ for some $m\in \{0,1,\ldots ,I_c +I_r\}$. Let $Y\equiv ((x_{1} (\cdot),v_{1}(\cdot)), \ldots ,(x_{N}(\cdot), v_{N} (\cdot)))$ be recovered from ${\mathbb {Y}}\equiv (x(0),v(0);v_{1} (\cdot),\ldots ,v_{N} (\cdot))$ by $x_{i} (\cdot) =x_{i}(0)+ \int_0^\cdot v_{i}(s) \, ds$, $1\le i\le N$. Then, because of $v\neq 0$ if $v\in V$, 
\begin{eqnarray*} 
|x(t)-x(\rho_m)|\equiv|x(t)(\mathbb{Y})-x(\rho_m)(\mathbb{Y})|>0
\end{eqnarray*}
and $v(t')-v(\rho_m)\equiv v(t')(\mathbb{Y})-v(\rho_m)(\mathbb{Y})=0$ for $t'\in [\rho_m ({\mathbb{Y}}),t]$. Taking over the notation $B_r (\mathbb{Y})$ from the proof of the previous lemma, there exists a sufficiently small $r>0$ such that for all $\t {\mathbb {Y}}\in B_r (\mathbb{Y})$ we have $t<\rho_{m+1}(\t {\mathbb {Y}})$, 
\begin{eqnarray*} 
\left|x(t)\left(\t{\mathbb{Y}}\right)-x(\rho_m)\left(\t{\mathbb{Y}}\right)\right| >\frac12|x(t)(\mathbb{Y})-x(\rho_m)(\mathbb{Y})|>0 
\end{eqnarray*}
as well as $v(t')(\t{\mathbb{Y}})=v(\rho_m)(\t{\mathbb{Y}})$ for $t'\in [\rho_m (\t {\mathbb{Y}}),t]$, and therefore ${\t{\mathbb {Y}}}\in D_{\mathbb {Y},I,m}$. By (the proof of) Lemma \ref{Lemma4.4}, $(\, \mathbb{Y}(0)\, ;\, (\Gamma, \1)\, ;\, ({V}_\mathcal {S},\1)\, ,\, ({V}_\mathcal {T},\1)) $ is therefore an inner point of $W_{\mathbb {Y},I;\1}$ with respect to the topology in $(D^N)_\beta\times V^N\times (0,1)^{|J_c|}\times V^{|J_c|}\times V^{|J_r|}$. This implies the lemma. \qed
\bigskip

%\nid 
{\sc Notation. (1) } Introduce 
\begin{eqnarray*}
\t D_{\mathbb{Y}}:=\bigcup_{I\in\mathcal{I},\, m\in\{0,1,\ldots ,I_c+I_r\}}D_{\mathbb{Y},I,m} \, .
\end{eqnarray*} 
\smallskip

\nid
{\sc (2) } Let $\t C^\infty_0(\t D_{\mathbb{Y}})$ the set of all 
\begin{eqnarray*}
&&\hspace{-0.5cm}F(\mathbb{Y}(\omega))\\ 
&&\hspace{0.5cm}=f\left( y(\omega)(0), \gamma_k(i,j), v(\sigma_k(i,j)), v(\tau_l(i'));\ 1\le i<j\le N, 1\le i'\le N,\ k,l\le I_c+I_r \vphantom{\dot{f}} \right)
\end{eqnarray*} 
such that we have the following. 
\begin{itemize}
\item[(i)] 
\begin{eqnarray*}
f\in C^\infty_b\left((D^N)_\beta\times V^N\times(0,1)^{ \frac12 N(N-1)\cdot (I_c+I_r)} \times V^{\frac12 N(N-1)\cdot (I_c+I_r)}\times V^{N\cdot (I_c+I_r)} \right)\, . 
\end{eqnarray*} 
\item[(ii)] If $n$ is the index of a coordinate corresponding to some $\tau_l(i')$ then we suppose that 
\begin{eqnarray*}
f(\ldots ,v_n,\ldots )=f(\ldots ,-v_n,\ldots )\, ,\quad v\in V. 
\end{eqnarray*} 
\item[(iii)] The function $f$ is supported by a compact subset of $\bigcup_{I \in\mathcal{I}\, m\in \{1,\ldots,I_c+I_r\}} W_{\mathbb{Y},I,m} $. 
\end{itemize}
\begin{proposition}\label{Proposition4.2} 
Let $t>0$ and let $G\in\t C_0^\infty(\t D_{\mathbb{Y}})$ as well as $\phi\in C_b^\infty({\mathbb R}^{N\cdot d})$. \\ 
(a) For all $i\in\{1,\ldots ,N\}$ and $r\in\{1,\ldots ,d\}$ the random variable $\Nab v_{(i,r)}(t)$ is well-defined by Proposition \ref{Proposition3.10} and we have 
\begin{eqnarray*}
E\left[\partial_{(i,r)}\phi(v(t))\, G\right]=E\left[\phi(v(t))\, H_{(i,r)}(v(t),G)\right]
\end{eqnarray*}
where 
\begin{eqnarray}\label{4.3}
&&\hspace{-.5cm}H_{(i,r)}(v(t),G)=\delta \left(G\cdot\Nab v_{(i,r)}(t)\right)\, .
\end{eqnarray}
In particular, $H_{(i,r)}(v(t),G)\in \t C_0^\infty(\t D_{\mathbb{Y}})$. \\ 
(b) For all $i\in\{1,\ldots ,N\}$ and $r\in\{1,\ldots ,d\}$ we have 
\begin{eqnarray*}
E\left[\partial_{(i,r)}\phi(x(t))\, G\right]=E\left[\phi(x(t))\, H_{(i,r)}(x(t),G)\right]
\end{eqnarray*}
where 
\begin{eqnarray*}
&&\hspace{-.5cm}H_{(i,r)}(x(t),G)=\sum_{I\in\mathcal{I},\, m\in \{0,1,\ldots ,I_c+I_r\}}\delta \left(\frac{\chi_{D_{\mathbb{Y},I,m}}\cdot G}{t-\rho_m}\cdot \Nab v_{(i,r)}(t)\right)\, .
\end{eqnarray*}
In particular, $H_{(i,r)}(x(t),G)\in \t C_0^\infty(\t D_{\mathbb{Y}})$. 
\end{proposition}

Proof. {\it Step 1 } The first two steps are devoted to (a). Recall that from a given reduced trajectory $\mathbb{Y}(\omega)$ the (full) trajectory 
\begin{eqnarray*}
Y(\omega) \equiv y(\omega)(t)=\left((x_1 (\omega) (t),v_1(\omega) (t)), \ldots ,(x_N(\omega) (t),v_N(\omega)(t))\vphantom{l^1}\right)\, ,\quad \omega\in\Omega, 
\end{eqnarray*}
can be reconstructed by $x_i (\omega)(t) =x_i (\omega)(0)+\int_0^t v_i(\omega)(s) \, ds$, $1\le i\le N$, $t\ge 0$. 

Recall also from Definition \ref{Definition2.3} (ii) and (iii) that for a fixed $\mathbb{Y} (\omega) \in D_{\mathbb Y}$ and all $1\le i'\le N$ and $l\le I_c+I_r$ it holds that $\langle v(\tau_l (i'))/ |v(\tau_l(i'))|, n(x(\tau_l(i'))) \rangle < -\ve$, where $\ve>0$ is given in Definition \ref{Definition2.3}. In particular,  
\begin{eqnarray*}
\min_{1\le i'\le N,\, l\le I_c+I_r} M\left(x_{i'} (\tau_l(i')(\omega))\, ,\, v_{i'} (\tau_l(i') (\omega))\vphantom{l^1}\right) >0\, .
\end{eqnarray*}
Consequently, on any subset of $\t D_{\mathbb{Y}}$ for which the corresponding subset 
\begin{eqnarray*}
\mathcal{M}\subseteq\bigcup_{I \in\mathcal{I},\, m\in \{1,\ldots, I_c+I_r\}}W_{\mathbb{Y},I,m} 
\end{eqnarray*}
is compact, by the respective continuities we have 
\begin{eqnarray*}
\inf_\mathcal{M}\min_{1\le i'\le N,\, l\le I_c+I_r}M\left(x_{i'} (\tau_l(i'))\, ,\, v_{i'} (\tau_l(i'))\vphantom{l^1}\right)>0 \, .
\end{eqnarray*}
From Proposition \ref{Proposition3.12} recall the term $z^{\sh}(\mathbb{Y} (\omega ))$. Similarly to the above, by the assumptions on positivity of the initial density $p_0$ and the collision kernel $B$ made in the beginning of the present section, and the positivity of the density $g$ at the outcomes of the random variables $\gamma_k (i,j)$ in Definition \ref{Definition2.2} (iii) and (iv), we may conclude 
\begin{eqnarray*}
\sup_\mathcal{M} z^{\sh}(\mathbb{Y}(\omega))<\infty \, ,\quad \h\in{\bf T}_{\mathbb{Y}}, 
\end{eqnarray*}
on any subset of $\t D_{\mathbb{Y}}$ for which the corresponding subset $\mathcal {M}$ of $\bigcup_{I \in\mathcal{I}\, m\in \{1,\ldots, I_c+I_r\}}$ $W_{\mathbb{Y}, I,m}$ is compact.
\medskip

\nid 
{\it Step 2 } Recalling Definition \ref{Definition3.2} introduce 
\begin{eqnarray*}
\alpha_r(\sigma_k(i,j)):=\left(\left((a_{i',j';k'})_{k'\in {\mathbb N}}\right)_{\{1 \le i'<j'\le N\}}\, ;\, \left((b_{i';l'})_{l'\in {\mathbb N}}\right)_{\{1\le i'\le N\}} \right)\in\mathcal{A} 
\end{eqnarray*}
where all $b_{i';l'}=0$ and $a_{i',j';k'}=0$ for all $(i',j',k')\neq (i,j,k)$, but $a_{i,j;k}$ is the unit vector of the $r$th coordinate in ${\mathbb R}^d$. Similarly, let  
\begin{eqnarray*}
\alpha_r(\tau_l(i)):=\left(\left((a_{i',j';k'})_{k'\in {\mathbb N}}\right)_{\{1\le 
i'<j'\le N\}}\, ;\, \left((b_{i';l'})_{l'\in {\mathbb N}}\right)_{\{1\le i'\le N\}} 
\right)\in\mathcal{A} 
\end{eqnarray*}
where all $a_{i',j';k'}=0$ and $b_{i';l'}=0$ for all $(i',l')\neq (i,l)$, but $b_{i;l}$ is the unit vector of the $r$th coordinate in ${\mathbb R}^d$. Keeping Definitions \ref{Definition3.2}, \ref{Definition3.5}, and \ref{Definition3.8} in mind let us introduce  
\begin{eqnarray*} 
&&\hspace{-.5cm}\h^{(m)}_{(i,r)}\equiv\h^{(m)}_{(i,r)}(\mathbb{Y})\vphantom{\int_0}\nonumber \\ 
&&\hspace{.5cm}:=\left\{
\begin{array}{l}
\left(0;0;h(\alpha_r(\sigma_k(i,j)),{\mathbb Y})\vphantom{\displaystyle l^1}\right) \vphantom{\displaystyle\int} \\  \hphantom{yy} {\rm if}\ \rho_{m} = \sigma_k (i,j)\ \mbox{\rm for some}\ 1\le i<j\le N\ \mbox{\rm and}\ k\in {\mathbb Z}_+\vphantom{\displaystyle\int} \\ 
\left(0;0;h(\alpha_r(\tau_l(i)),{\mathbb Y})\vphantom{\displaystyle l^1}\right) \vphantom{\displaystyle\int} \\   
\hphantom{yy} {\rm if}\ \rho_{m} = \tau_l(i)\ \mbox{\rm for some}\ 1\le i\le N\ \mbox{\rm and}\ l\in {\mathbb Z}_+\vphantom{\displaystyle\int} 
\end{array}
\right.\, ,  
\end{eqnarray*}
$r\in\{1,\ldots ,d\}$, where the first zero in $(0;0;\cdot)$ is the zero in ${\mathbb R}^{2N\cdot d}$ and the second zero in $(0;0;\cdot)$ is the zero in $\mathcal{G}$. Recalling Definition \ref{Definition3.5} (a), Theorem \ref{Theorem3.6}, and Propositions \ref{Proposition3.10} as well as \ref{Proposition3.13} we observe that for $\rho_m=\tau_l(i)$ we have  $v_{ (i,r)} (\rho_m) \in D(\Nab)$ and 
\begin{eqnarray*}
\Nab v_{ (i,r)}(\rho_m) =\h^{(m)}_{(i,r)}\in {\bf T}_{\mathbb Y}\, .    
\end{eqnarray*}
This is also obvious for $\rho_m=\sigma_k (i,j)$. Therefore 
\begin{eqnarray}\label{4.4}
&&\hspace{-.5cm}\left\langle\Nab v_{(i,r)}(\rho_m),\Nab\phi(v(\rho_m)) \right \rangle_{T_{\mathbb Y}} =\partial_{(i,r)}\phi(v(\rho_m))\, ,
\end{eqnarray}
$(i,r)\in\{1,\ldots ,N\}\times\{1,\ldots ,d\}$, see Proposition \ref{Proposition3.19} (a). 

For the next calculations keep in mind that $G\in\t C_0^\infty(\t D_{\mathbb{Y}})$ implies $(\chi_{D_{\mathbb{Y},I,m}}\cdot G)\in\t C_0^\infty(\t D_{\mathbb{Y}})$, i.e. 
\begin{eqnarray*}
(\chi_{D_{\mathbb{Y},I,m}}\cdot G)\cdot\, \Nab v_{(i,r)} (\rho_m)=(\chi_{D_{\mathbb{Y},I,m}}\cdot G)\cdot \h^{(m)}_{(i,r)} \in D_0 (\delta)\, ,  
\end{eqnarray*}
cf. also Corollary \ref{Corollary3.18}. Let us note that by property (iii) of the last Notation (2), we have  with $v_i(\rho_m)\in c\{G\}$, $i=1, \ldots N$, $m=1, \ldots ,I_c+I_r$. According to Corollary \ref{Corollary3.15} (a), relation (\ref{4.4}) yields  
\begin{eqnarray*}
&&\hspace{-.5cm}E\left[\chi_{D_{\mathbb{Y},I,m}}\cdot\partial_{(i,r)}\phi(v(t))\, G\right]=E\left[\chi_{D_{\mathbb{Y},I,m}}\cdot \partial_{(i,r)}\phi(v(\rho_m))\, G\right]\vphantom{\dot{f}}  \\ 
&&\hspace{.5cm}=E\left[(\chi_{D_{\mathbb{Y},I,m}}\cdot G) \left\langle\Nab v_{(i,r)} (\rho_m), \Nab \phi (v(\rho_m))\right \rangle_{ T_{\mathbb Y}}\right]  \\ 
&&\hspace{.5cm}=E\left[\phi(v(\rho_m))\, \delta\left((\chi_{D_{\mathbb{Y},I,m}} \cdot G)\cdot \Nab v_{(i,r)}(\rho_m)\right)\right] \vphantom{\left(\dot{f}\right)_{ T_{\mathbb Y}}} \, .
\end{eqnarray*}
Since $\chi_{D_{\mathbb{Y},I,m}} \cdot G\in \t C_0^\infty(\t D_{\mathbb{Y}})$ we have 
\begin{eqnarray}\label{4.6}
\delta\left( (\chi_{D_{\mathbb{Y},I,m}} \cdot G)\cdot \Nab v_{(i,r)} (\rho_m) \right)=\delta \left((\chi_{D_{\mathbb{Y},I,m}}\cdot G)\cdot \h^{(m)}_{(i,r)} \right)=0\quad\mbox{\rm on}\quad \complement D_{\mathbb{Y},I,m}
\end{eqnarray}
by Proposition \ref{Proposition3.19} (b) together with Proposition \ref{Proposition3.10} (a) and (b) adjusted to $G\in \t C_0^\infty (\t D_\mathbb{Y})$, and therefore 
\begin{eqnarray}\label{4.5}
&&\hspace{-.5cm}E\left[\chi_{D_{\mathbb{Y},I,m}}\cdot\partial_{(i,r)}\phi(v(t))\, G\right]=E\left[\left(\chi_{D_{\mathbb{Y},I,m}}\cdot \phi(v(\rho_m))\right)\, \delta \left((\chi_{D_{\mathbb{Y},I,m}}\cdot G)\cdot\h^{(m)}_{(i,r)} \right) \right] \vphantom{\left(\dot{f}\right)_{ T_{\mathbb Y}}} \nonumber \\ 
&&\hspace{.5cm}=E\left[\left(\chi_{D_{\mathbb{Y},I,m}}\cdot \phi(v(t))\right)\, \delta \left( (\chi_{D_{\mathbb{Y},I,m}}\cdot G)\cdot \h^{(m)}_{(i,r)}\right)\right]\, , \vphantom{\left(\dot{f}\right)_{ T_{\mathbb Y}}} \nonumber \\ 
&&\hspace{.5cm}=E\left[\phi(v(t))\, \delta \left( (\chi_{D_{\mathbb{Y},I,m}}\cdot G)\cdot \h^{(m)}_{(i,r)}\right)\right]\, , \vphantom{\left(\dot{f}\right)_{ T_{\mathbb Y}}}
\end{eqnarray}
$(i,r)\in\{1,\ldots ,N\}\times\{1,\ldots ,d\}$. For such $(i,r)$ it follows from Lemma \ref{Lemma4.5} and Proposition \ref{Proposition3.10} together with Theorem \ref{Theorem3.6} that $\chi_{D_{\mathbb {Y},I,m}}\cdot v_{(i,r)}(t)=\chi_{ D_{\mathbb {Y},I,m}}\cdot v_{(i,r)}(\rho_m)\in D(\Nab)$ and 
\begin{eqnarray}\label{4.8}
&&\hspace{-.5cm}\h^{(m)}_{(i,r)}=\Nab v_{(i,r)}(\rho_m)=\Nab v_{(i,r)}(t)\quad \mbox{\rm on}\quad D_{\mathbb{Y},I,m}\, ,
\end{eqnarray}
$I\in\mathcal{I}$, $m\in\{0,1,\ldots ,I_c+I_r\}$. For such $(i,r)$ introduce also 
\begin{eqnarray*}
&&\hspace{-.5cm}H_{(i,r)}(v(t),G):=\sum_{I\in\mathcal{I},\, m\in\{0,1,\ldots ,I_c+I_r\}}\delta \left((\chi_{D_{\mathbb{Y},I,m}}\cdot G)\cdot\Nab v_{(i,r)}(t)\right) \, .
\end{eqnarray*}

Now note that Proposition \ref{Proposition3.12} is still in force for $G\in\t C_0^\infty (\t D_\mathbb{Y})$. From (\ref{4.6})-(\ref{4.8}) and Proposition \ref{Proposition3.19} (b) together with the smoothness assumptions from the beginning of Section \ref{sec:4} we conclude 
\begin{eqnarray}\label{4.6**}
\delta \left((\chi_{D_{\mathbb{Y},I,m}}\cdot G)\cdot\Nab v_{(i,r)}(t)\right)\in \t C_0^\infty(\t D_{ \mathbb{Y}}) \, , \quad I\in\mathcal{I},\, m\in\{0,1,\ldots ,I_c+I_r\},
\end{eqnarray}
recall in particular that by Step 1 the function $z^{\sh^{(m)}_{(i,r)}}$ is bounded on any subset of $D_{ \mathbb{Y},I,m}$ for which the corresponding subset of $W_{\mathbb{Y},I,m}$ is compact. Thus, 
\begin{eqnarray*}
H_{(i,r)}(v(t),G)\in \t C_0^\infty(\t D_{ \mathbb{Y}})\, . 
\end{eqnarray*}
Furthermore, we find the representation (\ref{4.3}) of $H_{(i,r)}(v(t), G)$ as stated in the proposition, and 
\begin{eqnarray*}
E\left[\partial_{(i,r)}\phi(v(t))\, G\right]=E\left[\phi(v(t))\, H_{(i,r)}(v(t),G)\right]\, .\vphantom{\int}
\end{eqnarray*}

\nid
{\it Step 3 } Focus now on (b). Keeping Lemma \ref{Lemma4.5} and the arguments of Step 2 in mind, it follows from Definition \ref{Definition3.5} (a), Theorem \ref{Theorem3.6}, Proposition \ref{Proposition3.10}, and Proposition \ref{Proposition3.13} together with Definition \ref{Definition3.14} (a) that for all $(i',r')\in \{1,\ldots ,N\}\times \{1,\ldots ,d\}$ we have $x_{(i',r')} (\rho_m)\in D(\Nab)$. Moreover, we deduce from Proposition \ref{Proposition3.19} (a) that on $D_{\mathbb {Y},I,m}$
\begin{eqnarray*}
&&\hspace{-.5cm}\left\langle\Nab v_{(i,r)}(\rho_m),\Nab \phi(x(t)) \right \rangle_{T_{\mathbb Y}}\vphantom{\int} \\ 
&&\hspace{.5cm}=\left\langle\h^{(m)}_{(i,r)},\textstyle{\sum_{(i',r')}} \partial_{(i',r')}\phi(x(t)) \cdot\Nab \textstyle{\left(x_{(i',r')}(\rho_m)+(t-\rho_m) v_{ (i',r')}(\rho_m)\right)}\right \rangle_{T_{\mathbb Y}}\vphantom{\int} \\ 
&&\hspace{.5cm}=\left\langle\h^{(m)}_{(i,r)},\textstyle{\sum_{(i',r')}} \partial_{(i',r')}\phi(x(t))\cdot (t-\rho_m) \h^{(m)}_{(i',r')}\right \rangle_{T_{\mathbb Y}} \vphantom{\int} \\
&&\hspace{.5cm}=\partial_{(i,r)}\phi(x(t))\cdot (t-\rho_m(\omega))\, ,\vphantom{\int}
\end{eqnarray*}
$(i,r)\in\{1,\ldots ,N\}\times\{1,\ldots ,d\}$. Here we mention that, in order to proceed from the second to the third line, we have used  $\langle \h^{(m)}_{(i,r)}, \Nab x_{(i',r')}(\rho_m)\rangle_{T_{\mathbb Y}}=0$ since $x_{(i',r')}(\rho_m)$ is on $D_{\mathbb{Y},I,m}$ a smooth function of $x(0)$ and $v(\rho(0)),\ldots ,v(\rho (m-1))$. As in Step 2 arrive at 
\begin{eqnarray*}
&&\hspace{-.5cm}E\left[\chi_{D_{\mathbb{Y},I,m}}\cdot\partial_{(i,r)}\phi(x(t))\, G\right]\vphantom{\int}\nonumber \\ 
&&\hspace{.5cm}=E\left[\frac{\chi_{D_{\mathbb{Y},I,m}}\cdot G}{t-\rho_m} \left \langle\Nab v_{(i,r)} (\rho_m), \Nab \phi (x(t))\right \rangle_{ T_{\mathbb Y}}\right]\nonumber \\ 
&&\hspace{.5cm}=E\left[\phi(x(t))\, \delta\left(\frac{\chi_{D_{\mathbb{Y},I,m}}\cdot G}{t-\rho_m}\cdot \Nab v_{(i,r)} (\rho_m)\right)\right] \vphantom{\left(\dot{f}\right)_{ T_{\mathbb Y}}}\nonumber \\ 
&&\hspace{.5cm}=E\left[\phi(x(t))\, \delta\left(\frac{\chi_{D_{\mathbb{Y},I,m}}\cdot G}{t-\rho_m}\cdot \Nab v_{(i,r)} (t)\right)\right] \vphantom{\left(\dot{f}\right)_{ T_{\mathbb Y}}} \, .
\end{eqnarray*}
Note here that $t-\rho_m\ge c$ for some $c>0$ on every subset of $D_{\mathbb{Y}, I,m}$ for which the corresponding subset $W_{\mathbb{Y},I,m}$ is compact. Furthermore note that on such a subset of $D_{\mathbb{Y}, I,m}$ the random time $\rho_m$ is a smooth function of $x(0)$ and $v(\rho(0)),\ldots ,v(\rho(m-1))$. Therefore 
\begin{eqnarray*}
\frac{\chi_{D_{\mathbb{Y},I,m}}\cdot G}{t-\rho_m}\in \t C_0^\infty (\t D_{\mathbb{Y}})
\end{eqnarray*}
which together with (\ref{4.6**}) implies 
\begin{eqnarray*}
&&\hspace{-.5cm}H_{(i,r)}(x(t),G):=\sum_{I\in\mathcal{I},\, m\in \{0,1,\ldots ,I_c+I_r\}}\delta \left(\frac{\chi_{D_{\mathbb{Y},I,m}}\cdot G}{t-\rho_m}\cdot \Nab v_{(i,r)}(t)\right)\in \t C_0^\infty (\t D_{\mathbb{Y}})\, .
\end{eqnarray*}
\qed
\bigskip 

%\nid
{\sc Notation. (1) } For any multi-index $\alpha=(\alpha_{(1,1)},\ldots ,\alpha_{(N,d)}) \in {\mathbb Z}_+^{N\cdot d}$ let 
\begin{eqnarray*}
|\alpha|:=\sum_{i=1}^N\sum_{r=1}^d\alpha_{(i,r)}\quad\mbox{\rm and}\quad \partial^\alpha :=\left(\partial_{(1,1)}^{\alpha_{(1,1)}}\ldots\, \partial_{(N,d)}^{\alpha_{(N,d)}} \right)\, . 
\end{eqnarray*}
{\sc (2) } Let $e_{(i,r)}=(0\ldots 0\, 1\, 0\ldots 0)$ denote the multi-index corresponding to the coordinate $(i,r)\in\{1,\ldots ,N\}\times\{1,\ldots ,d\}$. 

\bigskip

Proposition \ref{Proposition4.2} immediately provides the following. 
\begin{corollary}\label{Corollary4.3} 
Let $G\in\t C_0^\infty(\t D_{\mathbb Y})$ and let $\phi\in C_b^\infty({\mathbb R}^{N\cdot d})$. Then for any multi-index $\alpha\in {\mathbb Z}_+^{N\cdot d}$ with $|\alpha|\ge 1$ we have 
\begin{eqnarray*}
E\left[\partial^\alpha\phi(v(t))\, G\right]=E\left[\phi(v(t))\, H_\alpha (v(t) 
,G)\right]
\end{eqnarray*}
as well as 
\begin{eqnarray*}
E\left[\partial^\alpha\phi(x(t))\, G\right]=E\left[\phi(x(t))\, H_\alpha (x(t) 
,G)\right]
\end{eqnarray*}
where, for $\beta=\alpha+e_{(i,r)}$, $(i,r)\in\{1,\ldots ,N\}\times\{1,\ldots ,d\}$, 
\begin{eqnarray*}
H_\beta(v(t),G)=H_{(i,r)}\left(v(t),H_\alpha(v(t),G)\right)\in\t C_0^\infty 
(\t D_{\mathbb Y})  
\end{eqnarray*}
as well as 
\begin{eqnarray*}
H_\beta(x(t),G)=H_{(i,r)}\left(x(t),H_\alpha(x(t),G)\right)\in\t C_0^\infty 
(\t D_{\mathbb Y})\, . 
\end{eqnarray*}
\end{corollary}
\bigskip

The following is the main result of the paper.
\begin{theorem}\label{Theorem4.5} 
Let $\mathbb{Y}\equiv \left(\mathbb{Y}(0); (v_1(s),\ldots ,v_N(s))_{s\ge 0} \vphantom {l^1}\right)$ be a reduced process with reduced trajectories $\omega \equiv \mathbb{Y}(\omega)$. Fix $t>0$. Under the hypotheses of this section we have the following. \\ 
(a) There exists a finite collection $\mathcal{C}=\{\mathcal{C}_i\}$ of disjoint open subsets of $V^N$ satisfying 
\begin{eqnarray}\label{4.6*}
\bigcup_{\mathcal{C}_i\in \mathcal{C}}\mathcal{C}_i\subseteq\{v(\omega) (t): \omega \in D_{\mathbb{Y}}\} \subseteq\bigcup_{\mathcal{C}_i\in \mathcal{C}}\overline{\mathcal {C}_i} 
\end{eqnarray}
such that the random variable $v(t)\equiv v(\omega)(t)$ possesses on every $\mathcal{C}_i\in \mathcal{C}$ a not necessarily bounded infinitely differentiable density $p_{v,i}$ with respect to the Lebesgue measure. \\ 
(b) There exists a finite collection $\mathcal{D}=\{\mathcal{D}_i\}$ of disjoint open subsets of $D^N$ satisfying 
\begin{eqnarray*}
\bigcup_{\mathcal{D}_i\in \mathcal{D}}\mathcal{D}_i\subseteq\{x(\omega) (t): \omega \in D_{\mathbb{Y}}\} \subseteq\bigcup_{\mathcal{D}_i\in \mathcal{D}}\overline{\mathcal {D}_i} 
\end{eqnarray*}
such that the random variable $x(t)\equiv x(\omega)(t)$ possesses on every $\mathcal {D}_i\in \mathcal{D}$ a not necessarily bounded infinitely differentiable density $p_{x,i}$ with respect to the Lebesgue measure. 
\end{theorem}
Proof. We prove (a). Likewise (b) can be proved.
\medskip

\nid
{\it Step 1 } As an immediate consequence of Lemmas \ref{Lemma4.4} and \ref{Lemma4.5}, for all $I\in \mathcal {I}$ and $m\in\{0,1,\ldots ,I_c+I_r\}$ the set $\{v(\omega)(t): \omega \equiv \mathbb {Y}(\omega)\in D_{\mathbb{Y},I,m}\} =\{v (\rho_m (\mathbb{Y} (\omega))): \mathbb {Y}(\omega)\in D_{\mathbb{Y},I,m}\}$ is an open subset of $V^N$. Recalling that $\t D_{\mathbb{Y}}=\bigcup_{I\in\mathcal{I},\, m\in\{0,1, \ldots ,I_c+I_r\}}D_{\mathbb{Y},I,m}$, it turns out that $\{v(\omega) (t): \omega \equiv \mathbb {Y}(\omega)\in \t D_{\mathbb{Y}}\setminus D_{\mathbb{Y}, I,m}\} $ is also an open subset of $V^N$. Consider all multiple intersections of the the sets 
\begin{eqnarray*}
\{v(\omega) (t): \omega \in D_{\mathbb{Y},I,m}\} \quad \mbox{\rm and} \quad \{v(\omega)(t):\omega \equiv \mathbb {Y}(\omega)\in \t D_{\mathbb{Y}}\setminus D_{\mathbb{Y},I,m}\} 
\end{eqnarray*}
for all $I\in \mathcal {I}$ and $m\in\{0,1,\ldots ,I_c+I_r\}$ and 
%
% and we have 
%
% \begin{eqnarray*}
% \{v(\omega) (t): \omega \in D_{\mathbb{Y}}\}\subset\bigcup_{I\in\mathcal{I},\, m\in\{0,1,\ldots , I_c+I_r\}}\overline{\{v(\omega) (t): \omega \in D_{\mathbb{Y},I,m}\}}\, .
% \end{eqnarray*}
%
% Consider the sets  
%
% \begin{eqnarray*}
% &&\hspace{-.5cm}\{v(\omega) (t): \omega \in D_{\mathbb{Y},I_1,m_1}\} \vphantom{\dot{f}} \\ 
% &&\hspace{-.5cm}\{v(\omega) (t): \omega \in D_{\mathbb{Y},I_1,m_1}\} \cap \{v (\omega) (t): \omega \in D_{\mathbb{Y},I_2,m_2}\}  \\ 
% &&\hspace{-.5cm}\{v(\omega) (t): \omega \in D_{\mathbb{Y},I_1,m_1}\} \cap \ldots \cap \{v (\omega) (t): \omega \in D_{\mathbb{Y},I_3,m_3}\}  \\ 
% &&\hspace{-.5cm}\{v(\omega) (t): \omega \in D_{\mathbb{Y},I_1,m_1}\} \cap \ldots \cap \{v (\omega) (t): \omega \in D_{\mathbb{Y},I_4,m_4}\}  \\ 
% &&\hspace{.5cm} \vphantom{\dot{f}} \cdots \, ,\\
% &&\hspace{.5cm} \vphantom{\dot{f}}I_1,I_2,\ldots\in\mathcal{I},\, m_1,m_2,\ldots\in\{0,1,\ldots , I_c+I_r\}\, , 
% \end{eqnarray*}
%
and let $\mathcal{C}$ be the collection of all of these intersections which are non-empty and do not contain any other of these sets, except for the empty set. Then $\mathcal{C}$ is a finite collection of disjoint open subsets of $V^N$ satisfying (\ref{4.6*}), recall again Lemmas \ref{Lemma4.4} and \ref{Lemma4.5}. 
\medskip

\nid
{\it Step 2 } The crucial observation is that $G\in\t C_0^\infty(\t D_{\mathbb Y})$ must be zero on $D_{\mathbb Y}\setminus \t D_{\mathbb Y}$
%$\{\mathbb{Y}(\omega):v(\omega)(t) \in\bigcup_{ \mathcal {C}_i \in \mathcal{C}}\overline{\mathcal {C}_i}\setminus \bigcup_{\mathcal {C}_i\in \mathcal{C}}\mathcal {C}_i \}$ 
but can take values different from zero on $\t D_{\mathbb Y}\supseteq \{\mathbb{Y} (\omega):v(\omega)(t)\in\bigcup_{\mathcal {C}_i\in \mathcal{C}} \mathcal {C}_i\}$. 

Looking at the statement of Proposition \ref{Proposition4.2} and recalling Corollary \ref{Corollary3.18} it turns out that for $\phi\in C_0^\infty(\mathbb {R}^{N \cdot d})$ supported by a compact subset of $\bigcup_{\mathcal {C}_i\in \mathcal{C}} \mathcal {C}_i$ and $G\in\t C_0^\infty(\t D_{\mathbb{Y}})$ 
\begin{eqnarray*}
&&\hspace{-.0cm}E\left[\partial_{(i,r)}\phi(v(t))\, G\right]=E\left[\phi(v(t))\, \sum_{I\in\mathcal{I},\, m\in\{0,1,\ldots ,I_c+I_r\}}\delta \left((\chi_{D_{\mathbb{Y},I,m}}\cdot G)\cdot\h^{(m)}_{(i,r)}\right)\right] \\ 
&&\hspace{1.0cm}=-\sum_{I\in\mathcal{I},\, m\in\{0,1,\ldots ,I_c+I_r\}}E\left[\phi(v(t))\, \left(\partial_{\sh_{(i,r)}^m}(\chi_{D_{\mathbb{Y},I,m}}\cdot G)+\chi_{D_{\mathbb{Y},I,m}}\cdot G\cdot z^{\sh^{(m)}_{(i,r)}}\right)\right]\, .
\end{eqnarray*}
Now choose $G\in\t C_0^\infty(\t D_{\mathbb{Y}})$ such that $G=1$ on some subset $\t D_{c, \mathbb{Y}}$ of $\t D_{ \mathbb{Y}}$ for which the corresponding subset of $\bigcup_{I\in \mathcal {I}, \, m\in\{0,1,\ldots, I_c+I_r\}}W_{\mathbb{Y},I,m}$ is compact such that $\{v(\omega)(t):\omega\equiv\mathbb{Y}(\omega)\in \t D_{c, \mathbb{Y}}\}\supseteq\, $supp$\, \phi$. 

Regarding $m$ as a random variable on $\t D_{\mathbb{Y}}$, for any $\phi\in C_0^\infty(\mathbb {R}^{N \cdot d})$ supported by a compact subset of $\bigcup_{\mathcal {C}_i\in \mathcal{C}} \mathcal {C}_i$ we get
\begin{eqnarray*}
&&\hspace{-.5cm}E\left[\partial_{(i,r)}\phi(v(t))\right]=-E\left[\phi(v(t))\, z^{\sh^{(m)}_{(i,r)}}\right] = E\left[\phi(v(t))\, H_{(i,r)}(v(t),1)\right]
\end{eqnarray*}
and similarly from Corollary \ref{Corollary4.3}
\begin{eqnarray}\label{4.8*}
E\left[\partial^\alpha\phi(v(t))\right]=E\left[\phi(v(t))\, H_\alpha (v(t) 
,1)\right]
\end{eqnarray}
for any multi-index $\alpha\in {\mathbb Z}_+^{N\cdot d}$ with $|\alpha|\ge 1$. This implies part (a) of the theorem. Note that, in order to derive the claim from (\ref{4.8*}), we cannot adopt the method of  e.g. \cite{Sa05}, Chapter I. We follow and slightly  modify to $\phi\in C_0^\infty(\mathbb {R}^{N \cdot d})$ the lines of e.g. \cite{Ku13}, Proposition 2.1.12.  
\qed

\end{document}